\documentclass{article}

\usepackage[left=1.25in,top=1.25in,right=1.25in,bottom=1.25in,head=1.25in]{geometry}
\usepackage{amsfonts,amsmath,amssymb,amsthm}
\usepackage{verbatim,float,url,dsfont}
\usepackage{graphicx,subfigure,psfrag}
\usepackage{algorithm,algorithmic}
\usepackage[colon]{natbib}
\usepackage[colorlinks=true,citecolor=blue,urlcolor=blue,linkcolor=blue]{hyperref}
\usepackage{rotating}

\newtheorem{theorem}{Theorem}
\newtheorem{lemma}{Lemma}

\newtheorem{proposition}{Proposition}
\theoremstyle{definition}
\newtheorem{remark}{Remark}

\newtheorem{assumption}{Assumption}

\newcommand{\argmin}{\mathop{\mathrm{argmin}}}

\def\R{\mathbb{R}}
\def\E{\mathbb{E}}
\def\P{\mathbb{P}}
\def\Cov{\mathrm{Cov}}

\def\half{\frac{1}{2}}

\def\tr{\mathrm{tr}}
\def\df{\mathrm{df}}

\def\col{\mathrm{col}}

\def\rank{\mathrm{rank}}

\def\sign{\mathrm{sign}}

\def\diag{\mathrm{diag}}

\def\hbeta{\hat{\beta}}

\def\htheta{\hat{\theta}}
\def\halpha{\hat{\alpha}}

\def\heta{{\hat{\eta}}}
\def\hs{\hat{s}}

\def\cD{\mathcal{D}}

\def\cH{\mathcal{H}}

\def\cN{\mathcal{N}}

\def\Risk{\mathrm{Risk}}
\def\Err{\mathrm{Err}}
\def\hRisk{\widehat{\mathrm{Risk}}}
\def\hErr{\widehat{\mathrm{Err}}} 
\def\hdf{\widehat{\df}}
\def\hTheta{\widehat{\Theta}}
\def\Opt{\mathrm{Opt}}
\def\ExOpt{\mathrm{ExOpt}}
\def\edf{\mathrm{edf}}
\def\hOpt{\widehat{\mathrm{Opt}}}

\title{Excess Optimism: How Biased is the Apparent Error of an
  Estimator Tuned by SURE?}    
\author{Ryan J. Tibshirani \and Saharon Rosset}
\date{}

\begin{document}
\maketitle

\begin{abstract}
Nearly all estimators in statistical prediction come with an
associated tuning parameter, in one way or another.
Common practice, given data, is to choose the tuning parameter value
that minimizes a constructed estimate of the prediction error of the
estimator;
% Of course, estimating prediction error has a long history
% in statistics, and many methods have been proposed for this problem;   
we focus on Stein's unbiased risk estimator, or SURE
\citep{stein1981estimation,efron1986biased},
which forms an unbiased estimate of the prediction error by
augmenting the observed training error with an
estimate of the degrees of freedom of the estimator.  
Parameter tuning via SURE minimization has been 
advocated by many
authors, in a wide variety of problem settings, and in general, it is
natural to ask: what is the prediction error of the SURE-tuned
estimator?  An obvious strategy would be simply use 
the apparent error estimate as reported by SURE, 
i.e., the value of the SURE criterion at its minimum, 
to estimate the prediction error of the SURE-tuned estimator. 
But this is no longer unbiased; in fact, we would expect that the 
minimum of the SURE criterion is systematically biased downwards  
for the true prediction error. In this paper, we define the excess
optimism to be the amount of this downward bias in the SURE 
minimum.  We argue that the following two properties motivate the
study of excess optimism: (i) an unbiased estimate of excess 
optimism, added to the SURE criterion at its minimum, gives  
an unbiased estimate of the prediction error of the SURE-tuned 
estimator; (ii) excess optimism serves as an upper bound on the excess
risk, i.e., the difference between the risk of the SURE-tuned
estimator and the oracle risk (where the oracle uses the best fixed
tuning parameter choice).  We study excess optimism in two common
settings: the families of shrinkage and subset regression estimators.
Our main results include a James-Stein-like property of
SURE-tuned shrinkage estimation, which is shown to dominate the MLE,
and both upper and lower bounds on excess optimism for 
SURE-tuned subset regression; when the collection of subsets here is 
nested, our bounds are particularly tight, and reveal that in the case
of no signal, the excess optimism is always in between 0 and 10
degrees of freedom, no matter how many models are being selected 
from. 
% We describe a general approach for characterizing excess optimism 
% using (extensions of) Stein's formula.
We also describe a bootstrap method for estimating excess optimism,
and outline some extensions of our framework beyond the standard
homoskedastic, squared error model that we consider throughout
majority of the paper. 
\end{abstract}

\section{Introduction}
\label{sec:intro}

Consider data $Y \in \R^n$, drawn from a generic model
\begin{equation}
\label{eq:data_model}
Y \sim F, \quad \text{where} \;
\E(Y) = \theta_0, \; \Cov(Y) = \sigma^2 I. 
\end{equation} 
The mean $\theta_0 \in \R^n$ is unknown, and the variance 
$\sigma^2>0$ is assumed to be known. Let \smash{$\htheta \in \R^n$}
denote an estimator of the mean. Define the prediction error, also
called test error or just error for short, of \smash{$\htheta$} by     
\begin{equation}
\label{eq:err}
\Err(\htheta) = \E\|Y^*-\htheta(Y)\|_2^2,
\end{equation}
where $Y^* \sim F$ is independent of $Y$ and the 
expectation is taken over all that is random (over both $Y,Y^*$).  
A remark about notation: we write \smash{$\htheta$} to  
denote an {\it estimator} (also called a rule, procedure, or   
algorithm), and \smash{$\htheta(Y)$} to denote an 
{\it estimate} (a particular realization given data $Y$).  
Hence it is perfectly well-defined to write the error as 
\smash{$\Err(\htheta)$}; this is indeed a fixed (i.e., nonrandom) 
quantity, because \smash{$\htheta$} represents a rule, not  
a random variable. %While this may sound pedantic, it 
This will be helpful to keep in mind when our notation becomes a bit
more complicated. %shortly. 

Estimating prediction error as in \eqref{eq:err} is a 
classical problem in statistics. One convenient method that does  
not require the use of held-out data stems from the 
{\it optimism theorem}, which says that 
\begin{equation}
\label{eq:opt_thm}
\Err(\htheta) = \E\|Y-\htheta(Y)\|_2^2 +  2 \sigma^2
\df(\htheta),  
\end{equation}
where \smash{$\df(\htheta)$}, called the {\it degrees of freedom}
of \smash{$\htheta$}, is defined as
\begin{equation}
\label{eq:df}
\df(\htheta) = \frac{1}{\sigma^2} \tr\big(\Cov 
(\htheta(Y),Y)\big) = \frac{1}{\sigma^2}  
\sum_{i=1}^n \Cov(\htheta_i(Y), Y_i).
\end{equation}
% The optimism theorem \eqref{eq:opt_thm} is easy to verify: we
% begin with \smash{$\E\|Y^*-\htheta(Y)\|_2^2=n\sigma^2 +
%   \E\|\theta_0-\htheta(Y)\|_2^2$}, add and subtract $Y$  
% inside the norm in the second term, and expand. 
Let us define the {\it optimism} of \smash{$\htheta$} as 
\smash{$\Opt(\htheta) = \E\|Y^* - \htheta(Y)\|_2^2 - \E\|Y -
  \htheta(Y)\|_2^2$}, the difference in prediction and training
errors. Then, we can rewrite \eqref{eq:opt_thm} as 
\begin{equation}
\label{eq:opt_thm_2}
\Opt(\htheta) = 2 \sigma^2 \df(\htheta),
\end{equation}
which explains its name. A nice treatment of the optimism theorem can
be found in \citet{efron2004estimation}, though the idea can be found
much earlier, e.g., 
\citet{mallows1973comments,stein1981estimation,efron1986biased}. 
In fact, \citet{efron2004estimation} developed more general versions
of the optimism theorem in \eqref{eq:opt_thm}, beyond
the standard setup in \eqref{eq:data_model}, \eqref{eq:err}; we discuss
extensions along these lines in Section \ref{sec:efron_q}.

The optimism theorem in \eqref{eq:opt_thm} suggests an
estimator for the error in \eqref{eq:err}, defined by
\begin{equation}
\label{eq:sure}
\hErr(Y) = 
\|Y-\htheta(Y)\|_2^2 + 2\sigma^2 \hdf(Y), 
\end{equation}
where \smash{$\hdf$} is any unbiased estimator
of the degrees of freedom of $\htheta$, as defined in \eqref{eq:df}, 
i.e., it satisfies \smash{$\E[\hdf(Y)] = \df(\htheta)$}.  
Clearly, from \eqref{eq:sure} and \eqref{eq:opt_thm}, we see that
\begin{equation}
\label{eq:sure_unbiased} 
\E[\hErr(Y)]= \Err(\htheta),
\end{equation}
i.e., \smash{$\hErr$} is an unbiased estimator of
the prediction error of \smash{$\htheta$}. We will call the estimator  
\smash{$\hErr$} in \eqref{eq:sure} {\it Stein's unbiased 
  risk estimator}, or SURE, in honor of
\citet{stein1981estimation}. This is somewhat of an abuse of
notation, as \smash{$\hErr$} is actually an estimate of prediction
error, \smash{$\Err(\htheta)$} in \eqref{eq:err}, and not risk,  
\begin{equation}
\label{eq:risk}
\Risk(\htheta) = \E\|\theta_0-\htheta(Y)\|_2^2.
\end{equation}
However, the two are essentially equivalent 
notions, because \smash{$\Err(\htheta) = n\sigma^2 + 
  \Risk(\htheta)$}, noted above.  Also, the term
SURE is already in wide use in the literature, so we stick with it
here.  

We note that, when \smash{$\htheta$} is a linear regression 
estimator (onto a fixed and full column rank design matrix), the
degrees of freedom of \smash{$\htheta$} is simply $p$, the number  
of predictor variables in the regression, and SURE reduces to Mallows'
well-known $C_p$ formula \citep{mallows1973comments}.

\subsection{Stein's formula}

\citet{stein1981estimation} studied a risk decomposition, as 
in \eqref{eq:sure}, with the specific degrees of freedom estimator 
\begin{equation}
\label{eq:stein_div}
\hdf(Y) = (\nabla \cdot \htheta)(Y) = 
\sum_{i=1}^n \frac{\partial \htheta_i}{\partial Y_i} (Y),
\end{equation}
called the divergence of the map 
\smash{$\htheta : \R^n \to \R^n$}. Assuming a normal distribution 
$F=N(\theta_0,\sigma^2 I)$ for 
the data in \eqref{eq:data_model} and regularity conditions on
\smash{$\htheta$} (i.e., continuity, weak differentiability, and
essential boundedness of the weak derivative), Stein showed that the
divergence estimator defined by \eqref{eq:stein_div} is unbiased for 
\smash{$\df(\htheta)$}; to be explicit
\begin{equation}
\label{eq:stein_formula}
\df(\htheta) = \E\bigg[
\sum_{i=1}^n \frac{\partial \htheta_i}{\partial Y_i} (Y) \bigg].
\end{equation}
This elegant and important result has had significant a following in
statistics (e.g., see the 
references given in the next subsection). 

% Throughout, we will use the term SURE to refer to the estimator
% \smash{$\hErr$} in \eqref{eq:sure} for {\it any} unbiased estimator
% of degrees of freedom, not just that coming from Stein's divergence 
% formula in \eqref{eq:stein_div}.

\subsection{Parameter tuning via SURE}

Here and henceforth, we write \smash{$\htheta_s$} for the
estimator of interest, where the subscript $s$ highlights the 
dependence of this estimator on a tuning parameter,
taking values in a set $S$. The term ``tuning parameter''
is used loosely, and we do not place any restrictions on
$S$ (e.g., this can be a continuous or a discrete collection of 
tuning parameter values). Abstractly, we can just think of 
\smash{$\{\htheta_s : s \in S\}$} as a family of estimators under 
consideration.  
We use \smash{$\hErr_s$} to denote the prediction error estimator in
\eqref{eq:sure} for 
\smash{$\htheta_s$}, and \smash{$\hdf_s$} to denote an unbiased
degrees of freedom estimator for \smash{$\htheta_s$}. 

One sensible strategy for choosing the tuning parameter $s$,
associated with our estimator \smash{$\htheta_s$}, is to select the    
value minimizing SURE in \eqref{eq:sure}, denoted   
\begin{equation}
\label{eq:sure_tuning}
\hs(Y) = \argmin_{s \in S} \,
\hErr_s(Y).
\end{equation}
We can think of \smash{$\hs$} as an estimator of some optimal
tuning parameter value, namely, an estimator of 
\begin{equation}
\label{eq:oracle_tuning}
s_0 = \argmin_{s \in S} \,
\Err(\htheta_s),
\end{equation}
the tuning parameter value minimizing error.  
When \smash{$\htheta_s$} is the linear regression estimator onto a set   
of predictor variables indexed by the parameter $s$, the rule in 
\eqref{eq:sure_tuning} 
encompasses model selection via $C_p$ minimization, which is 
a classical topic in statistics. In general, tuning parameter
selection via SURE minimization has been widely advocated by authors
across various problem settings, e.g.,     
\citet{donoho1995adapting,johnstone1999wavelet,
  zou2007degrees,zou2008regularized,
  tibshirani2011solution,tibshirani2012degrees,
  candes2013unbiased,ulfarsson2013tuning1,ulfarsson2013tuning2,
  chen2015degrees}, 
just to name a few.
% TODO: references from all the French papers on df?

\subsection{What is the error of the SURE-tuned estimator?}
\label{sec:err_sure}

Having decided to use \smash{$\hs$} as a rule for choosing the
tuning parameter, it is natural to ask: what is 
the error of the subsequent SURE-tuned estimator
\smash{$\htheta_{\hs}$}?  To  be explicit, this estimator produces the
estimate \smash{$\htheta_{\hs(Y)}(Y)$} given data $Y$, where 
\smash{$\hs(Y)$} is the tuning parameter value minimizing the SURE 
criterion, as in \eqref{eq:sure_tuning}. Initially, it might seem 
reasonable to use the apparent error estimate given to us SURE,  
i.e., \smash{$\hErr_{\hs(Y)}(Y)$}, to estimate the prediction error of
\smash{$\htheta_{\hs}$}.  To be explicit, this gives
\begin{equation*}
\hErr_{\hs(Y)}(Y) = \|Y-\htheta_{\hs(Y)}(Y)\|_2^2 + 
2\sigma^2 \hdf_{\hs(Y)}(Y) 
\end{equation*}
at each given data realization $Y$.  However, even
though \smash{$\hErr_s$} is unbiased for \smash{$\Err(\htheta_s)$} for
each fixed $s \in S$, the estimator \smash{$\hErr_{\hs}$} is no longer  
generally unbiased for \smash{$\Err(\htheta_{\hs})$}, and commonly, 
it will be too optimistic, i.e., we will commonly observe that
\begin{equation}
\label{eq:sure_bias}
\E[\hErr_{\hs(Y)}(Y)] < \Err(\htheta_{\hs}) = \E\|Y^* -  
\htheta_{\hs(Y)}(Y)\|_2^2.
\end{equation}
After all, for each data instance $Y$, the value \smash{$\hs(Y)$} is
specifically chosen to minimize \smash{$\hErr_s(Y)$} over all $s \in
S$, and thus we would expect \smash{$\hErr_{\hs}$} to be biased  
downwards as an estimator of the error of \smash{$\htheta_{\hs}$}.  
Of course, the optimism of training error, as displayed in
\eqref{eq:opt_thm}, \eqref{eq:df}, \eqref{eq:opt_thm_2}, is by now a
central principle in statistics and (we believe) nearly all
statisticians are aware of and account for this optimism in applied
statistical modeling. But the optimism of the optimized SURE criterion 
itself, as suggested in \eqref{eq:sure_bias}, is more subtle and has
received less attention.  

\subsection{Excess optimism}
\label{sec:ex_opt}

In light of the above discussion, we define the {\it excess optimism}
associated with \smash{$\htheta_{\hs}$} by\footnote{The excess
  optimism here is not only associated with \smash{$\htheta_{\hs}$}
  itself, but also with the the SURE family \smash{$\{\hErr_s : s \in
    S\}$}, used to define \smash{$\hs$}. This is meant to be implicit
  in our language and our notation.} 
\begin{equation}
\label{eq:ex_opt}
\ExOpt(\htheta_{\hs}) = \Err(\htheta_{\hs}) - \E[\hErr_{\hs(Y)}(Y)].
\end{equation}
We similarly define the {\it excess degrees of freedom} of 
\smash{$\htheta_{\hs}$} by
\begin{equation}
\label{eq:ex_df}
\edf(\htheta_{\hs}) = \df(\htheta_{\hs}) - 
\E[\hdf_{\hs(Y)}(Y)].
\end{equation}
The same motivation for excess optimism can be retold from
the perspective of degrees of freedom: even though the degrees of 
freedom estimator \smash{$\hdf_s$} is unbiased for
\smash{$\df(\htheta_s)$} for each fixed $s \in S$, we should not
expect \smash{$\hdf_{\hs}$} to be unbiased for
\smash{$\df(\htheta_{\hs})$}, and again it will commonly biased
downwards, i.e., excess degrees of freedom in \eqref{eq:ex_df}
will be commonly positive. 

It should be noted that the two perspectives---excess optimism and
excess degrees of freedom---are equivalent, as the optimism theorem in
\eqref{eq:opt_thm} (which holds for any estimator) applied to
\smash{$\htheta_{\hs}$} tells us that  
\begin{equation*}
\Err(\htheta_{\hs}) = \E\|Y-\htheta_{\hs(Y)}(Y)\|_2^2  
+  2 \sigma^2 \df(\htheta_{\hs}).
\end{equation*}
Therefore, we have
\begin{equation*}
\ExOpt(\htheta_{\hs}) = 2\sigma^2 \edf(\htheta_{\hs}),
\end{equation*}
analogous to the usual relationship between optimism and degrees of
freedom.  

It should also be noted that the focus on prediction error, rather
than risk, is a decision based on ease of exposition, and that excess 
optimism can be equivalently expressed in terms of risk, i.e., 
\begin{equation}
\label{eq:ex_opt_risk}
\ExOpt(\htheta_{\hs}) = \Risk(\htheta_{\hs}) - \E[\hRisk_{\hs(Y)}(Y)],  
\end{equation}
where we define \smash{$\hRisk_s=\hErr_s - n\sigma^2$}, an unbiased     
estimator of \smash{$\Risk(\htheta_s)$} in \eqref{eq:risk}, for each
$s \in S$. 

Finally, a somewhat obvious but important point is the following:
an unbiased estimator \smash{$\widehat{\edf}$} of excess degrees
of freedom \smash{$\edf(\htheta_{\hs})$} leads to an unbiased
estimator of prediction error \smash{$\Err(\htheta_{\hs})$},
i.e., \smash{$\hErr_{\hs}+2\sigma^2\widehat{\edf}$}, by
construction of excess degrees of freedom in 
\eqref{eq:ex_df}.  Likewise,
\smash{$\hRisk_{\hs}+2\sigma^2\widehat{\edf}$} is an unbiased
estimator of the risk \smash{$\Risk(\htheta_{\hs})$}.  

\subsection{Is excess optimism always nonnegative?}
\label{sec:ex_opt_nneg}

Intuitively, it seems reasonable to believe that excess optimism
should be always nonnegative, i.e., in any setting, the expectation
of the SURE criterion at its minimum should be no more than the
actual error rate of the SURE-tuned estimator.  However, we are not
able to give a general proof of this phenomenon. 

In each setting that we study in this work---shrinkage estimators,
subset regression estimators, and soft-thresholding estimators---we
prove that the excess degrees of freedom is nonnegative,
abeit with different proof techniques.  We have not seen evidence,
theoretical or empirical, to suggest that excess degrees of freedom
can be negative for certain classes 
of estimators; but of course, without a
general proof of nonnegativity, we cannot rule out the
possibility that it is negative in some (likely pathological)
situations. 

\subsection{Summary of contributions}
\label{sec:summary}

The goal of this work is to understand excess optimism, or
equivalently, excess degrees of freedom, associated with estimators
that are tuned by optimizing SURE. Below, we provide a outline of our 
results and contributions.

\begin{itemize}
\item In Section \ref{sec:oracle}, we develop further motivation for
  the study of excess optimism, by showing that it upper bounds the 
  excess risk, i.e., the difference between the risk of the estimator
  in question and the oracle risk, in Theorem \ref{thm:oracle_bd}.

\item In Section \ref{sec:shrinkage}, we precisely characterize (and
  give an unbiased estimator for) the
  excess degrees of freedom of the SURE-tuned shrinkage estimator,
  both in a classical normal means problem setting and in a regression 
  setting, in \eqref{eq:ex_df_shrink} and \eqref{eq:ex_df_shrink_reg}, 
  respectively.  This shows that the excess degrees of freedom in both
  of these settings always nonnegative, and at most 2.  Our analysis
  also reveals an   
  interesting connection between SURE-tuned shrinkage  
  estimation and James-Stein estimation.  

\item In Sections \ref{sec:subset_reg} and \ref{sec:subset_reg_mh}, we
  derive bounds on the excess degrees of freedom of the SURE-tuned
  subset regression estimator (or equivalently, the $C_p$-tuned subset 
  regression estimator), using different approaches. Theorem
  \ref{thm:ex_df_subset_reg} shows from first principles that, under  
  reasonable conditions on the subset regression models being  
  considered, the excess degrees of freedom of SURE-tuned  
  subset regression is small compared to the oracle risk.
  Theorems \ref{thm:ex_df_subset_reg_mh} and  
  \ref{thm:ex_df_subset_reg_nested} are derived using a 
  more refined general result, from \citet{mikkelsen2016degrees},
  and present exact (though not always explicitly computable)
  expressions for excess degrees of freedom.  Some implications for
  excess degrees of freedom in SURE-tuned subset regression
  estimator: we see that it is always nonnegative, and 
  is (perhaps) surprisingly small for nested collections of 
  subsets, e.g., it is at most 10 for any nested collection
  (no matter the number of predictors) when $\theta_0=0$.  

\item In Section \ref{sec:stein}, we consider strategies for 
  characterizing the excess degrees of freedom of generic estimators
  using Stein's formula, and extensions of Stein's formula for
  discontinuous mappings from 
  \citet{tibshirani2015degrees,mikkelsen2016degrees}.
  We use the extension from \citet{tibshirani2015degrees} in Section 
  \ref{sec:soft_thresh} to prove that excess degrees of freedom in
  SURE-tuned soft-thresholding is always nonnegative. We use 
  that from \citet{mikkelsen2016degrees} in Section
  \ref{sec:subset_reg_mh} to prove results on subset regression,
  already described.

\item In Section \ref{sec:bootstrap}, we study a simple bootstrap
  procedure for estimating excess degrees of freedom, which appears to
  work reasonably well in practice.

\item In Section \ref{sec:discussion}, we wrap up with a short
  discussion, and briefly describe extensions of our work to
  heteroskedastic data, and alternative loss functions (other than
  squared loss). 
\end{itemize}

\subsection{Related work}  
\label{sec:related}

There is a lot of work related to the topic of this paper.  In
addition to the classical contributions of
\citet{mallows1973comments,stein1981estimation,   
efron1986biased,efron2004estimation},
 on optimism and degrees of 
freedom, that have already been
discussed, it is worth mentioning \citet{breiman1992little}.  In
Section 2 of this work, the author warns precisely of the downward
bias of SURE for estimating prediction error in regression
models, when the former is evaluated at the model that minimizes SURE
(or here, $C_p$). Breiman was thus keenly aware of excess optimism;
he roughly calculated, for all subsets regression 
with $p$ orthogonal variables, that the SURE-tuned subset   
regression estimator has an approximate excess optimism of 
$0.84p\sigma^2$, in the null case when $\theta_0=0$. 

Several authors have addressed the problem of characterizing the risk  
of an estimator tuned by SURE (or a similar method) by
uniformly controlling the deviations of SURE from its mean over
all tuning parameter values $s\in S$, i.e., by establishing that a
quantity like \smash{$\sup_{s \in  S}|\hRisk_s(Y)-\Risk(\htheta_s)|$},
in our notation, converges to zero in a suitable sense.
Examples of this uniform control strategy are found in
\citet{li1985stein,li1986asymptotic,li1987asymptotic,kneip1994ordered},  
who study linear smoothers; \citet{donoho1995adapting}, who  
study wavelet smoothing; \citet{cavalier2002oracle}, who study
linear inverse problems in sequence space; and \citet{xie2012sure}, 
who study a family of shrinkage estimators in a heteroskedastic
model. Notice that the idea of uniformly controlling the deviations
of SURE away from its mean is quite different in spirit than our
approach, in which we directly seek to understand the gap between 
\smash{$\E[\hRisk_{\hs(Y)}(Y)]$} and \smash{$\Risk(\htheta_{\hs})$}.
It is not clear to us that uniform control of SURE deviations can be
used to precisely understand this gap, i.e., to precisely understand
excess optimism. 

Importantly, the strategy of uniform control can often be used 
to derive so-called oracle inequalities of the form     
\begin{equation}
\label{eq:oracle_ineq}
\Risk(\htheta_{\hs}) \leq (1 + o(1)) \Risk(\htheta_{s_0}),
\end{equation}
Such oracle inequalities are derived in
\citet{li1985stein,li1986asymptotic,li1987asymptotic,
kneip1994ordered,donoho1995adapting,cavalier2002oracle,xie2012sure}.
In Section \ref{sec:oracle}, we will return to the oracle
inequality \eqref{eq:oracle_ineq}, and will show that
\eqref{eq:oracle_ineq} can be established in some cases via a
bound on excess optimism. 

When the data are normally distributed, i.e., when
$F=N(\theta_0,\sigma^2 I)$ in \eqref{eq:data_model}, one might think
to use Stein's formula on the SURE-tuned estimator
\smash{$\htheta_{\hs}$} itself, in order to compute its proper degrees
of freedom, and hence excess optimism. 
This idea is pursued in Section \ref{sec:stein}, where we also show
that implicit differentiation 
%(similar to the work of
%\citet{dijkstra2013ridge} on ridge regression) 
can be applied in order to characterize the excess degrees of freedom,
under some assumptions.  
We must emphasize, however, that these assumptions  
are very strong. Stein's original
work,  \citet{stein1981estimation}, established the result in
\eqref{eq:stein_formula}, when the estimator \smash{$\htheta$} is
continuous and weakly differentiable, as a function of 
$Y$. But,
even when \smash{$\htheta_s$} is itself continuous in $Y$ for each
$s \in S$, it is possible for the SURE-tuned estimator 
\smash{$\htheta_{\hs}$} to be discontinuous in $Y$, and
in these cases, Stein's formula does not apply.
\citet{tibshirani2015degrees} and  
\citet{mikkelsen2016degrees} derive extensions of 
Stein's formula to deal with estimators having (specific types of)  
discontinuities. We leverage these extensions in
Section \ref{sec:stein}.  

A parallel problem is to study the excess optimism associated  
with parameter tuning by cross-validation, considered in
\citet{varma2006bias,tibshirani2009bias,
bernau2013correcting,krstajic2014cross,tsamardinos2015performance}. 
Since it is difficult to study cross-validation mathematically, 
these works do not develop formal characterizations or corrections and
are mostly empirically-driven.  
% An exception is \citet{li1987asymptotic},
% who provides theory for leave-one-out cross-validation applied to  
% subset regression and nearest-neighbor estimators. 
% RJT: This is NOT excess optimism, rather uniform bounds

Lastly, it is worth mentioning that some of the motivation of 
\citet{efron2014estimation} is similar to that in our paper, though
the focus is different: Efron focuses on constructing proper estimates
of standard error (and confidence intervals) for estimators 
that are defined with inherent parameter tuning (he uses the term 
``model selection'' rather than parameter tuning).  Discontinuities 
play a major role in \citet{efron2014estimation}, as they do in
ours (i.e., in our Section \ref{sec:stein}); Efron proposes to
replace parameter-tuned estimators with bagged (bootstrap
aggregated) versions, as the latter estimators are smoother
and can deliver shorter standard errors (or confidence
intervals).  More broadly, post-selection inference,
as studied in \citet{berk2013valid,lockhart2014significance, 
lee2016exact,tibshirani2016exact,fithian2014optimal} and several 
others, is also related in spirit to our work, though our focus is on 
prediction error rather than inference.  While post-selection 
prediction can also be studied from the conditional perspective
that is often used in post-selection inference, this seems to
be less common.  A notable exception is 
\citet{tianharris2016prediction}, who proposes a clever 
randomization scheme for constructing estimates of prediction error
that are conditionally valid on a model selection event, in a
regression setting.

\section{An upper bound on the oracle gap}
\label{sec:oracle}

We derive a simple inequality that relates the error of the estimator  
\smash{$\htheta_{\hs}$}
% where recall \smash{$\hs$} is the rule that selects the the tuning 
% parameter value minimizing SURE, as in \eqref{eq:sure_tuning}, 
to the error of what we may call the {\it oracle} estimator
\smash{$\htheta_{s_0}$}, where $s_0$ is the tuning parameter
value minimizing the (unavailable) 
true prediction error, as in \eqref{eq:oracle_tuning}.  Observe that 
\begin{equation}
\label{eq:simple_bd}
\E[\hErr_{\hs(Y)}(Y)] = \E\Big(\min_{s\in S} \, \hErr_s(Y) \Big)
\leq \min_{s \in S} \, \E[\hErr_s(Y)] = \min_{s \in S} \,
\Err(\htheta_s) =  \Err(\htheta_{s_0}). 
\end{equation}
By subtracting the left- and right-most expressions from
\smash{$\Err(\htheta_{\hs})$}, the true prediction error of
\smash{$\htheta_{\hs}$}, we have established the following result.

\begin{theorem}
\label{thm:oracle_bd}
For any family of estimators \smash{$\{\htheta_s : s \in S\}$}, it 
holds that 
\begin{equation}
\label{eq:oracle_bd}
\Err(\htheta_{\hs}) \leq \Err(\htheta_{s_0}) + \ExOpt(\htheta_{\hs}).
\end{equation}
Here, \smash{$\hs$} is the tuning parameter rule defined by
minimizing SURE, as in \eqref{eq:sure_tuning}, $s_0$ is the oracle
tuning parameter value minimizing prediction error, as in
\eqref{eq:oracle_tuning}, and \smash{$\ExOpt(\htheta_{\hs})$} is the
excess optimism, as defined in \eqref{eq:ex_opt}. 
\end{theorem}

Theorem \ref{thm:oracle_bd} says that the excess optimism, which is 
a quantity that we can in principle calculate (or at least, estimate),
serves as an upper bound for the gap between the prediction error of
\smash{$\htheta_{\hs}$} and the oracle error.  This gives an
interesting, alternative motivation for  excess optimism to that given
in the introduction: excess optimism tells us how far the SURE-tuned
estimator \smash{$\htheta_{\hs}$} can be from the best member of the
class \smash{$\{\htheta_s : s \in   S\}$}, in terms of prediction
error.  A few remarks are in order. 

\begin{remark}[\textbf{Risk inequality}]
Recalling that excess optimism can be equivalently posed in terms of   
risk, as in \eqref{eq:ex_opt_risk}, the bound in
\eqref{eq:oracle_bd} can also be written in terms of risk, namely,  
\begin{equation}
\label{eq:oracle_bd_risk}
\Risk(\htheta_{\hs}) \leq \Risk(\htheta_{s_0}) +
  \ExOpt(\htheta_{\hs}),
\end{equation}
which says the excess risk \smash{$\Risk(\htheta_{\hs}) -
  \Risk(\htheta_{s_0})$} of the SURE-tuned estimator is upper bounded
by its excess optimism, \smash{$\ExOpt(\htheta_{\hs})$}.
If we can show that this excess optimism is small compared to the
oracle risk, in particular, if we can show that
\smash{$\ExOpt(\htheta_{\hs}) =  o(\Risk(\htheta_{s_0}))$}, then 
\eqref{eq:oracle_bd_risk} implies the oracle inequality
\eqref{eq:oracle_ineq}. We will revisit this idea in Sections
\ref{sec:shrinkage} and \ref{sec:subset_reg}.
\end{remark}

\begin{remark}[\textbf{Beating the oracle}]
If \smash{$\ExOpt(\htheta_{\hs}) < 0$}, then \eqref{eq:oracle_bd}
implies \smash{$\htheta_{\hs}$} outperforms the
oracle, in terms of prediction error (or risk). Technically this is
not impossible, as \smash{$\theta_{s_0}$} is the optimal
fixed-parameter estimator, in the class \smash{$\{\theta_s : s \in
  S\}$}, whereas \smash{$\htheta_{\hs}$} is tuned in a data-dependent
fashion. But it seems unlikely to us that excess optimism can be
negative, recall Section \ref{sec:ex_opt_nneg}. 
\end{remark}

\begin{remark}[\textbf{Beyond SURE}]
The argument in \eqref{eq:simple_bd} and thus the validity of Theorem   
\ref{thm:oracle_bd} only used the fact that \smash{$\hs$} was
defined by minimizing an unbiased estimator of prediction error, and 
SURE is not the only such estimator.  For example, the 
result in Theorem \ref{thm:oracle_bd} applies to the standard
hold-out estimator of prediction error, when hold-out data $Y^* 
\sim F$ (independent of $Y$) is available. While the result does not 
exactly carry over to cross-validation (since the standard
cross-validation estimator of prediction error is not unbiased in
finite samples, at least not without additional corrections and
assumptions), we can think of it as being true in some approximate 
sense.  
\end{remark}

\section{Shrinkage estimators}
\label{sec:shrinkage}

In this section, we focus on shrinkage estimators, and consider normal 
data, $Y \sim F=N(\theta_0, \sigma^2 I)$ in \eqref{eq:data_model}. 
Due to the simple form of the family of shrinkage estimators (and 
the normality assumption), we can compute an (exact) unbiased 
estimator of excess degrees of freedom, and excess optimism.  

\subsection{Shrinkage in normal means} 
\label{sec:shrink_means}

First, we consider the simple family of shrinkage estimators
\begin{equation}
\label{eq:est_shrink}
\htheta_s(Y) = \frac{Y}{1+s},
\quad \text{for} \; s \geq 0.
\end{equation}
In this case, we can see that SURE in \eqref{eq:sure} is 
\begin{equation}
\label{eq:sure_shrink}
\hErr_s(Y) = \|Y\|_2^2 \frac{s^2}{(1+s)^2} + 
{2\sigma^2}\frac{n}{1+s}.
\end{equation}
Here we have used exact calculation (rather than an unbiased
estimate) for the degrees of freedom,
\smash{$\df(\htheta_s)=n/(1+s)$}.  The next lemma characterizes
\smash{$\hs$}, the mapping defined by the minimizer of
\eqref{eq:sure_shrink}.  The proof is elementary and delayed until
the appendix.   

\begin{lemma}
\label{lem:sure_min_shrink}
Define $g(x) = a x^2/(1+x)^2 + 2b/(1+x)$, where $a,b > 0$.  Then the 
minimizer of $g$ over $x \geq 0$ is
\begin{equation*}
x^* = \begin{cases}
\frac{b}{a-b} & \text{if} \; a \geq b \\
\infty & \text{if} \; a < b.
\end{cases}
\end{equation*}
\end{lemma}

According to Lemma \ref{lem:sure_min_shrink}, the rule
\smash{$\hs$} defined by minimizing \eqref{eq:sure_shrink} is 
\begin{equation*}
\hs(Y) = \begin{cases}
\displaystyle 
\frac{n\sigma^2}{\|Y\|_2^2 - n\sigma^2} & 
\text{if} \; \|Y\|_2^2 \geq n\sigma^2 \\ 
\infty & \text{if} \; \|Y\|_2^2 < n\sigma^2. 
\end{cases}
\end{equation*}
Plugging this in gives the SURE-tuned shrinkage estimate
\smash{$\htheta_{\hs(Y)}(Y) = Y/(1+\hs(Y))$}.  This is continuous and   
weakly differentiable as a function of $Y$, and hence by Stein's
formula \eqref{eq:stein_formula}, we can form an unbiased estimator    
of its degrees of freedom by computing its divergence.  When
\smash{$\hs(Y) < \infty$}, the divergence of \smash{$\htheta_{\hs}$}
at $Y$ is
\begin{align}
\nonumber
\frac{n}{1+\hs(Y)} -
\sum_{i=1}^n \frac{Y_i}{(1+\hs(Y))^2} 
\frac{\partial \hs}{\partial Y_i}(Y)  
&= \frac{n}{1+\hs(Y)} +
\sum_{i=1}^n \frac{Y_i}{(1+\hs(Y))^2} 
\frac{n\sigma^2}{(\|Y\|_2^2 - n\sigma^2)^2} 2Y_i \\
\label{eq:sure_div_shrink}
&= \frac{n}{1+\hs(Y)} + 
\frac{2\hs(Y)}{1+\hs(Y)}.
\end{align}
When \smash{$\hs(Y)=\infty$}, the divergence is 0.

Hence, we can see directly that for the SURE-tuned shrinkage
estimator \smash{$\htheta_{\hs}$}, we have the excess degrees of
freedom bound
\begin{equation}
\label{eq:ex_df_shrink}
\edf(\htheta_{\hs}) = 
\E\bigg(\frac{2\hs(Y)}{1+\hs(Y)} \, ; \, \hs(Y) < \infty    
\bigg) \leq 2,
\end{equation}
and so \smash{$\ExOpt(\htheta_{\hs}) \leq 4\sigma^2$}. 
% We can alternatively express the excess degrees of freedom as  
% \begin{equation*}
% \edf(\htheta_{\hs}) =  2 \E(n/W \,;\, W \geq n),
% \end{equation*}
% where $W$ is a noncentral $\chi^2$ random variate with $n$
% degrees of freedom, and mean $\E\|W\|_2^2 = \|\theta_0\|_2^2 + 
% n$ (noncentrality parameter $\|\theta_0\|_2^2$).
A lot is known about shrinkage estimators in the current normal  
means problem that we are considering, dating back to the  
seminal work of \citet{james1961estimation}; some excellent recent
references are Chapter 1 of \citet{efron2010large}, and Chapter 2 of
\citet{johnstone2015gaussian}. It is easy to show that
the oracle choice of tuning parameter in the current setting is 
\smash{$s_0=n\sigma^2/\|\theta_0\|_2^2$}, thus    
\begin{equation}
\label{eq:oracle_risk_shrink}
\Risk(\htheta_{s_0}) = \frac{n\sigma^2 \|\theta_0\|_2^2}
{n\sigma^2 + \|\theta_0\|_2^2}. 
\end{equation}
By our excess optimism bound of $4\sigma^2$, and Theorem
\ref{thm:oracle_bd} (actually, \eqref{eq:oracle_bd_risk}, the risk 
version of the result in the theorem), the risk of the 
SURE-tuned shrinkage estimator \smash{$\htheta_{\hs}$} satisfies  
\begin{equation}
\label{eq:risk_bd_shrink}
\Risk(\htheta_{\hs}) \leq \frac{n\sigma^2 \|\theta_0\|_2^2} 
{n\sigma^2 + \|\theta_0\|_2^2} + 4\sigma^2. 
\end{equation}

\begin{remark}[\textbf{Oracle inequality for SURE-tuned shrinkage}]
For large \smash{$\|\theta_0\|_2^2$}, the risk gap 
of $4\sigma^2$ for the SURE-tuned shrinkage estimator is negligible
next to the oracle risk in \eqref{eq:oracle_risk_shrink}. Specifically,
if \smash{$\|\theta_0\|_2^2 \to \infty$} as $n \to \infty$ (with
$\sigma^2$ held constant), then we see that
\eqref{eq:risk_bd_shrink} implies the oracle inequality
\eqref{eq:oracle_ineq} for the SURE-tuned shrinkage estimator. 
\end{remark}

\subsection{Interlude: James-Stein estimation}  

The SURE-tuned shrinkage estimator of the last subsection can be 
written as 
\begin{equation*}
\htheta_{\hs(Y)}(Y) = \begin{cases}
\displaystyle 
\frac{1}{1+\frac{n\sigma^2}{\|Y\|_2^2 - n\sigma^2}} 
Y & \text{if} \; \|Y\|_2^2 \geq n\sigma^2 \\ 
0 & \text{if} \; \|Y\|_2^2 < n\sigma^2,
\end{cases}
\end{equation*}
or more concisely, as
\begin{equation}
\label{eq:sure_tuned_shrink}
\htheta_{\hs(Y)}(Y) = 
\bigg(1 - \frac{n\sigma^2}{\|Y\|_2^2}\bigg)_+ Y,
\end{equation}
where we write $x_+=\max\{x,0\}$ for the positive part of $x$.   
Meanwhile, the positive part James-Stein estimator
\citep{james1961estimation,baranchik1964multiple} is defined as   
\begin{equation}
\label{eq:james_stein}
\htheta^{\mathrm{JS+}}(Y) = 
\bigg(1 - \frac{(n-2)\sigma^2}{\|Y\|_2^2}\bigg)_+ Y,
\end{equation}
so the two estimators \eqref{eq:sure_tuned_shrink} and
\eqref{eq:james_stein} only differ by the appearance of $n$ versus
$n-2$ in the shrinkage factor.  This connection---between SURE-tuned 
shrinkage estimation and positive part James-Stein estimation---seems
to be not very well-known, and was a surprise to us; after writing an
initial draft of this paper, we found that this fact was mentioned in
passing in \citet{xie2012sure}.   We now give a few remarks.

\begin{remark}[\textbf{Dominating the MLE}]
It can be shown that the SURE-tuned shrinkage estimator 
in \eqref{eq:sure_tuned_shrink} 
dominates the MLE, i.e., \smash{$\htheta^{\mathrm{MLE}}(Y)=Y$},
just like the positive part James-Stein estimator in
\eqref{eq:james_stein}.  For this to be true of the former estimator,
we require $n \geq 5$, while the latter estimator only requires $n
\geq 3$.  

Our proof of \smash{$\htheta_{\hs}$} dominating
\smash{$\htheta^{\mathrm{MLE}}$} mimicks Stein's elegant proof for the
James-Stein estimator, \citep{stein1981estimation}.  Consider
SURE for \smash{$\htheta_{\hs}$}, which gives an unbiased
estimator of the risk of \smash{$\htheta_{\hs}$}, provided we
compute its divergence properly, as in \eqref{eq:sure_div_shrink}.
Write \smash{$\hat{R}$} for this unbiased risk estimator.
If \smash{$\hs(Y) < \infty$}, i.e., $\|Y\|_2^2 \geq n\sigma^2$, then 
\begin{align*}
\hat{R}(Y) &= -n\sigma^2 +
\frac{\hs(Y)^2}{(1+\hs(Y))^2} \|Y\|_2^2 + 
2\sigma^2 \bigg(\frac{n}{1+\hs(Y)} +  
\frac{2\hs(Y)}{1+\hs(Y)}\bigg) \\
&= -n\sigma^2 + \frac{(n\sigma^2)^2}{\|Y\|_2^2} +
 2n\sigma^2 \frac{\|Y\|_2^2-n\sigma^2}{\|Y\|_2^2} +
4\sigma^2 \frac{n\sigma^2}{\|Y\|_2^2} \\
&= n\sigma^2 - (n-4)\sigma^2 \frac{n\sigma^2}{\|Y\|_2^2} 
< n\sigma^2.
\end{align*}
If $\hs(Y)=\infty$, i.e., $\|Y\|_2^2 < n\sigma^2$, then
we have \smash{$\hat{R}(Y) = -n\sigma^2 + \|Y\|_2^2 < 0$}. Taking an 
expectation, we thus see that
\smash{$\Err(\htheta_{\hs}) = \E[\hat{R}(Y)]< n\sigma^2$}, which
establishes the result, as $n\sigma^2$ is the risk of the MLE.
\end{remark}

\begin{remark}[\textbf{Risk of positive part James-Stein}] 
A straightforward calculation, similar to that given above for
\smash{$\htheta_{\hs}$} (see also Theorem 5 of
\citet{donoho1995adapting}) shows that the risk of the positive part
James-Stein estimator satisfies
\begin{equation}
\label{eq:risk_bd_js}
\Risk(\htheta^{\mathrm{JS+}}) \leq \frac{n\sigma^2 \|\theta_0\|_2^2} 
{n\sigma^2 + \|\theta_0\|_2^2} + 2\sigma^2, 
\end{equation}
so it admits an even tighter gap to the oracle risk than does the
SURE-tuned shrinkage estimator, recalling \eqref{eq:risk_bd_shrink}.   
\end{remark}

\begin{remark}[\textbf{Inadmissibility of the SURE-tuned shrinkage 
    estimator}]
Comparing \eqref{eq:risk_bd_js} and \eqref{eq:risk_bd_shrink} suggests  
that the positive part James-Stein estimator might have better risk
than the SURE-tuned shrinkage estimator.  This is indeed true, in the
strongest sense possible, as it can be shown that
\smash{$\htheta^{\mathrm{JS+}}$} dominates \smash{$\htheta_{\hs}$};
the proof simply follows the same arguments as those given above 
for the proof of \smash{$\htheta_{\hs}$} dominating the MLE.  (Also,
the positive part James-Stein estimator is itself dominated by others,
see, e.g., \citet{shao1994improving}.)

It is worth noting that \smash{$\htheta_s$} is itself admissible, for
every fixed tuning parameter value $s \geq 0$, since it is the unique
Bayes estimator under the normal prior $\theta_0 \sim
N(0,s^{-1}I)$. That \smash{$\htheta_{\hs}$} is
inadmissible---which is defined at each $Y$ by minimizing an
unbiased estimate of risk over the family of admissible estimators
\smash{$\{\htheta_s : s \geq 0\}$}---is therefore perhaps surprising.      
\end{remark}

\subsection{Shrinkage in regression}
\label{sec:shrink_reg}

Now, we consider the family of regression shrinkage estimators  
\begin{equation}
\label{eq:est_shrink_reg}
\htheta_s(Y) = \frac{P_X Y}{1+s},
\quad \text{for} \; s \geq 0,
\end{equation}
where we write $P_X \in \R^{n\times n}$ for the projection matrix onto the 
column space of a predictor matrix $X \in \R^{n\times p}$, i.e.,
$P_X=X(X^T X)^{-1} X^T$ if $X$ has full column rank, and $P_X=X(X^T
X)^+ X^T$ otherwise (here and throughout, $A^+$ denotes the
pseudoinverse of a matrix $A$). 

Treating $X$ as fixed (nonrandom), it is easy to check that SURE
\eqref{eq:sure} for our regression shrinkage estimator is
\begin{equation}
\label{eq:sure_shrink_reg}
\hErr_s(Y) = \|P_X Y\|_2^2 \frac{s^2}{(1+s)^2} + 
{2\sigma^2}\frac{r}{1+s},
\end{equation}
where $r=\rank(X)$, the rank of $X$.  This is directly analogous to
\eqref{eq:sure_shrink} in the normal means setting, and Lemma
\ref{lem:sure_min_shrink} shows that the minimizer
\smash{$\hs$} of \eqref{eq:sure_shrink_reg} is defined by
\begin{equation*}
\hs(Y) = \begin{cases}
\displaystyle 
\frac{r\sigma^2}{\|P_X Y\|_2^2 - r\sigma^2} & 
\text{if} \; \|P_X Y\|_2^2 \geq r\sigma^2 \\ 
\infty & \text{if} \; \|P_X Y\|_2^2 < r\sigma^2. 
\end{cases}
\end{equation*}
The same arguments as in Section \ref{sec:shrink_means} then lead to 
the same excess degrees of freedom bound 
\begin{equation}
\label{eq:ex_df_shrink_reg}
\edf(\htheta_{\hs}) = 
\E\bigg(\frac{2\hs(Y)}{1+\hs(Y)} \, ; \, \hs(Y) < \infty  
\bigg) \leq 2,
\end{equation}
thus \smash{$\ExOpt(\htheta_{\hs}) \leq 4\sigma^2$}.
By direct calculation, the oracle tuning parameter is
\smash{$s_0 = r\sigma^2 / \|P_X \theta_0\|_2^2$}, and now
\begin{equation}
\label{eq:oracle_risk_shrink_reg}
\Risk(\htheta_{s_0}) = \frac{r\sigma^2 \|\theta_0\|_2^2
+ \|P_X \theta_0\|_2^2 (\|\theta_0\|_2^2 - \|P_X \theta_0\|_2^2)} 
{r\sigma^2 + \|P_X \theta_0\|_2^2}. 
\end{equation}
Combining our excess optimism bound of $4\sigma^2$ with  
Theorem \ref{thm:oracle_bd} (i.e., combining
it with \eqref{eq:oracle_bd_risk},  
the risk version of the result in the theorem), we have 
\begin{equation}
\label{eq:risk_bd_shrink_reg}
\Risk(\htheta_{\hs}) \leq \frac{r\sigma^2 \|\theta_0\|_2^2 +
 \|P_X \theta_0\|_2^2 (\|\theta_0\|_2^2 - \|P_X 
  \theta_0\|_2^2)}{r\sigma^2 + \|P_X \theta_0\|_2^2} + 4\sigma^2.   
\end{equation}

\begin{remark}[\textbf{Oracle inequality for SURE-tuned regression
    shrinkage}] The risk gap of $4\sigma^2$, 
  for the SURE-tuned regression  
  shrinkage estimator, will be negligible next to the oracle risk
  \eqref{eq:oracle_risk_shrink_reg} under various sufficient
  conditions.  For example, if 
  \smash{$\|\theta_0\|_2^2 \to \infty$} and {$\|P_X\theta_0\|_2^2
    |\|\theta_0\|_2^2 - \|P_X\theta_0\|_2^2|= O(r)$} as $n,r \to
  \infty$ (and $\sigma^2$ is held constant), then it is not hard to
  check that \eqref{eq:risk_bd_shrink_reg} implies the oracle
  inequality \eqref{eq:oracle_ineq} for the SURE-tuned regression
  shrinkage estimator.   
\end{remark}

\subsection{Interlude: James-Stein and ridge regression}

The SURE-tuned regression shrinkage estimator of the previous
subsection can be expressed as
\begin{equation}
\label{eq:sure_tuned_shrink_reg}
\htheta_{\hs}(Y) =
\bigg(1 - \frac{r \sigma^2}{\|P_X Y\|_2^2}\bigg)_+ P_X Y,
\end{equation}
which resembles the positive part James-Stein regression estimator 
\begin{equation}
\label{eq:james_stein_reg}
\htheta^{\mathrm{JS+}}(Y) =
\bigg(1 - \frac{(r-2) \sigma^2}{\|P_X Y\|_2^2}\bigg)_+ P_X Y.
\end{equation}
The same properties as before, of \eqref{eq:sure_tuned_shrink_reg}
dominating the MLE (i.e., the least squares regression estimator), 
\smash{$\htheta^{\mathrm{MLE}}(Y) = P_X Y$}, and also
\eqref{eq:james_stein_reg} dominating
\eqref{eq:sure_tuned_shrink_reg}, carry over to the current setting.   

We point out a connection to
penalized regression.  For any fixed tuning parameter value $s \geq
0$, we can express the estimate in \eqref{eq:est_shrink_reg} as 
\smash{$\htheta_s (Y)= X\hbeta_s(Y)$}, where \smash{$\hbeta_s(Y)$}
solves the convex (though not necessarily strictly convex) penalized 
regression problem, 
\begin{equation}
\label{eq:x_ridge}
\hbeta_s(Y) \in \argmin_{\beta \in \R^p} \; 
\half \|Y - X\beta\|_2^2 + s \|X\beta\|_2^2.
\end{equation}
Hence an alternative interpretation for the estimator
\smash{$\htheta_{\hs}$} in \eqref{eq:sure_tuned_shrink_reg} (whose 
close cousin is the positive part James-Stein regression estimator
\smash{$\htheta^{\mathrm{JS+}}$} in \eqref{eq:james_stein_reg}) is
that we are using SURE to select the tuning parameter over 
the family of penalized regression estimators in \eqref{eq:x_ridge},
for $s \geq 0$.  This has the precise risk guarantee in
\eqref{eq:risk_bd_shrink_reg} (and \smash{$\htheta^{\mathrm{JS+}}$} 
enjoys an even stronger guarantee, with $2\sigma^2$ in place of
$4\sigma^2$).

Compared to \eqref{eq:x_ridge}, a more familiar penalized regression
problem to most statisticians is perhaps the ridge regression problem 
\citep{hoerl1970ridge}, 
\begin{equation}
\label{eq:ridge}
\hbeta^{\mathrm{ridge}}_s(Y) = \argmin_{\beta \in \R^p} \; 
\half \|Y - X\beta\|_2^2 + s \|\beta\|_2^2.
\end{equation}
Several differences between \eqref{eq:x_ridge} and 
\eqref{eq:ridge}  
can be enumerated; one interesting difference is that the solution in
the former problem shrinks uniformly across all dimensions
$1,\ldots,p$, whereas that in the latter problem shrinks less in 
directions of high variance and more in directions of low
variance, defined with respect to the predictor variables (i.e.,
shrinks less in the top eigendirections of $X^T X$).

It is generally accepted that neither regression shrinkage estimator, 
in \eqref{eq:x_ridge} and \eqref{eq:ridge}, is better than the
other.\footnote{It is worth pointing out that the former problem
\eqref{eq:x_ridge} does not give a well-defined, i.e., unique solution
for the coefficients when $\rank(X)<p$, and the latter problem 
\eqref{eq:ridge} does, when $s>0$.}  But, we have seen that  
SURE-tuning in the first problem \eqref{eq:x_ridge} provides us with 
an estimator  
\smash{$\htheta_{\hs}=X\hbeta_{\hs}$} that has a definitive risk
guarantee \eqref{eq:risk_bd_shrink_reg} and provably dominates the  
MLE.  The story for ridge regression is less clear; to quote
\citet{efron2016computer}, Chapter 7.3: ``{\it There is no 
  [analogous] guarantee for ridge regression, and no foolproof way to
  choose the ridge parameter.}''  Of course, if we could bound the
excess degrees of freedom for SURE-tuned ridge regression, then this
could lead (depending on the size of the bound) to a useful risk
guarantee, providing some rigorous backing to SURE tuning for ridge
regression.  However, characterizing excess degrees of freedom 
for ridge regression is far from straightforward, as we remark next. 

\begin{remark}[\textbf{Difficulties in analyzing excess degrees
    of freedom for SURE-tuned ridge regression}]
While it may seem tempting to analyze the risk of the SURE-tuned 
ridge regression estimator, \smash{$\htheta^{\mathrm{ridge}}_{\hs} = 
  X\hbeta^{\mathrm{ridge}}_{\hs}$} (where \smash{$\hs$} is the 
SURE-optimal ridge parameter map), using arguments that mimick    
those we gave above for the SURE-tuned shrinkage estimator 
\smash{$\htheta_{\hs} =  
  X\hbeta_{\hs}$}, this is not an easy task. 
When $X$ is orthogonal, the two estimators
\smash{$\htheta_s$}, \smash{$\htheta^{\mathrm{ridge}}_s$} are
exactly the same, 
for all $s \geq 0$, hence our previous analysis already covers the
SURE-tuned ridge regression estimator
\smash{$\htheta^{\mathrm{ridge}}_{\hs}$}. But for a general $X$,
the story is far more complicated, for two reasons: (i) the
SURE-optimal tuning parameter map \smash{$\hs$} is not available in
closed form for ridge regression, and (ii) the SURE-tuned ridge
estimator \smash{$\htheta^{\mathrm{ridge}}_{\hs}$} is not necessarily    
continuous with respect to the data $Y$, so Stein's formula   
cannot be used to compute an unbiased estimator of its
degrees of freedom. (Specifically, it is unclear
if the SURE-optimal ridge parameter map \smash{$\hs$} is itself  
continuous with respect to $Y$, as it is defined by the minimizer of a 
possibly multimodal SURE criterion; see Figure \ref{fig:s_jump}.)   

It is really the second reason, i.e., (possible) disconinuities in 
\smash{$\htheta^{\mathrm{ridge}}_{\hs}$}, that makes the
analysis so complicated. Even when \smash{$\hs$} cannot be expressed
in closed form, implicit differentiation can be used to compute the 
divergence of \smash{$\htheta^{\mathrm{ridge}}_{\hs}$}, as we explain
in Section \ref{sec:stein_char}; in the presence of discontinuities,
however, this divergence will not be enough to characterize the
degrees of freedom (and thus excess degrees of freedom) of 
\smash{$\htheta^{\mathrm{ridge}}_{\hs}$}. Extensions of Stein's
divergence formula from \citet{tibshirani2015degrees} and
\citet{mikkelsen2016degrees} can be used to characterize degrees
of freedom for estimators having certain types of discontinuities,
which we review in Section \ref{sec:stein_ext_char}. Generally 
speaking, these extensions require complicated
calculations. Later, in Section \ref{sec:hetero}, we revisit the
ridge regression problem, and we compute the divergence of the
SURE-tuned ridge estimator via implicit differentiation, but we leave
proper treatment of discontinuties to future work.  
\end{remark}

% Further extensions can be considered, for example, one can define a
% family of regression estimators 
% \begin{equation*}
% \htheta_s(Y) = \frac{1}{1+s} P_X Y
% + \frac{s}{1+s} \bar{Y} \mathds{1}, 
% \quad \text{for} \; s \geq 0,
% \end{equation*}
% where \smash{$\bar{Y} = \sum_{i=1}^n Y_i /n$}, and 
% $\mathds{1} = (1,1,\ldots,1) \in \R^n$ is the vector of all 1s.  In
% words, we are shrinking now towards the sample mean, instead of
% towards zero. 
% Essentially all results carry over to this setting;
% \smash{$\htheta_{\hs}$} now dominates the MLE whenever $p>5$. We 
% omit details. 

\section{Subset regression estimators}
\label{sec:subset_reg}

Here we study subset regression estimators, and again
consider normal data, $Y \sim F=N(\theta_0, \sigma^2 I)$ in
\eqref{eq:data_model}.  Our family of estimators is defined by
regression onto subsets of the columns of a predictor matrix $X \in
\R^{n\times p}$, i.e.,  
\begin{equation}
\label{eq:est_subset_reg}
\htheta_s(Y) = P_{X_s} Y \quad \text{for} \; s \in S,
\end{equation}
where each \smash{$s=\{j_1,\ldots,j_{p_s}\}$} is an arbitrary subset
of $\{1,\ldots,p\}$ of size $p_s$, \smash{$X_s \in \R^{n\times p_s}$}
denotes the columns of $X$ indexed by elements of $s$,
\smash{$P_{X_s}$} denotes the projection matrix onto the column space
of $X_s$, and $S$ denotes a collection of subsets of 
$\{1,\ldots,p\}$. We will abbreviate \smash{$P_s=P_{X_s}$}, and we
will assume, without any real loss of generality, that for each $s\in
S$, the matrix $X_s$ has full column rank (otherwise, simply replace
each instance of $p_s$ below with $r_s = \rank(X_s)$).     

% We do not place assumptions up front on the data distribution in  
% \eqref{eq:data_model}, but will present specific distributional
% assumptions as needed below.  For now, under the generic data 
% model in \eqref{eq:data_model},  
SURE in \eqref{eq:sure} is now the familiar $C_p$ criterion 
\begin{equation}
\label{eq:sure_subset_reg}
\hErr_s(Y) = \|Y - P_s Y\|_2^2 + 2\sigma^2 p_s.
\end{equation}
As $S$ is discrete, it is not generally possible to express the
minimizer \smash{$\hs(Y)$} of the above criterion in closed form, and
so, unlike the previous section, not generally possible to
analytically characterize the excess degrees of freedom of the
SURE-tuned subset regression estimator \smash{$\htheta_{\hs}$}. In
what follows, we derive an upper bound on the excess degrees  
of freedom, using elementary arguments. Later in 
Section \ref{sec:subset_reg_mh}, we give a lower bound and a
more sophisticated upper bound, by leveraging a powerful
tool from \citet{mikkelsen2016degrees}.    

\subsection{Upper bounds for excess degrees of freedom  
in subset regression} 
\label{sec:subset_reg_upper}

Note that we can write the excess degrees of freedom as
\begin{equation}
\label{eq:ex_df_subset_reg}
\edf(\htheta_{\hs}) = \frac{1}{\sigma^2}\E\big[
(P_{\hs(Y)}(Y))^T (Y-\theta_0) \big] - \E(p_{\hs(Y)}) 
= \frac{1}{\sigma^2}\E \|P_{\hs(Y)} Z\|_2^2 - \E(p_{\hs(Y)}),  
\end{equation}
where $Z=Y-\theta_0$ has mean zero and covariance $\sigma^2
I$. Furthermore, by defining $W_s=\|P_s Z\|_2^2/\sigma^2$ for $s \in
S$, we have
\begin{equation}
\label{eq:ex_df_subset_reg_ub_1}
\edf(\htheta_{\hs}) = \E (W_{\hs(Y)} - p_{\hs(Y)}) \leq
\E \Big[ \max_{s \in S} \, (W_s - p_s) \Big].
\end{equation}
As $Y \sim N(\theta_0,\sigma^2 I)$, we have 
$W_s \sim \chi_{p_s}^2$ 
for each $s\in S$, and the next lemma provides a useful
upper bound for the right-hand side above.
Its proof is given in the appendix.   

\begin{lemma}
\label{lem:chi_sq_max}
Let \smash{$W_s \sim \chi^2_{p_s}$}, $s \in S$.  This collection need
not be independent. Then for any $0 \leq \delta < 1$, 
\begin{equation}
\label{eq:chi_sq_max}
\E \Big[ \max_{s \in S} \, (W_s - p_s) \Big] \leq 
\frac{2}{1-\delta} \log \sum_{s \in S} (\delta e^{1-\delta})^{-p_s/2}.
\end{equation}
\end{lemma}

It is worth noting that the proof of the Lemma \ref{lem:chi_sq_max}   
relies only on the moment generating function of the
chi-squared distribution, and therefore our assumption of  
normality for the data $Y$ could be weakened.
%by directly assuming a certain moment generating function bound 
%on $W_s$, $s \in S$.  
For example, it a similar result to that in Lemma \ref{lem:chi_sq_max}
can be derived when each $W_s$, $s \in S$ is subexponential
(generalizing the chi-squared assumption).  For simplicity, we do not
pursue this.     

Combining \eqref{eq:ex_df_subset_reg_ub_1}, \eqref{eq:chi_sq_max}  
gives an upper bound on the excess
degrees of freedom of \smash{$\htheta_{\hs}$},
\begin{equation}
\label{eq:ex_df_subset_reg_ub_2}
\edf(\htheta_{\hs}) \leq \frac{2}{1-\delta} \log \sum_{s \in S}  
(\delta e^{1-\delta})^{-p_s/2}.
\end{equation}
To make this more explicit, we denote by $|S|$ the size of $S$, and
\smash{$p_{\max}=\max_{s\in S}\, p_s$}, and consider a simple
upper bound for the right-hand side in
\eqref{eq:ex_df_subset_reg_ub_2},  
\begin{equation}
\label{eq:ex_df_subset_reg_ub_3}
\edf(\htheta_{\hs}) \leq \frac{2}{1-\delta} \log|S| + 
p_{\max}\bigg(\frac{\log(1/\delta)}{1-\delta} - 1\bigg).
\end{equation}
This simplification should be fairly tight, i.e., the right-hand
side in \eqref{eq:ex_df_subset_reg_ub_3} should be close to that  
in \eqref{eq:ex_df_subset_reg_ub_2}, when $|S|$ and  
\smash{$\max_{s\in S}\, p_s - \min_{s\in S}\, p_s$} are both not very
large. Now, any choice of $0 \leq \delta < 1$ can be used to
give a valid bound in \eqref{eq:ex_df_subset_reg_ub_3}.
As an example, taking $\delta=9/10$ gives
\begin{equation*}
\edf(\htheta_{\hs}) \leq 20 \log|S| + 0.054 p_{\max}. 
\end{equation*}
By \eqref{eq:oracle_bd_risk}, the risk reformulation of the result in
Theorem \ref{thm:oracle_bd}, we get the finite-sample risk bound   
\begin{equation*}
\Risk(\htheta_{\hs}) \leq \|(I-P_{s_0})\theta_0\|_2^2 + 
\sigma^2 (p_{s_0} + 0.108 p_{\max}) + 20 \log|S|,
\end{equation*}
where we have explicitly written the oracle risk as
\smash{$\Risk(\htheta_{s_0}) = \|(I-P_{s_0})\theta_0\|_2^2 + \sigma^2
  p_{s_0}$}.

\subsection{Oracle inequality for SURE-tuned subset regression}  
\label{sec:subset_reg_oracle} 

The optimal choice of $\delta$, i.e., the choice
giving the tightest bound in \eqref{eq:ex_df_subset_reg_ub_3} (and
so, the tightest risk bound), will 
depend on $|S|$ and \smash{$p_{\max}$}. The 
analytic form of such a value of $\delta$ is not clear, given the
somewhat complicated nature of the bound in 
\eqref{eq:ex_df_subset_reg_ub_3}. But, we can adopt an
asymptotic   
perpsective: if $\log|S|$ is small compared to the oracle risk
\smash{$\Risk(\htheta_{s_0})$}, and \smash{$p_{\max}$} is not too 
large compared to the oracle risk, then 
\eqref{eq:ex_df_subset_reg_ub_3} implies 
\smash{$\edf(\htheta_{\hs}) = o(\Risk(\htheta_{s_0}))$}. We state this
formally next, leaving the proof to the appendix. 

\begin{theorem}
\label{thm:ex_df_subset_reg}
Assume that $Y \sim N(\theta_0, \sigma^2 I)$, and that
there is a sequence $a_n > 0$, $n=1,2,3,\ldots$ with $a_n \to   
0$ as $n \to \infty$, such that the risk of the oracle subset
regression estimator \smash{$\htheta_{s_0}$} satisfies
\begin{equation}
\label{eq:ex_df_subset_reg_ass}
\frac{1}{a_n} \frac{\log|S|}{\Risk(\htheta_{s_0})} \to 0
\quad \text{and} \quad 
a_n \frac{p_{\max}}{\Risk(\htheta_{s_0})} \to 0
\quad \text{as} \; n \to \infty.
\end{equation}
Then there is a sequence $0 \leq \delta_n < 1$, $n=1,2,3,\ldots$
with 
$\delta_n \to 1$ as $n \to \infty$, such that 
\begin{equation*}
\bigg[\frac{2}{1-\delta_n} \log|S| + 
p_{\max} \bigg( \frac{\log(1/\delta_n)}{1-\delta_n} - 1 \bigg) \bigg]
/ \Risk(\htheta_{s_0}) \to 0
\quad \text{as} \; n \to \infty.
\end{equation*}
Plugging this into the bound in \eqref{eq:ex_df_subset_reg_ub_3}
shows that \smash{$\edf(\htheta_{\hs})/\Risk(\htheta_{s_0}) \to 0$}, so 
\smash{$\ExOpt(\htheta_{\hs})/\Risk(\htheta_{s_0}) \to 0$} as well,
establishing the oracle inequality \eqref{eq:oracle_ineq} for the
SURE-tuned subset regression estimator. 
\end{theorem}

The assumptions in \eqref{eq:ex_df_subset_reg_ass} may look abstract,
but are not strong and satisfied under fairly simple conditions.  For
example, if we assume that
\smash{$\|(I-P_{s_0})\theta_0\|_2^2=0$} (which means there is no  
bias), and as $n \to \infty$ (with $\sigma^2$ constant) it holds
that $(\log|S|)/p_{s_0} \to 0$ and $p_{\max}/p_{s_0}=O(1)$ (which
means the number $|S|$ of candidate models is much smaller than  
\smash{$2^{p_0}$}, and we are not searching over much larger models
than the oracle),
then it is easy to check \eqref{eq:ex_df_subset_reg_ass} is
satisfied, say, with
 \smash{$a_n=\sqrt{(\log|S|)/p_{s_0}}$}.  The
assumptions in \eqref{eq:ex_df_subset_reg_ass} can accomodate more
general settings, e.g., in which there is bias, or in which
\smash{$p_{\max}/p_{s_0}$} diverges, as long as these quantities scale
at appropriate rates.

Theorem \ref{thm:ex_df_subset_reg} establishes the classical oracle
inequality \eqref{eq:oracle_ineq} for the SURE-tuned subset regression
estimator, which is nothing more than the $C_p$-tuned (or 
AIC-tuned, 
as $\sigma^2$ is assumed to be known) subset  
regression estimator.  This of course is not really a new result; cf.\
classical theory on model selection in regression, as in Corollary 2.1
of \citet{li1987asymptotic}.  This 
author established a result similar to \eqref{eq:oracle_ineq} for
the $C_p$-tuned subset regression estimator, chosen over a family of
nested regression models, and showed asymptotic equivalence of the
attained loss to the oracle loss (rather than the attained and oracle
risks), in probability. 

We remark that a similar analysis to that above, where we
upper bound the excess degrees of freedom and risk, should be possible
for a general  
discrete family of linear smoothers, beyond linear regression
estimators.  This would cover, e.g., $s$-nearest neighbor
regression estimators across various choices $s=1,2,3,\ldots,|S|$.
The linear smoother setting is studied by \citet{li1987asymptotic},
and would make for another demonstration of our excess optimism
theory, but we do not pursue it.

\section{Characterizing excess degrees of freedom with (extensions of)
  Stein's formula}
\label{sec:stein}

In this section, we keep the normal assumption,  
$Y \sim F=N(\theta_0, \sigma^2 I)$ in \eqref{eq:data_model}, 
and we move beyond individual families of estimators, 
by studying the use of Stein's 
formula (and extensions thereof) for calculating excess  
degrees of freedom, in an effort to understand this quantity in 
some generality. 

\subsection{Stein's formula, for smooth estimators} 
\label{sec:stein_char}

We consider the case in which the set $S \subseteq \R$ is an
interval, i.e., in which the estimator \smash{$\htheta_s$} is
defined over a continuously-valued (rather than a discrete) tuning
parameter $s \in S$. We make the following assumption.  

\begin{assumption}
\label{ass:s_smooth}
The map \smash{$\hs : \R^n \to S$} is continuously differentiable.  
\end{assumption}

\begin{figure}[htb]
\centering
\includegraphics[width=0.55\textwidth]{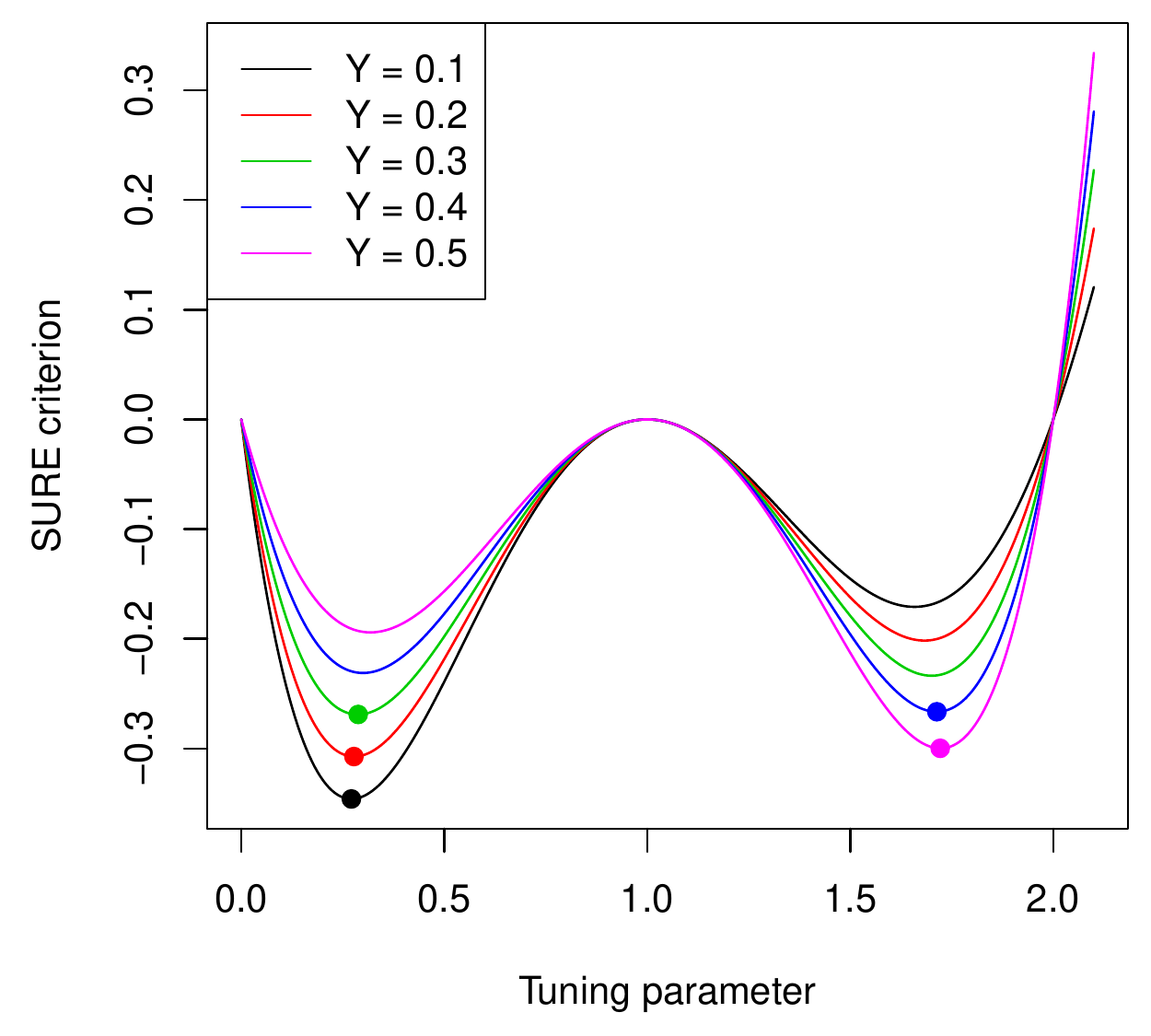}
\caption{\it\small An illustration of a discontinuous mapping 
  \smash{$\hat{s}$}. Each curve represents the SURE criterion 
    $G(Y,\cdot)$, as a function of the tuning parameter $s$, at
    nearby values of the (one-dimensional) data realization $Y$.  As
    $Y$ varies, $G(Y,\cdot)$ changes smoothly, but its minimizer 
    \smash{$\hat{s}(Y)$} jumps discontinuously, from about $0.75$ 
    at $Y=0.3$ (green curve) to $1.75$ at $Y=0.4$ (blue curve).} 
\label{fig:s_jump}
\end{figure}

It is worth noting that Assumption \ref{ass:s_smooth} seems
strong. In particular, it is not implied by the SURE criterion in
\eqref{eq:sure} being smooth in $(Y,s)$ jointly, i.e., by the map  
$G : \R^n \times S \to \R$, defined by   
\begin{equation}
\label{eq:sure_g}
G(Y,s) = \|y-\htheta_s(Y)\|_2^2 + 2\sigma^2 \hdf_s(Y),
\end{equation}
being smooth. When $G(Y,\cdot)$ is multimodal over $s \in
S$, its minimizer \smash{$\hs(Y)$} can jump discontinuously as $Y$
varies, even if $G$ itself varies smoothly.  Figure \ref{fig:s_jump}
provides an illustration of this phenomenon.
Notably, the SURE criterion for the 
family of shrinkage estimators we considered in Section
\ref{sec:shrink_means} (as well as Section \ref{sec:shrink_reg}) was  
unimodal, and Assumption \ref{ass:s_smooth} held in this
setting; however, we see no reason for this to be true in general.  
Thus, we will use Assumption \ref{ass:s_smooth} to develop a
characterization of excess degrees of freedom, shedding light on
the nature of this quantity, but should keep in mind that our
assumptions may represent a somewhat restricted setting.   

With appropriate regularity conditions placed on the family 
\smash{$\{\htheta_s : s \in S\}$}, the smoothness of \smash{$\hs$}  
guaranteed in Assumption \ref{ass:s_smooth} will imply smoothness 
of the SURE-tuned estimator \smash{$\htheta_{\hs}$}.  To state 
these regularity conditions precisely, we introduce the following
notation. 
Define the ``parent'' mapping \smash{$\widehat\Theta : \R^n \times
  S \to \R^n$} by \smash{$\htheta_s = \hTheta(\cdot, s)$} for
each $s \in S$.  Also define \smash{$h: \R^n \to \R^n \times S$} by  
\smash{$h(Y)=(Y,\hs(Y))$}. Note that, in this notation, the SURE-tuned 
estimator is given by the composition \smash{$\htheta_{\hs} = \hTheta 
  \circ h$}.  The following are our assumptions on \smash{$\hTheta$}.   

\begin{assumption}
\label{ass:theta_smooth}
The function \smash{$\hTheta : \R^n \times S \to \R^n$} is continuous 
and weakly differentiable in its first $n$ components---meaning that
it is differentiable on (Lebesgue) almost every line segment
parallel to one of the first $n$ coordinate axes.  In addition, 
\smash{$\E[\sup_{s\in S} \sum_{i=1}^n | \partial \hTheta_i(Y)/\partial
  Y_i |] < \infty$}.
\end{assumption}

The definition of weak differentiability used in
Assumption  \ref{ass:theta_smooth} is slightly stronger 
than the usual    
definition---which requires absolute continuity (instead of
differentiability) on almost every line segment parallel to the
coordinate axes.  We use the slightly stronger notion for simplicity;
together with Assumption \ref{ass:s_smooth}, it is  
easy to check that Assumption \ref{ass:theta_smooth} 
implies that the map \smash{$\htheta_{\hs} = \hTheta \circ h$} is 
continuous and weakly differentiable, and also
\smash{$\E[\sum_{i=1}^n | \partial \htheta_{\hs,i}(Y)/\partial  Y_i |]
  < \infty$}. 

Therefore we may apply Stein's formula \eqref{eq:stein_formula},
along with the chain rule, to compute the degrees of freedom of
\smash{$\htheta_{\hs}$}: 
\begin{align*}
\df(\htheta_{\hs}) &= \E\bigg(\sum_{i=1}^n
\frac{\partial (\hTheta_i \circ h)}{\partial Y_i}(Y) \bigg) \\
&= \E\bigg[\sum_{i=1}^n \bigg(\frac{\partial \hTheta_i}{\partial Y_i}    
(h(Y)) + \frac{\partial \hTheta_i}{\partial s} (h(Y)) 
\frac{\partial \hs}{\partial Y_i}(Y) \bigg) \bigg] \\
&= \E[\hdf_{\hs(Y)}(Y)] + 
\E\bigg(\sum_{i=1}^n \frac{\partial \hTheta_i}{\partial Y_i}
(Y,\hs(Y)) \frac{\partial \hs}{\partial Y_i}(Y) \bigg).
\end{align*}
Note that the Stein divergence \smash{$\hdf_s(Y)=\sum_{i=1}^n \partial
  \hTheta_i (Y,s)/\partial Y_i$} is an unbiased estimator of
\smash{$\df(\htheta_s)$}, for each $s\in S$, under Assumption 
\ref{ass:theta_smooth}.  Hence, comparing the 
last line above to the definition of excess degrees of
freedom in \eqref{eq:ex_df}, we find that
\begin{equation}
\label{eq:ex_df_explicit}
\edf(\htheta_{\hs}) = \E\bigg(\sum_{i=1}^n \frac{\partial
  \hTheta_i}{\partial s} (Y,\hs(Y)) \frac{\partial \hs}{\partial
  Y_i}(Y) \bigg). 
\end{equation}

The above expression provides an explicit characterization of excess  
degrees of freedom, and in principle, it even gives an unbiased 
estimator of excess degrees of freedom, i.e., the   
quantity inside the expectation in \eqref{eq:ex_df_explicit}.  
Note that the strategy for analyzing the families of shrinkage  
estimators in Sections \ref{sec:shrink_means} and
\ref{sec:shrink_reg} was precisely the same as that used to arrive at
\eqref{eq:ex_df_explicit} (i.e., simply employing the chain rule), and
so it is easy to check that \eqref{eq:ex_df_explicit} reproduces the   
results from these sections on excess degrees of freedom.

Unfortunately, the unbiased excess degrees of freedom estimator
suggested by
\eqref{eq:ex_df_explicit} is not always tractable. Computing  
\smash{$\partial \hTheta_i/\partial s$}, $i=1,\ldots,n$ in
\eqref{eq:ex_df_explicit} is often easy, at least when the 
estimator \smash{$\htheta_s$} (for fixed $s$) is available in
closed-form.  But computing \smash{$\partial \hs / \partial Y_i$}, 
$i=1,\ldots,n$ in \eqref{eq:ex_df_explicit} is typically much harder;  
even for simple problems, the SURE-optimal tuning parameter
\smash{$\hs$} often cannot be written in closed-form. 
Fortunately, we can use implicit differentiation to rewrite
\eqref{eq:ex_df_explicit} in more useable form. We require
the following assumption on the SURE criterion, which recall, we
denote by $G$ in \eqref{eq:sure_g}.  

\begin{assumption}
\label{ass:g_smooth}
The map $G : \R^n \times S \to \R$ is continuously differentiable.
Furthermore, for each $Y \in \R^n$, the minimizer 
\smash{$\hs(Y)$} of $G(Y,\cdot)$ is the unique point satisfying 
\begin{align}
\label{eq:g_deriv_1} 
\frac{\partial G}{\partial s} (Y,\hs(Y)) &= 0, \\
\label{eq:g_deriv_2}
\frac{\partial^2 G}{\partial s^2} (Y,\hs(Y)) &> 0.
\end{align}
\end{assumption}

As in our comment following Assumption \ref{ass:s_smooth}, we must
point out that Assumption \ref{ass:g_smooth} seems quite strong, and
as far as we can tell, in a generic problem setting there seems to be
nothing preventing $G(Y,\cdot)$ from being multimodal, which would 
violate Assumption \ref{ass:g_smooth}.  Still, we will use it to
develop insight on the nature of excess degrees of freedom. 
Differentiating \eqref{eq:g_deriv_1} with respect to $Y_i$ and using
the chain rule gives
\begin{equation*}
\frac{\partial^2 G}{\partial Y_i \partial s} (Y,\hs(Y)) + 
\frac{\partial^2 G}{\partial s^2} (Y,\hs(Y)) 
\frac{\partial \hs}{\partial Y_i}(Y) = 0,
\end{equation*}
and after rearranging,
\begin{equation*}
\frac{\partial \hs}{\partial Y_i}(Y) = 
-\bigg(\frac{\partial^2 G}{\partial s^2} (Y,\hs(Y))\bigg)^{-1} 
 \frac{\partial^2 G}{\partial Y_i \partial s} (Y,\hs(Y)).
\end{equation*}
Plugging this into \eqref{eq:ex_df_explicit}, for each $i=1,\ldots,n$, 
we have established the following result. 

\begin{theorem}
\label{thm:ex_df_implicit}
Under $Y \sim N(\theta_0, \sigma^2 I)$, and Assumptions  
\ref{ass:s_smooth}, \ref{ass:theta_smooth}, \ref{ass:g_smooth}, the
excess degrees of freedom of the SURE-tuned estimator
\smash{$\htheta_{\hs}$} is given by
\begin{equation}
\label{eq:ex_df_implicit}
\edf(\htheta_{\hs}) = -\E\bigg[
\bigg(\frac{\partial^2 G}{\partial s^2} (Y,\hs(Y))\bigg)^{-1} 
\sum_{i=1}^n \bigg( \frac{\partial
  \hTheta_i}{\partial s} (Y,\hs(Y)) 
 \frac{\partial^2 G}{\partial Y_i \partial s} (Y,\hs(Y))\bigg)
 \bigg]. 
\end{equation}
\end{theorem}

A straightforward
calculation shows that, for the classes of shrinkage estimators
in Sections \ref{sec:shrink_means} and
\ref{sec:shrink_reg}, the expression \eqref{eq:ex_df_implicit} 
matches the excess degrees of freedom results derived in these
sections. In principle, whenever Assumptions \ref{ass:s_smooth}, 
\ref{ass:theta_smooth}, \ref{ass:g_smooth} hold, 
Theorem \ref{thm:ex_df_implicit} gives an explicitly computable
unbiased estimator for excess degrees of freedom, i.e., the quantity
inside the expectation in \eqref{eq:ex_df_implicit}.  It is unclear 
to us (as we have already discussed) to what extent these assumptions  
hold in general, but we can still use \eqref{eq:ex_df_implicit} to
derive some helpful intuition on excess degrees of freedom.
Roughly speaking:
\begin{itemize}
\item
if (on average) \smash{$(\partial^2 G/\partial s^2)(Y,\hs(Y))$} is 
large, i.e., $G(Y,\cdot)$ is sharply curved around its minimum, 
i.e., SURE sharply identifies the optimal tuning parameter value
\smash{$\hs(Y)$} given $Y$, then this drives the excess degrees of  
freedom to be smaller; 
\item if (on average) \smash{$|(\partial^2 G/\partial Y_i \partial      
  s)(Y,\hs(Y))|$} is large, i.e., \smash{$|(\partial G/\partial
  s)(Y,\hs(Y))|$} varies quickly with $Y_i$, i.e., the function
whose root in \eqref{eq:g_deriv_1} determines \smash{$\hs(Y)$}
changes quickly with $Y_i$, then this drives the excess degrees of
freedom to be larger;
\item the pair of terms in the summand in \eqref{eq:ex_df_implicit}
tend to have opposite signs (their specific signs are a reflection
of the tuning parametrization associated with $s \in S$), 
which cancels out the $-1$ in front, and makes the excess degrees of  
freedom positive.    
\end{itemize}

% \begin{itemize}
% \item Perhaps if $G(Y,\cdot)$ were unimodal for each $Y$,
% then we could use some pertubation analysis type arguments?   
% \item Maybe the implicit function theorem can be applied to
%   \smash{$H = \partial G / \partial s : \R^n \times S \to
%     \R$}. Assume that $H$ is continuously differentiable, and a 
%   nonsingular Jacobian condition at $(Y,s)$.  Then the implicit
%   function theorem says that, locally, in an open neighborhood of 
%   $(Y,s)$, we can write
% \begin{equation*}
% H(Y,s) = 0 \iff s = g(Y),
% \end{equation*}
% for a unique continuously differentiable function $g$.  The
% nonsingular Jacobian condition is actually very simple for us: it
% just says that
% \begin{equation*}
% 0 \not= \frac{\partial H}{\partial s}(Y,s) = \frac{\partial^2
%   G}{\partial s^2}(Y,s).
% \end{equation*}
% So we just need the SURE criterion to be strictly curved at its
% minimum. Can't we ``piece'' together this implicit function theorem 
% and say that we have a single continuous function defining all of
% \smash{$\hs$}, i.e., \smash{$\hs(Y) = g(Y)$} for all $Y$?
% \end{itemize}
% RJT: NO, this doesn't make sense, because we would just be tracking
% a local min, not the global min, with the implicit function theorem 

\subsection{Extensions of Stein's formula, for nonsmooth  
  estimators}  
\label{sec:stein_ext_char}

When the estimator in question does not satisfy the requisite
smoothness conditions, i.e., continuity and weak
differentiability, Stein's formula 
\eqref{eq:stein_formula} is not directly applicable.
This is especially relevant to the topic
of our paper, as
the SURE-tuned estimator \smash{$\htheta_{\hs}$} can itself be
discontinuous in $Y$ even if each member of the family 
\smash{$\{\htheta_s : s \in S\}$} is continuous in $Y$ (due to
discontinuities in the SURE-optimal tuning parameter map
\smash{$\hs$}).  This will necessarily be the case for a discrete
tuning parameter set $S$, and it can also be the case for a continuous
tuning parameter set $S$, recall Figure \ref{fig:s_jump}.

Fortunately, extensions of Stein's formula have been recently
developed, to account for discontinuities of certain types.  
\citet{tibshirani2015degrees} established an extension for estimators
that are piecewise smooth. To define this notion of piecewise
smoothness precisely, we must introduce some notation. Given an
estimator \smash{$\htheta : \R^n \to \R^n$}, we 
write \smash{$\htheta_i(\,\cdot\,,Y_{-i}) : \R \to  
  \R$} for the $i$th component function \smash{$\htheta_i$} of  
\smash{$\htheta$} acting on the $i$th coordinate of the input
alone, with all other $n-1$ coordinates fixed at $Y_{-i}$.  
We also write \smash{$\cD(\htheta_i(\,\cdot\,,Y_{-i}))$} to denote the
set of  dicontinuities of the map
\smash{$\htheta_i(\,\cdot\,,Y_{-i})$}.  In this notation, the
estimator \smash{$\htheta$} is said to be {\it p-almost
  differentiable} if, for each $i=1,\ldots,n$ and (Lebesgue) almost
every $Y_{-i} \in \R^{n-1}$, the map
\smash{$\htheta_i(\,\cdot\,,Y_{-i}) : \R \to \R$} 
is absolutely continuous on each of the open intervals  
$(-\infty, \delta_1), (\delta_2, \delta_3), \ldots, (\delta_m,
\infty)$, where $\delta_1 < \delta_2 < \ldots < \delta_m$ are the 
sorted elements of
\smash{$\cD(\htheta_i(\,\cdot\,,Y_{-i}))$}, assumed to be a
finite set. For p-almost differentiable \smash{$\htheta$},
\citet{tibshirani2015degrees} proved that
\begin{equation}
\label{eq:tibs_formula}
\df(\htheta) = \E\bigg[ \sum_{i=1}^n \frac{\partial \htheta_i}
{\partial Y_i}(Y)\bigg] + \frac{1}{\sigma} \E\Bigg[
\sum_{i=1}^n \sum_{\delta \in \cD(\htheta_i(\,\cdot\,,Y_{-i}))}
\phi\bigg(\frac{\delta-\theta_{0,i}}{\sigma}\bigg)
\big[\htheta_i(\delta,Y_{-i})_+ -
\htheta_i(\delta,Y_{-i})_-\big]\Bigg],
\end{equation}
under some regularity conditions that ensure the second term on the 
right-hand side is well-defined.
%i.e., the ``extra'' term compared to the usual
%Stein's formula \eqref{eq:stein_formula}, accounting for
%discontinuities in \smash{$\htheta$}---is well-defined.  
Above, we denote one-sided limits from above and from below by  
\smash{$\htheta_i(\delta,Y_{-i})_+=\lim_{t \downarrow \delta}
  \htheta_i(t,Y_{-i})$} and
\smash{$\htheta_i(\delta,Y_{-i})_-=\lim_{t \uparrow \delta}
  \htheta_i(t,Y_{-i})$}, respectively, for the map
\smash{$\htheta_i(\cdot,Y_{-i})$}, $i=1,\ldots,n$, and we denote by
$\phi$ the univariate standard normal density. 

A difficulty with \eqref{eq:tibs_formula} is that it is often hard to
compute or characterize the extra term on the right-hand
side.  \citet{mikkelsen2016degrees}
derived an alternate extension of Stein's formula for piecewise
Lipschitz estimators.  While this setting is more restricted than 
that in \citet{tibshirani2015degrees}, the resulting
characterization is more ``global'' (instead of being based on
discontinuities along the coordinate axes), and thus it can be
more tractable in some cases. Formally, \citet{mikkelsen2016degrees}   
consider an estimator \smash{$\htheta :  \R^n \to \R^n$} with
associated regular open sets $U_j \subseteq \R^n$,  
$j=1,\ldots,J$ whose closures cover $\R^n$
(i.e., \smash{$\cup_{j=1}^J \bar{U}_j = \R^n$}), such that each 
map \smash{$\htheta^j := \htheta |_{U_j}$} (the restriction
of \smash{$\htheta$} to $U_j$) is locally Lipschitz continuous. The
authors proved that, for such an estimator \smash{$\htheta$},   
\begin{equation}
\label{eq:mh_formula}
\df(\htheta) = \E\bigg[\sum_{i=1}^n \frac{\partial \htheta_i} 
{\partial Y_i}(Y)\bigg] + \frac{1}{2} \sum_{j \not= k}  
\int_{\bar{U}_j \cap \bar{U}_k} \Big\langle 
\htheta^k(y) - \htheta^j(y),  \eta_j(y) \Big\rangle  
 \phi_{\theta_0,\sigma^2 I}(y) \, d\cH^{n-1}(y), 
\end{equation}
again under some further regularity conditions that ensure the second
term on the right-hand side is well-defined.
Above, $\eta_j(y)$ denotes the outer unit normal vector to
$\partial U_j$   
(the boundary of $U_j$) at a point $y$, $j=1,\ldots,J$,  
\smash{$\phi_{\theta_0,\sigma^2 I}$} is the density of a
normal variate with mean $\theta_0$ and covariance $\sigma^2   
I$, and $\cH^{n-1}$ denotes the $(n-1)$-dimensional Hausdorff
measure.  
 
Our interest in \eqref{eq:tibs_formula},  
\eqref{eq:mh_formula} is in applying these extensions 
to \smash{$\htheta = \htheta_{\hs}$},
the SURE-tuned estimator defined from a family \smash{$\{\htheta_s :
  s \in S\}$}. A general formula for excess degrees of freedom,
following from \eqref{eq:tibs_formula} or \eqref{eq:mh_formula}, would
be possible, but also complicated in terms of the required regularity
conditions.  Here is a high-level discussion, to reiterate
motivation for \eqref{eq:tibs_formula},
\eqref{eq:mh_formula} and outline their applications.  
We discuss the discrete
and continuous tuning parameter settings separately.
\begin{itemize}
\item 
When the tuning parameter $s$ takes discrete values
(i.e., $S$ is a discrete set), extensions such as
\eqref{eq:tibs_formula}, \eqref{eq:mh_formula} are needed to 
characterize excess degrees freedom, because the estimator
\smash{$\htheta_{\hs}$} is generically discontinuous and Stein's
original formula cannot be used. In the discrete
setting, the first term on the right-hand side of both 
\eqref{eq:tibs_formula}, \eqref{eq:mh_formula} (when
\smash{$\htheta=\htheta_{\hs}$}) is \smash{$\E[\hdf_{\hs(Y)}(Y)]$}, in
the notation of \eqref{eq:ex_df},
and thus the second term on the right-hand side of either       
\eqref{eq:tibs_formula}, \eqref{eq:mh_formula} (when
\smash{$\htheta=\htheta_{\hs}$}) gives precisely the
excess degrees of freedom. 

\item When $s$ takes continuous values (i.e., $S$ is a connected 
subset of Euclidean space),
extensions as in \eqref{eq:tibs_formula}, \eqref{eq:mh_formula} are not
strictly speaking always needed, though it seems likely to us that 
they will be needed in many cases, because the SURE-tuned estimator 
\smash{$\htheta_{\hs}$} can inherit discontinuities from the
SURE-optimal parameter map \smash{$\hs$} (recall Figure
\ref{fig:s_jump}). In the continous tuning parameter case, both the
first and second terms on the right-hand sides of
\eqref{eq:tibs_formula}, \eqref{eq:mh_formula} (when
\smash{$\htheta=\htheta_{\hs}$}) can contribute to excess degrees of
freedom; i.e., excess degrees of freedom is given by the second term 
plus any terms left over from applying the
chain-rule for differentiation in the first term.
\end{itemize}

Over the next two subsections, we
demonstrate the usefulness of the extensions in 
\eqref{eq:tibs_formula},     
\eqref{eq:mh_formula} by applying them in two specific
settings.

\subsection{Soft-thresholding estimators}
\label{sec:soft_thresh}

Consider the family of soft-thresholding estimators with component 
functions 
\begin{equation}
\label{eq:est_soft_thresh}
\htheta_{s,i}(Y) = \sign(Y_i) (|Y_i| - s)_+, \quad i=1,\ldots,n,
\quad \text{for} \; s \geq 0.
\end{equation}
In this setting, SURE in \eqref{eq:sure} is 
\begin{equation}
\label{eq:sure_soft_thresh}
\hErr_s(Y) = \sum_{i=1}^n \min\{Y_i^2,s^2\} +  
2\sigma^2 |\{i : |Y_i| \geq s \}|.
\end{equation}
Soft-thresholding estimators, like the shrinkage estimators of
Section \ref{sec:shrink_means}, have been studied extensively in
the statistical literature; some key references that study risk
properties of soft-thresholding estimators are
\citet{donoho1994ideal,donoho1995adapting,donoho1998minimax}, and
Chapters 8 and 9 of \citet{johnstone2015gaussian} give a thorough 
summary. 

The extension of Stein's formula from \citet{tibshirani2015degrees},
as given in \eqref{eq:tibs_formula}, can be used to prove that the
excess degrees of freedom of the SURE-tuned soft-thresholding
estimator is nonnegative.  The key realization is as follows:
if a component function \smash{$\htheta_{\hs,i}$} of the SURE-tuned 
soft-thresholding estimator jumps discontinuously as we move
$Y$ along the $i$th coordinate axes, then the sign of
this jump must match the direction in which $Y_i$ is moving, 
thus the latter term on the right-hand side of
\eqref{eq:tibs_formula} is always nonnegative.  The proof is given in 
the appendix. 

\begin{theorem}
\label{thm:ex_df_soft_thresh_tibs}
The SURE-tuned soft-thresholding estimator \smash{$\htheta_{\hs}$} is
p-almost differentiable. Moreover, for each $i=1,\ldots,n$, each
$Y_{-i} \in \R^{n-1}$, and each discontinuity point $\delta$ of
\smash{$\htheta_{\hs(\cdot,Y_{-i}),i}(\cdot, Y_{-i}) : \R \to \R$}, it
holds that    
\begin{equation}
\label{eq:jump_nonneg}
\big[\htheta_{\hs(\delta,Y_{-i}),i}(\delta,Y_{-i})\big]_+ -
\big[\htheta_{\hs(\delta,Y_{-i}),i}(\delta,Y_{-i})\big]_- \geq 0. 
\end{equation}
Therefore, when $Y \sim N(\theta_0,\sigma^2 I)$, we have from
\eqref{eq:tibs_formula} that 
\smash{$\edf(\htheta_{\hs}) \geq 0$} and  
\begin{equation}
\label{eq:df_soft_thresh_lb}
\df(\htheta_{\hs}) \geq \E \big|\big\{i : |Y_i| \geq \hs(Y)
\big\}\big|. 
\end{equation}
\end{theorem}

The proof of Theorem \ref{thm:ex_df_soft_thresh_tibs} provides a
precise description of the discontinuities in the SURE-tuned  
soft-thresholding estimator, which we might be able to use to give a
tight upper bound the excess degrees of freedom (second term on the
right-hand side in \eqref{eq:tibs_formula}) and upper
bound on the risk of the SURE-tuned soft-thresholding
estimator, as well. We do not pursue this.

\subsection{Subset regression estimators, revisited}
\label{sec:subset_reg_mh}

We return to the setting of Section \ref{sec:subset_reg}, i.e., we
consider the family of subset regression estimators in
\eqref{eq:est_subset_reg}, which we can abbreviate by
\smash{$\htheta_s(Y) = P_s Y$}, $s \in S$, using the notation
of the latter section. In Section \ref{sec:subset_reg_upper}, 
recall, we derived upper bounds on the excess degrees of freedom of
the SURE-tuned subset regression estimator
\smash{$\edf(\htheta_{\hs})$}.  Here we apply the extension of Stein's 
formula from 
\citet{mikkelsen2016degrees}, as stated in \eqref{eq:mh_formula},
to represent excess degrees of freedom for
SURE-tuned subset regression in an alternative and (in principle)
exact form.  The calculation of the second-term
on the right-hand side in \eqref{eq:mh_formula} for the
SURE-tuned subset regression estimator, which yields the result 
\eqref{eq:ex_df_subset_reg_mh} in the next theorem, can already
be found in \citet{mikkelsen2016degrees} (in their study of best
subset selection). A complete proof is given in the appendix
nonetheless.       

\begin{theorem}[\citealt{mikkelsen2016degrees}]
\label{thm:ex_df_subset_reg_mh}
The SURE-tuned subset regression estimator \smash{$\htheta_{\hs}$} is
piecewise Lipschitz (in fact, piecewise linear) over regular open sets 
$U_s$, $s \in S$, whose closures cover $\R^n$.  For $s,t \in S$,
the outer unit normal vector $\eta_s(y)$ to  $\partial U_s$
at a point \smash{$y \in \bar{U}_s \cap \bar{U}_t$} is given by 
\begin{equation}
\label{eq:unit_normal_subset_reg}
\eta_s(y) = \frac{(P_t-P_s)y}{\|(P_t-P_s)y\|_2}.
\end{equation}
Therefore, when $Y \sim 
N(\theta_0,\sigma^2 I)$, we have from 
\eqref{eq:mh_formula} that
\begin{equation}
\label{eq:ex_df_subset_reg_mh}
\edf(\htheta_{\hs}) = 
\frac{1}{2} \sum_{s \not= t}  
\int_{\bar{U}_s \cap \bar{U}_t} \|(P_t-P_s)y\|_2 \,
 \phi_{\theta_0,\sigma^2 I}(y) \, d\cH^{n-1}(y).
\end{equation}
\end{theorem}

An important implication of the result in
\eqref{eq:ex_df_subset_reg_mh} is the nonnegativity of 
excess degrees of freedom in SURE-tuned subset regression,
\smash{$\edf(\htheta_{\hs}) \geq 0$}, which implies that
\smash{$\df(\htheta_{\hs}) \geq \E(p_{\hat{s}(Y)})$}. 
% with no assumptions on $\theta_0$.  
% Recall that in Section \ref{sec:subset_reg_lower}, we 
% established nonnegativity of excess degrees of freedom under the  
% assumption that 
% $P_s\theta_0=\theta_0$ for all $s \in S$ (but, without assuming 
% normality). 

While the integral \eqref{eq:ex_df_subset_reg_mh} is hard to
evaluate in general, it is somewhat more tractable in the case of
nested regression models. In the present setting each $s \in S$,
recall, is identified with a subset of $\{1,\ldots,p\}$. We say
the collection $S$ is {\it nested} if for each pair $s,t \in S$, 
we have either $s \subseteq t$ or $t \subseteq s$. The next result
shows that for a nested collection of regression models, the integral
expression \eqref{eq:ex_df_subset_reg_mh} for 
excess degrees of freedom simplifies considerably, and can be upper
bounded in terms of surface areas of balls under an appropriate
Gaussian probability measure.  

Before stating the result, it helps to introduce some notation.
For a matrix $A$, we write $A_{j:k}$ as shorthand for 
\smash{$A_{\{j,j+1,\ldots,k\}}$}, i.e., the submatrix given by
extracting columns $j$ through $k$.  Likewise, for a vector $a$, we 
write $a_{j:k}$ as shorthand for $(a_j,a_{j+1},\ldots,a_k)$.  When $s$
is identified with a nonempty subset $\{1,\ldots,j\}$, we 
write $P_s,U_s,\eta_s$ as $P_j,U_j,\eta_j$ respectively, and use
\smash{$P_j^\perp$} for the orthogonal projector to $P_j$.   Lastly, we
refer to the {\it Gaussian surface measure} $\Gamma_d$, defined over
(Borel) sets $A \subseteq \R^d$ as   
\begin{equation*}
\Gamma_d(A)=\liminf_{\delta \to 0} \frac{\P(Z \in A_\delta 
\setminus A)}{\delta},
\end{equation*}
where $Z \sim N(0,I)$ denotes a $d$-dimensional standard 
normal variate, and $A_\delta=A+B_d(0,\delta)$ is the Minkowski sum of 
$A$ and the $d$-dimensional ball $B_d(0,\delta)$ centered at the
origin with radius $\delta$.  For a set $A$ with smooth boundary
$\partial A$, an equivalent definition is
\smash{$\Gamma_d(A)=\int_{\partial A}  
\phi_{0,I}(x) \, d\cH^{d-1}(x)$}, where $\phi_{0,I}$ is the density of 
$Z$, and $\cH^{d-1}$ is the $(d-1)$-dimensional Hausdorff measure.
Helpful references on Gaussian surface area include
\citet{ball1993reverse,nazarov2003maximal,klivans2008learning}.  We
now state our main result of this subsection, whose proof is given in
the appendix.  

\begin{theorem}
\label{thm:ex_df_subset_reg_nested}
Assume that $Y \sim N(\theta_0,\sigma^2 I)$, and that all models in
the collection $S$ are nested.  Then the excess degrees of freedom of
the SURE-tuned subset regression estimator \smash{$\htheta_{\hs}$} is 
\begin{equation}
\label{eq:ex_df_subset_reg_mh_nested} 
\edf(\htheta_{\hs}) = 
\sqrt{2}\sigma \sum_{s \subseteq t}  
\sqrt{p_t-p_s} \int_{\bar{U}_s \cap \bar{U}_t} 
 \phi_{\theta_0,\sigma^2 I}(y) \, d\cH^{n-1}(y).
\end{equation}
Now, without a loss of generality (otherwise, the only real 
adjustment is notational), let us identify each $s$ with a subset 
$\{1,\ldots,j\}$.  Then the excess degrees of freedom is upper bounded
by    
\begin{equation}
\label{eq:ex_df_subset_reg_mh_balls}
\edf(\htheta_{\hs}) \leq 
\sum_{d=1}^p \sqrt{2d}(d+1) 
\max_{j=1,\ldots,d} \, 
\Lambda_d \Big( B_d \big(\mu_{(j+1):(j+d)}, \sqrt{2d} \big) \Big),  
\end{equation}
where $\mu=V^T \theta_0/\sigma$, and $V \in \R^{n\times p}$ is an
orthogonal matrix with columns \smash{$v_j = P_{j-1}^\perp
  X_j/\|P_{j-1}^\perp X_j\|_2$},   
$j=1,\ldots,p$ (where we let $P_0=0$ for notational convenience). Also, 
recall that $\Lambda_d(B_d(u,r))$ denotes the $d$-dimensional Gaussian 
surface area of a ball $B_d(u,r)$ centered at $u$ with radius $r$.   
When $\theta_0=0$, the result in
\eqref{eq:ex_df_subset_reg_mh_balls} can be sharpened and 
simplified, giving 
\begin{equation}
\label{eq:ex_df_subset_reg_mh_balls_null}
\edf(\htheta_{\hs}) \leq 
\sum_{d=1}^p \sqrt{2d}\bigg(1+\frac{1}{d}\bigg) 
\Lambda_d \big( B_d (0, \sqrt{2d}) \big) < 10.
\end{equation}
\end{theorem}

Though it is established in a restricted setting, $\theta_0=0$, the
result in \eqref{eq:ex_df_subset_reg_mh_balls_null} seems quite
strong, as it shows that the excess degrees of freedom of the
SURE-tuned subset regression is bounded by the constant $10$,
and therefore its excess optimism is bounded by the constant 
$20\sigma^2$, regardless of the number of predictors $p$ in the
regression problem.  

The derivation of \eqref{eq:ex_df_subset_reg_mh_balls_null} from
\eqref{eq:ex_df_subset_reg_mh_balls} relies on two key facts: (i) the
null case, $\theta_0=0$, admits a kind of symmetry that allows
us to apply a classic result in combinatorics (the gas stations
problem) to compute the exact probability of a collection of   
chi-squared inequalities, which leads to a reduction in the factor of
$d+1$ in each summand of \eqref{eq:ex_df_subset_reg_mh_balls} to a
factor of $1+1/d$ in each summand of
\eqref{eq:ex_df_subset_reg_mh_balls_null}; and (ii) the balls in the
null case, in the summands of
\eqref{eq:ex_df_subset_reg_mh_balls_null},  are centered at the
origin, so their Gaussian surface areas can be explicitly computed as
in \citet{ball1993reverse,klivans2008learning}. 

Neither fact is true in the nonnull case, $\theta_0\not=0$, 
making it more difficult to derive a sharp upper bound on excess
degrees of freedom. We finish with a couple remarks on the nonnull
setting; more serious investigation of explicitly bounding and/or
improving \eqref{eq:ex_df_subset_reg_mh_balls} is left to future work.   
%observation about the nonnull case with
%two models, separated by one variable.

\begin{remark}[\textbf{Nonnull case: two models}] 
When our collection is composed of just two nested models that are 
separated by a single variable, i.e., 
$S=\{\{1,\ldots,p-1\},\{1,\ldots,p\}\}$, straightforward
inspection of the proof of Theorem
\ref{thm:ex_df_subset_reg_mh} reveals that
\eqref{eq:ex_df_subset_reg_mh_balls} becomes 
\smash{$\edf(\htheta_{\hs})=\sqrt{2}\Lambda_1(B_1(v_2^T  
  \theta_0/\sigma, \sqrt{2}))$} (i.e., note the equality), where   
\smash{$v_2=P_{p-1}^\perp  
  X_p/\|P_{p-1}^\perp X_p\|_2$}. The Gaussian surface measure is 
trivial to compute here (under an arbitrary mean $\theta_0$)
because it reduces to two evaluations of the Gaussian density, and
thus we see that    
\begin{equation*}
\edf(\htheta_{\hs}) =
\sqrt{2} \phi(\sqrt{2} - v_2^T \theta_0 / \sigma) + 
\sqrt{2} \phi(\sqrt{2} + v_2^T \theta_0 / \sigma),
\end{equation*}
where $\phi$ is the standard (univariate) normal density.
When $\theta_0=0$, the excess degrees of freedom is  
\smash{$2\sqrt{2} \phi(\sqrt{2}) \approx 0.415$}.  For general
$\theta_0$, it is upper bounded by 
{$\max_{u \in \R} \, \sqrt{2} \phi(\sqrt{2} - u) + 
\sqrt{2} \phi(\sqrt{2} + u) \approx 0.575$}.
\end{remark}

\begin{remark}[\textbf{Nonnull case: general bounds}]  
For an arbitrary collection $S$ of nested models and abitrary mean  
$\theta_0$, a very loose upper bound on the right-hand side in 
\eqref{eq:ex_df_subset_reg_mh_balls} is
\smash{$\sqrt{2p}p(p+1)$}, which follows
as the Gaussian surface measure of any ball is at most 1, as shown in     
\citet{klivans2008learning}.  Under restrictions on $\theta_0$, 
tighter bounds on the Gaussian surface measures of the appropriate
balls should be possible.  Furthermore, the multiplicative factor of
$d+1$ in each summand of \eqref{eq:ex_df_subset_reg_mh_balls} is also 
likely larger than it needs to be; we note that an alternate excess
degrees of freedom bound to that in
\eqref{eq:ex_df_subset_reg_mh_balls} (following from similar
arguments) is  
\begin{multline}
\label{eq:ex_df_subset_reg_mh_balls_alternate}
\edf(\htheta_{\hs}) \leq \sqrt{2} \sum_{j < k} \sqrt{k-j} \, 
\P\big(W_j(\|\mu_{1:j}\|_2^2) > 2(j-1)\big) 
\P\big(W_{p-k}(\|\mu_{(k+1):p}\|_2^2) < 2(p-k)\big) \cdot{} \\ 
\Lambda_{k-j} \Big( B_{k-j} \big(\mu_{(j+1):k}, 
\sqrt{2(k-j)} \big) \Big), 
\end{multline}
where $W_d(\lambda)$ denotes a chi-squared random variable, with $d$  
degrees of freedom and noncentrality parameter $\lambda$.  Sharp  
bounds on the noncentral chi-squared tails could deliver a useful
upper bound on the right-hand side in
\eqref{eq:ex_df_subset_reg_mh_balls_alternate}; we  
do not expect the final bound reduce to a constant (independent of
$p$) as it did in \eqref{eq:ex_df_subset_reg_mh_balls_null} in the 
null case, but it could certainly improve on the results in Section
\ref{sec:subset_reg_upper}, i.e., the bound in
\eqref{eq:ex_df_subset_reg_ub_3}, which is on the order of
\smash{$p_{\max}$} (the largest subset size in $S$).   
\end{remark}

% \begin{remark}[\textbf{General (nonnested) subsets}]  For general 
%   nonnested collections of subsets $S$, bounds on the
%   right-hand side in the excess degrees of freedom expression 
%   \eqref{eq:ex_df_subset_reg_mh} are possible, though it is not
%   clear to us how tight they will be.  For \smash{$y \in \bar{U}_s
%     \cap \bar{U}_t$}, we can expand
% \begin{equation*}
% \|(P_t-P_s)y\|_2^2 = y^T (P_t+P_s) y - y^T (P_sP_t + P_tP_s) y 
% = 2\sigma^2(p_t-p_s) + y^T(P_s^\perp P_t + P_t P_s^\perp) y, 
% \end{equation*}
% where we have used \smash{$y^T(P_t-P_s)y =
%   2\sigma^2(p_t-p_s)$}, which comes from equating the SURE  
% criterions for $s$ and $t$. The above implies the bound
% \begin{equation*}
% \|(P_t-P_s)y\|_2 \leq\sigma\sqrt{2|p_t-p_s|} + 
% \sqrt{2 \|P_s^\perp y\|_2 \|P_t y\|_2},
% \end{equation*}
% or, the weaker but perhaps more manageable bound 
% \smash{$\|(P_t-P_s)y\|_2 \leq\sigma\sqrt{2|p_t-p_s|} +
%   \sqrt{2}\|P_X y\|_2$}, where $P_X$ denotes the projection matrix 
% onto the column space of $X$). Integrating this upper bound over 
% \smash{$\bar{U}_s \cap \bar{U}_t$}
% \end{remark}

\section{Estimating excess degrees of freedom with the bootstrap} 
\label{sec:bootstrap}

We discuss bootstrap methods for estimating excess degrees of 
freedom.  As we have thus far, we assume normality, $Y \sim F = 
N(\theta_0,\sigma^2 I)$ in \eqref{eq:data_model}, but in what follows
this assumption is used mostly for convenience,and can be relaxed
(we can of course replace the normal distribution in the parameteric
bootstrap  with any known data distribution, or in general, use the
residual bootstrap). The main ideas in this section are fairly simple,
and follow naturally from standard ideas for estimating optimism using  
the bootstrap, e.g.,
\citet{breiman1992little,ye1998measuring,efron2004estimation}.  

\subsection{Parametric bootstrap procedure}

First we descibe a parametric bootstrap procedure.  We draw
\begin{equation}
\label{eq:boot_model}
Y^{*,b} \sim N(\htheta_{\hs(Y)}(Y), \sigma^2 I), \quad
b=1,\ldots,B,
\end{equation}
where $B$ is some large number of bootstrap repetitions, e.g.,
$B=1000$. Our bootstrap estimate for the excess degrees of freedom 
\smash{$\edf(\htheta_{\hs})$} is then 
\begin{equation}
\label{eq:ex_df_boot}
\widehat\edf(Y) = 
\frac{1}{B}\sum_{b=1}^B \frac{1}{\sigma^2}
\sum_{i=1}^n \htheta_{\hs(Y^{*,b}),i} (Y^{*,b}) 
(Y_i^{*,b} - \bar{Y}_i^*) - 
\frac{1}{B}\sum_{b=1}^B \hdf_{\hs(Y^{*,b})} (Y^{*,b}),
\end{equation}
where we write \smash{$\bar{Y}_i^*=(1/B)\sum_{b=1}^B Y_i^{*,b}$} for 
$i=1,\ldots,n$, and $\hdf_s$ is our estimator for the degrees of
freedom of \smash{$\htheta_s$}, unbiased for each $s \in S$.
 Note that in \eqref{eq:ex_df_boot}, for each
bootstrap draw $b=1,\ldots B$, we compute the SURE-optimal tuning
parameter value \smash{$\hs(Y^{*,b})$} for the given bootstrap data 
\smash{$Y^{*,b}$}, and we compare the sum of empirical covariances 
(first term) to the plug-in degrees of freedom estimate (second term).
We can express the definition of excess degrees of freedom in
\eqref{eq:ex_df} as    
\begin{equation}
\label{eq:ex_df_2}
\edf(\htheta_{\hs}) = 
\E \bigg( \frac{1}{\sigma^2}
\sum_{i=1}^n \htheta_{\hs(Y),i}(Y) (Y_i-\theta_{0,i})\bigg)
- \E[\hdf_{\hs(Y)}(Y)],
\end{equation}
making it clear that \eqref{eq:ex_df_boot} estimates
\eqref{eq:ex_df_2}. 
Fortuituously, the validity of the bootstrap approximation 
\eqref{eq:ex_df_boot}, as noted by \citet{efron2004estimation}, 
does not depend on the smoothness of \smash{$\htheta_{\hs}$} as a
function of $Y$.  This makes it appropriate for
estimating excess degrees of freedom, even when
\smash{$\htheta_{\hs}$} is discontinuous (e.g., due to 
discontinuities in the SURE-optimal parameter mapping \smash{$\hs$}), 
which can be difficult to handle analytically (recall 
Sections \ref{sec:stein_ext_char}, \ref{sec:soft_thresh},
\ref{sec:subset_reg_mh}).  

It should be noted, however, that typical applications of the
bootstrap for estimating optimism, as reviewed in
\citet{efron2004estimation}, consider low-dimensional problems, 
and it is not clear that \eqref{eq:ex_df_boot} will be
appropriate for high-dimensional problems.  Indeed, 
we shall see in the examples in Section \ref{sec:boot_examples} that   
the bootstrap estimate for the degrees of freedom
\smash{$\df(\htheta_{\hs})$}, 
\begin{equation}
\label{eq:df_boot}
\hdf(Y) = 
\frac{1}{B}\sum_{b=1}^B \frac{1}{\sigma^2}
\sum_{i=1}^n \htheta_{\hs(Y^{*,b}),i} (Y^{*,b}) 
(Y_i^{*,b} - \bar{Y}_i^*), 
\end{equation}
can be poor in the high-dimensional settings being considered, which
is not unexpected. But (perhaps) unexpectedly, in these same settings 
we will also see that the {\it difference} between \eqref{eq:df_boot}
and the baseline estimate \smash{$(1/B)\sum_{b=1}^B
  \hdf_{\hs(Y^{*,b})}(Y^{*,b})$}, i.e., the 
bootstrap excess degrees of freedom estimate,
\smash{$\widehat\edf(Y)$} in \eqref{eq:ex_df_boot}, can still be
reasonably accurate.     

\subsection{Alternative bootstrap procedures}

Many alternatives to the parametric bootstrap procedure of the last
subsection are possible.  These alternatives change the sampling
distribution in \eqref{eq:boot_model}, but leave the estimate in
\eqref{eq:ex_df_boot} the same. We only describe the alternatives
briefly here, and refer to the \citet{efron2004estimation} and
references therein for more details. 

In the parametric bootstrap, the mean for the sampling distribution 
in \eqref{eq:boot_model} does not have to be
\smash{$\htheta_{\hs(Y)}(Y)$}; it can be an estimate that comes
from a bigger model (i.e., from an estimator with more degrees of
freedom), believed to have low bias.  The estimate from the 
``ultimate'' bigger model, as \citet{efron2004estimation} calls it, is
$Y$ itself.  This gives rise to the alternative bootstrap sampling
procedure 
\begin{equation}
\label{eq:boot_model_2}
Y^{*,b} \sim N(Y, c\sigma^2 I), \quad b=1,\ldots,B,
\end{equation}
for some $0 < c \leq 1$, as proposed in
\citet{breiman1992little,ye1998measuring}.  The choice of
sampling distribution in \eqref{eq:boot_model_2} might work well in
low dimensions, but we found that it grossly overestimated the degrees
of freedom \smash{$\df(\htheta_{\hs})$} in the high-dimensional
problem settings considered in Section \ref{sec:boot_examples}, and
led to erratic estimates for the excess degrees of freedom  
\smash{$\edf(\htheta_{\hs})$}.  For this reason, we preferred the
choice in \eqref{eq:boot_model}, which gave more stable
estimates.\footnote{Recall, by definition, that 
  \smash{$\htheta_{\hs(Y)}(Y)$} minimizes a risk estimate 
(SURE) at $Y$, over \smash{$\htheta_s(Y)$}, $s\in S$,
so intuitively it seems reasonable to use it
in place of the mean $\theta_0$ in \eqref{eq:boot_model}.
Further, in many  
high-dimensional families of estimators, e.g., the shrinkage and
thresholding families considered in Section \ref{sec:boot_examples},
we recover the saturated estimate \smash{$\htheta_s(Y)=Y$} for one 
``extreme'' value $s$ of the tuning parameter $s$, so the mean for the
sampling distribution in \eqref{eq:boot_model} will be $Y$ if this is 
what SURE determines is best, as an estimate for $\theta_0$.}

Another alternative bootstrap sampling procedure is the residual
bootstrap, 
\begin{equation}
\label{eq:boot_model_3}
Y^{*,b} \sim \htheta_{\hs(Y)}(Y) +
\mathrm{Unif}\big(\{r_1(Y),\ldots,r_n(Y)\}\big), \quad  
b=1,\ldots,B, 
\end{equation}
where we denote by \smash{$\mathrm{Unif}(T)$} the uniform 
distribution over a set $T$, and by
\smash{$r_i(Y)=Y_i-\htheta_{\hs(Y),i}(Y)$}, $i=1,\ldots,n$ the
residuals.  The residual bootstrap \eqref{eq:boot_model_3} is
appealing because it moves us 
away from normality, and does not require knowledge of
$\sigma^2$.  Our assumption throughout this paper is that $\sigma^2$
is known---of course, under this assumption, and under a normal  
data distribution, the parametric sampler \eqref{eq:boot_model} 
outperforms the residual sampler \eqref{eq:boot_model_3},
which is we why used the parametric bootstrap in the
experiments in Section \ref{sec:boot_examples}.  
A more realistic take on the problem of estimating
optimism and excess optimism would treat $\sigma^2$ as 
unknown, and allow for nonnormal data; for such a setting, the
residual bootstrap is an important tool and deserves more
careful future study.\footnote{If estimating excess optimism is our
  goal, instead of estimating excess degrees of freedom, then we can
  craft an estimate similar to \eqref{eq:ex_df_boot} 
  that does not depend on $\sigma^2$.  Combining this
  with the residual bootstrap, we have an estimate of excess optimism
  that does not require knowledge of $\sigma^2$ in any way.}  

% \begin{equation*}
% \widehat\edf(Y) = 
% \frac{1}{\sigma^2B}\sum_{b=1}^B  \sum_{i=1}^n 
% \Big(\htheta_{\hs(Y^{*,b}),i}
% (Y^{*,b}) -\htheta_{\hs(Y),i} (Y^{*,b}) \Big) 
% (Y_i^{*,b} - \bar{Y}_i^*),
% \end{equation*}

\subsection{Simulated examples}
\label{sec:boot_examples}

We empirically evaluate the excess degrees of freedom of the
SURE-tuned shrinkage estimator and the SURE-tuned soft-thresholding
estimator, across different configurations for the data generating
distribution, and evaluate the performance of the parametric bootstrap
estimator for excess degrees of freedom.  Specifically, our simulation
setup can be described as follows.

\begin{itemize}
\item We consider 10 sample sizes $n$, log-spaced in between 10 and 
  5000.   
\item We consider 3 settings for the mean parameter $\theta_0$: 
  the null setting, where we set $\theta_0=0$; the weak sparsity
  setting, where \smash{$\theta_{0,i}=4i^{-1/2}$} for $i=1,\ldots,n$;
  and the strong sparsity setting, where \smash{$\theta_{0,i}=4$} for 
  \smash{$i=1,\ldots,\lfloor\log{n}\rfloor$} and
  \smash{$\theta_{0,i}=0$} for 
  \smash{$i=\lfloor\log{n}\rfloor+1,\ldots,n$}. 
\item For each sample size $n$ and mean $\theta_0$, we draw
  observations $Y$ from the normal data model in \eqref{eq:data_model}
  with $\sigma^2=1$, for a total of 5000 repetitions.  
  %and draw test observations $Y^*$ from the same data model
  %(independently of $Y$),   
\item For each $Y$, we compute the SURE-tuned
  estimate over the shrinkage family in \eqref{eq:est_shrink}, and the
  SURE-tuned estimate over the soft-thresholding family in  
  \eqref{eq:est_soft_thresh}. 
\item For each SURE-tuned estimator \smash{$\htheta_{\hs}$}, we
  record various estimates of degrees of freedom, excess
  degrees of freedom, and prediction error (details given below). 
% On excess degrees of freedom: 
%  \begin{itemize}
% \item Monte Carlo excess degrees of freedom estimate, defined by
%   computing \smash{$\df(\htheta_{\hs})$} using empirical covariances
%   over the 5000 repetitions, and \smash{$\E[\hdf_{\hs(Y)}(Y)]$} using
%   an empirical average over the 5000 repetitions; 
% \item unbiased excess degrees of freedom estimate (from 
%   Stein's formula, but only available for the shrinkage family),
%   defined by \smash{$2\hs(Y)/(1+\hs(Y))$};   
% \item bootstrap excess degrees of freedom estimate, defined 
%   by \eqref{eq:ex_df_boot}, using $B=500$ bootstrap draws;
% \item excess optimism-based estimate, defined by  
%   \smash{$(\|Y^*-\htheta_{\hs(Y)}(Y)\|_2^2-\hErr_{\hs(Y)}(Y))/(2\sigma^2)$}, 
%   i.e., the difference between the observed test error and the
%   observed value of the SURE criterion at its minimum, scaled
%   appropriately. 
% \end{itemize}
% On prediction error: 
% \begin{itemize}
% \item naive error estimate, defined by \smash{$\hErr_{\hs(Y)}(Y))$},
%   the minimum value of the SURE criterion;
% \item unbiased error estimate (only available for the shrinkage family),
%   defined by $2\sigma^2$ times the unbiased excess degrees of
%   freedom estimate plus the naive error estimate; 
% \item bootstrap error estimate, defined by $2\sigma^2$ times
%   the bootstrap excess degrees of freedom estimate plus the naive 
%   error estimate; 
% \item observed training error, defined by
%   \smash{$\|Y-\htheta_{\hs(Y)}(Y)\|_2^2$};
% \item observed test error, defined by
%   \smash{$\|Y^*-\htheta_{\hs(Y)}(Y)\|_2^2$}.
% \end{itemize}
\end{itemize}

The simulation results are displayed in 
Figures \ref{fig:shrink} and \ref{fig:soft_thresh};
for brevity, we only report on the null and weak sparsity 
settings for the shrinkage family, and the null and strong sparsity
settings for the soft-thresolding 
family. All degrees of freedom, excess degrees of freedom, and
prediction error estimates (except the Monte Carlo estimates)
were averaged over the 5000 repetitions; the plots all display the 
averages along with $\pm 1$ standard error bars.  

Figure \ref{fig:shrink} shows the results for the shrinkage family,  
with the first row covering the null setting, and 
the second row the weak sparsity setting. The left
column shows the excess degrees of freedom of the
SURE-tuned shrinkage estimator, for growing $n$.
Four types of estimates of excess degrees of freedom are considered:   
Monte Carlo, computed from the 5000 repetitions (drawn in     
black); the unbiased estimate from Stein's formula, i.e.,
\smash{$2\hs(Y)/(1+\hs(Y))$} (in red); the bootstrap estimate 
\eqref{eq:ex_df_boot} (in green); and the observed (scaled) excess 
optimism, i.e.,
\smash{$(\|Y^*-\htheta_{\hs(Y)}(Y)\|_2^2-\hErr_{\hs(Y)}(Y))/(2\sigma^2)$},
where $Y^*$ is an independent copy of $Y$ (in gray).  The middle
column shows similar estimates, but for degrees of
freedom; here, the naive estimate is
\smash{$\hdf_{\hs(Y)}(Y)=n/(1+\hs(Y))$}; the unbiased estimate is 
\smash{$n/(1+\hs(Y))+2\hs(Y)/(1+\hs(Y))$}; 
the naive bootstrap estimate is
the second term in \eqref{eq:ex_df_boot}; and the bootstrap estimate 
is the first term in \eqref{eq:ex_df_boot}, i.e., as given in
\eqref{eq:df_boot}.  Lastly, the right column shows the analogous
quantities, but for estimating prediction error.  The error metric is
normalized by the sample size $n$ for visualization purposes.

We can see that the unbiased estimate of excess degrees of freedom is
quite accurate (i.e., close to the Monte Carlo gold standard)
throughout. The bootstrap estimate is also accurate in the null    
setting, but somewhat less accurate in the weak sparsity setting,
particularly for large $n$.  However, comparing it to the observed
(scaled) excess  optimism---which relies on test data and thus may not
be available in practice---the bootstrap estimate still appears
reasonable accurate, and more stable.  
While all estimates of degrees of freedom are 
quite accurate in the null setting, we can see that the two bootstrap 
degrees of freedom estimates are far too small in the weak sparsity
setting. This can be attributed to the high-dimensionality of the
problem (estimating $n$ means from $n$ observations). Fortuituously,
we can see that the {\it difference} between the bootstrap and naive 
bootstrap degrees of freedom estimates, i.e., the bootstrap excess
degrees of freedom estimate, is still relatively accurate even when   
the original two are so highly inaccurate.  Lastly, the error plots
show that the correction for excess optimism is more significant
(i.e., the gap between the naive error estimate and observed test
error is larger) in the null setting than in the weak sparsity
setting.  

Figure \ref{fig:soft_thresh} shows the results for the
soft-thresholding family.  The layout of plots is the same as that for
the shrinkage family (note that the unbiased estimates of
excess degrees of freedom and of degrees of freedom are not
available for soft-thresholding). The summary of results is also
similar: we can see that the bootstrap excess degrees of
freedom estimate is fairly accurate in general, and less accurate in
the nonnull case with larger $n$.  One noteworthy difference between
Figures \ref{fig:shrink} and \ref{fig:soft_thresh}: for the
soft-thresholding family, we can see that the excess degrees of
freedom estimates appear to be growing with $n$ (perhaps even
linearly), rather than remaining upper bounded by 2, as they are for
the shrinkage family (recall also that this is clearly implied by 
the characterization in \eqref{eq:ex_df_shrink}). 

\def\wid{0.32\textwidth}
\begin{figure}[htbp]
\centering
\begin{tabular}{cccc}
\hspace{-10pt} 
\parbox[c]{\wid}{
\includegraphics[width=\wid]{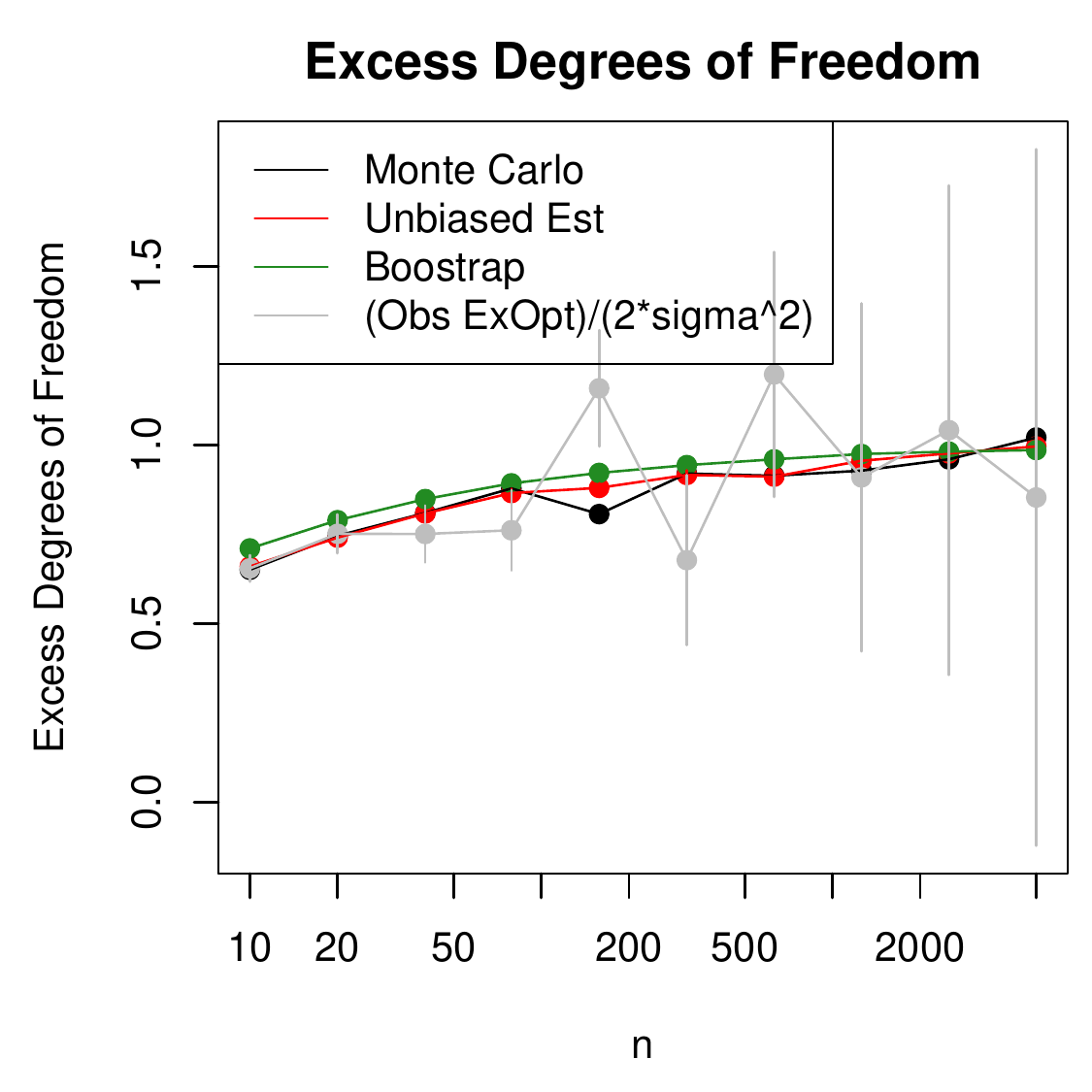}}
\hspace{-10pt} & 
\parbox[c]{\wid}{
\includegraphics[width=\wid]{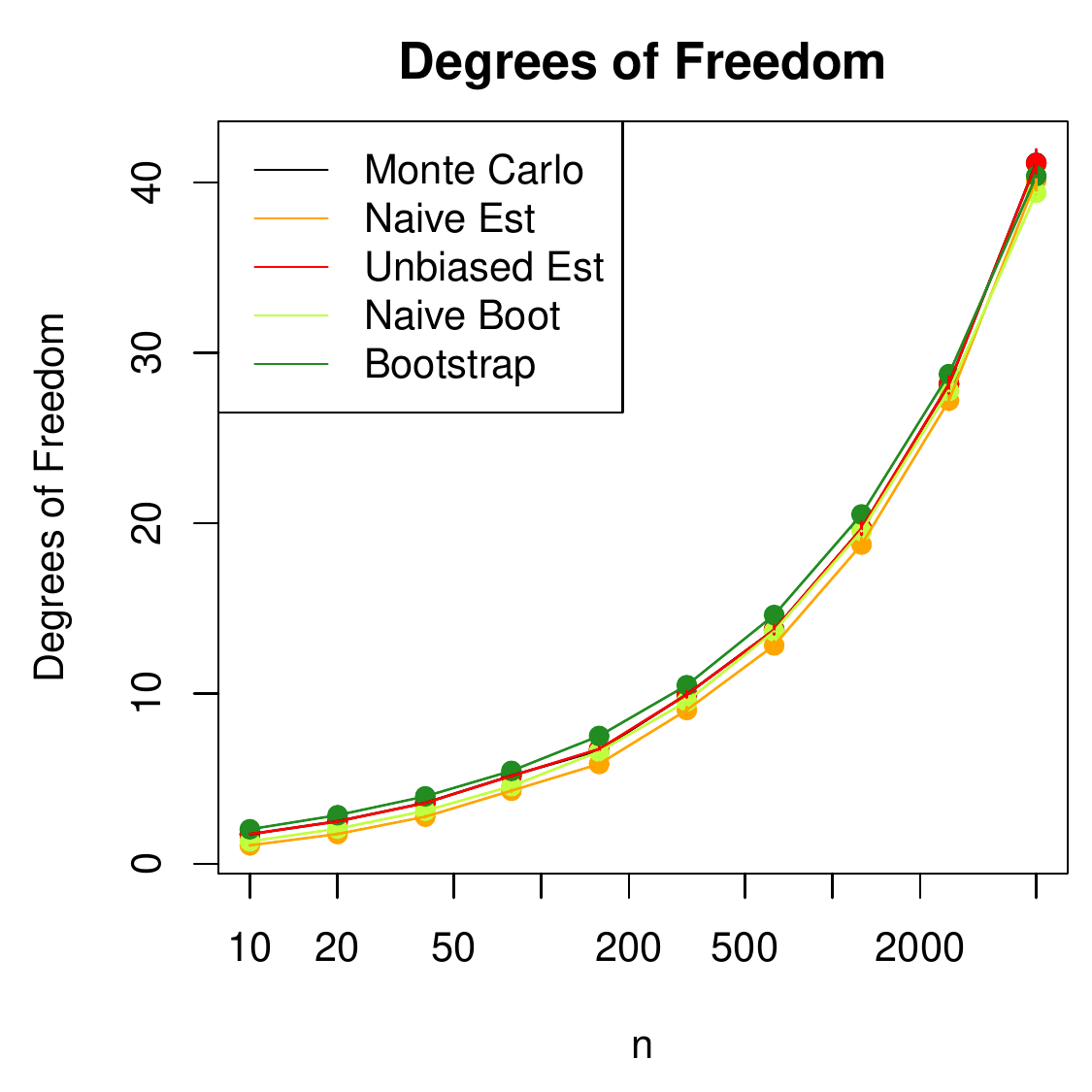}}
\hspace{-10pt} & 
\parbox[c]{\wid}{
\includegraphics[width=\wid]{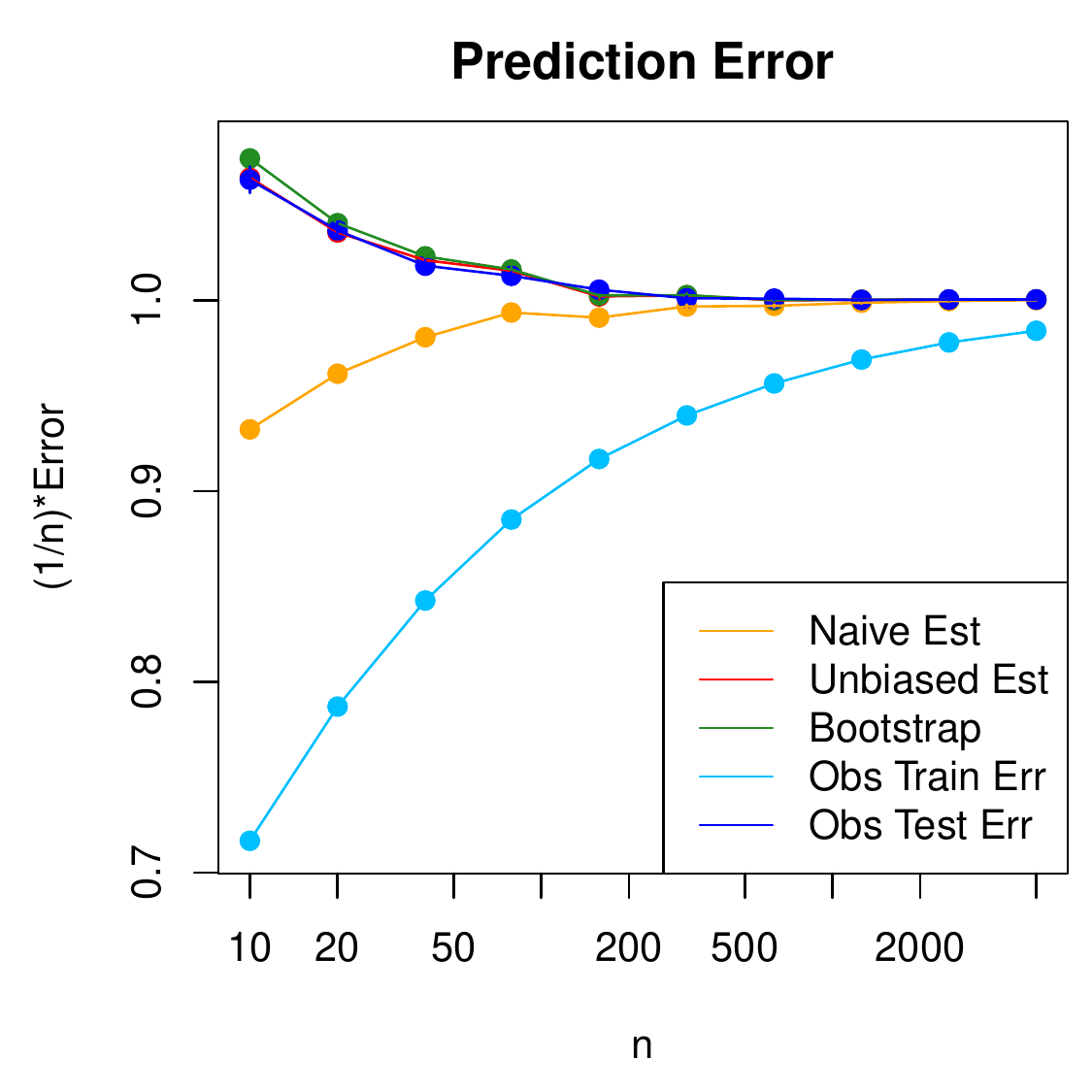}}
\hspace{-10pt} & 
\parbox[c]{0.5cm}{
\begin{turn}{-90} \small Null setting \end{turn}} \\
\parbox[c]{\wid}{
\includegraphics[width=\wid]{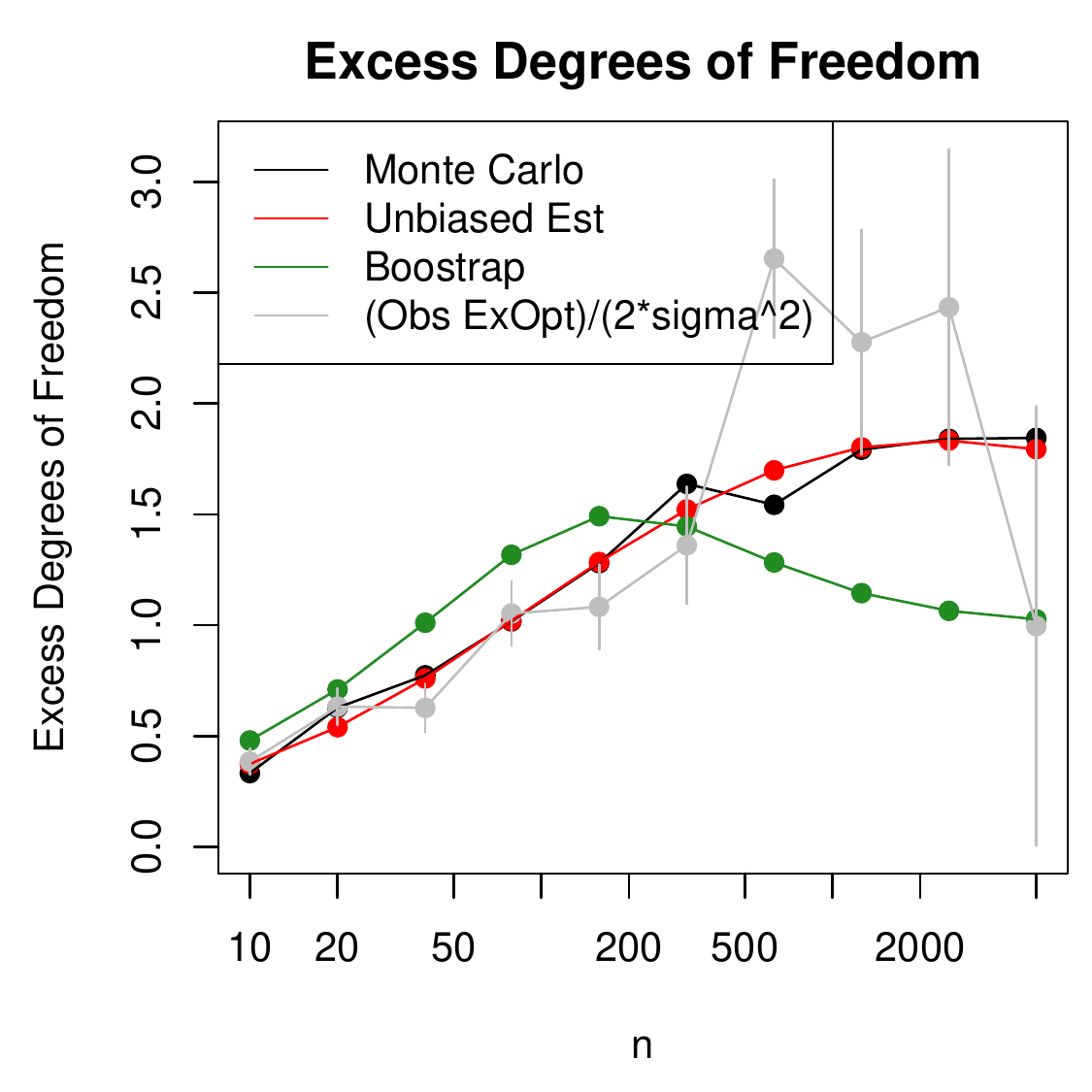}}
\hspace{-10pt} & 
\parbox[c]{\wid}{
\includegraphics[width=\wid]{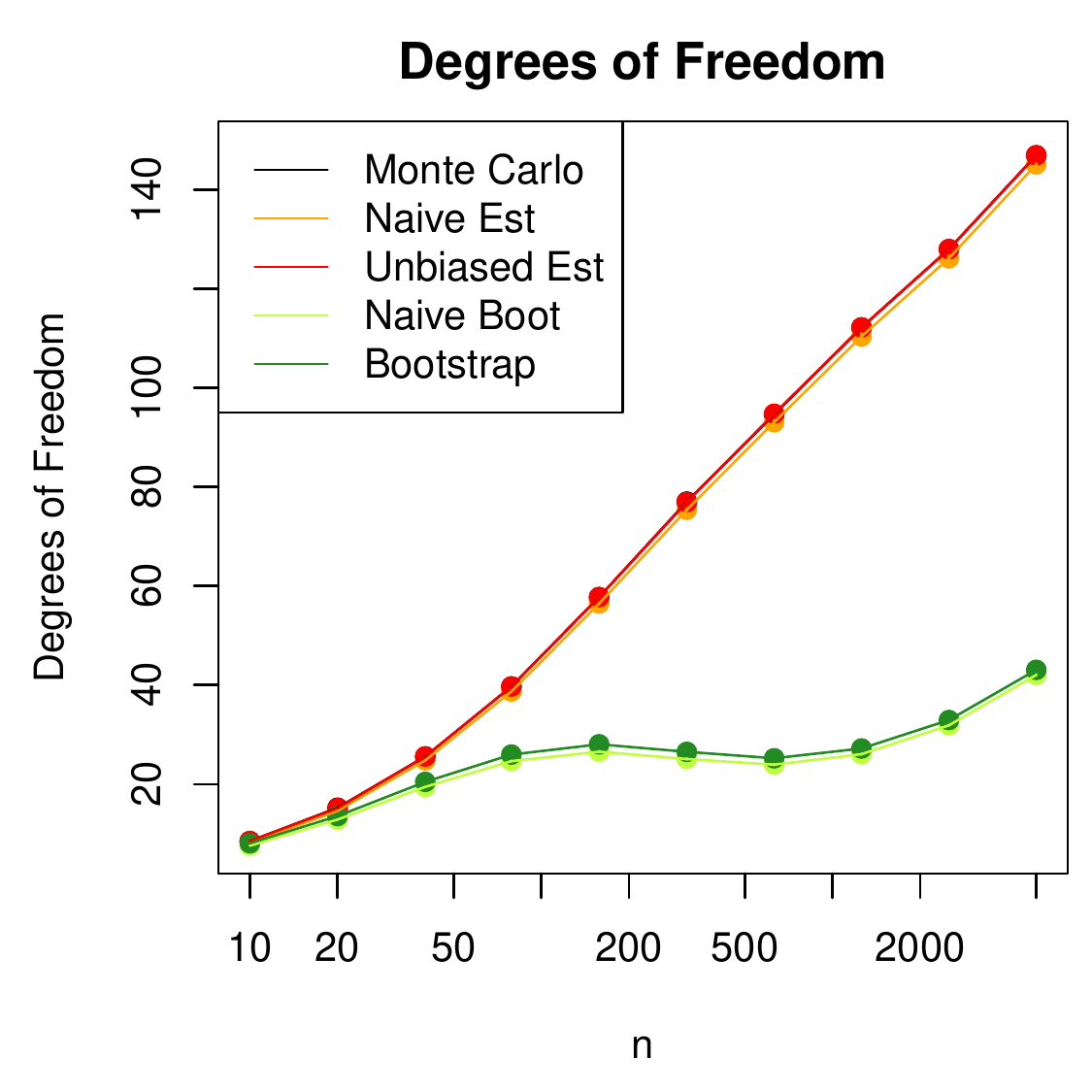}}
\hspace{-10pt} & 
\parbox[c]{\wid}{
\includegraphics[width=\wid]{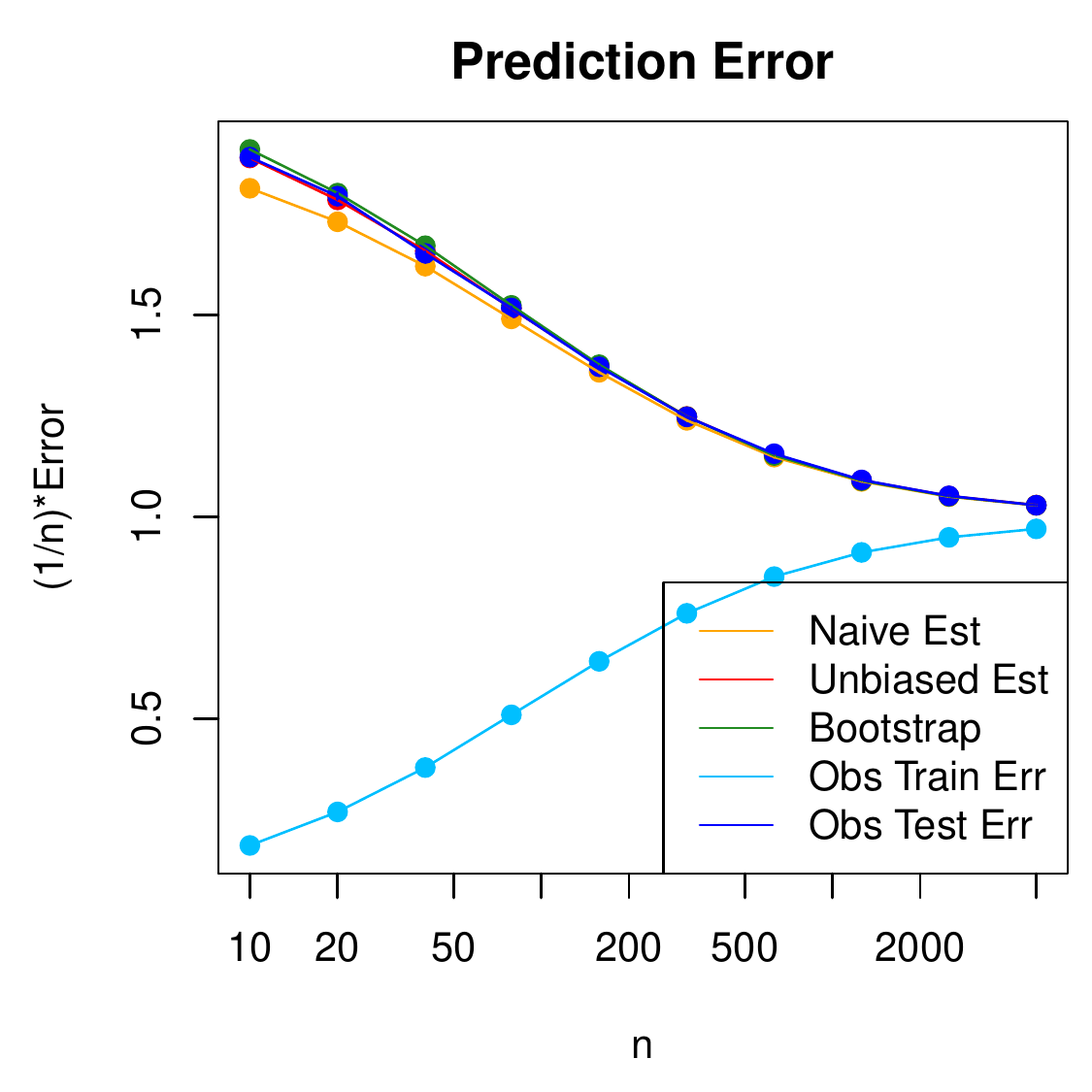}}
\hspace{-10pt} & 
\parbox[c]{0.5cm}{
\begin{turn}{-90} \small Weak sparsity setting \end{turn}}
\end{tabular}
\vspace{-10pt}
\caption{\it Simulation results for SURE-tuned shrinkage.} 
\label{fig:shrink}

\bigskip\smallskip\smallskip
\begin{tabular}{cccc}
\hspace{-10pt} 
\parbox[c]{\wid}{
\includegraphics[width=\wid]{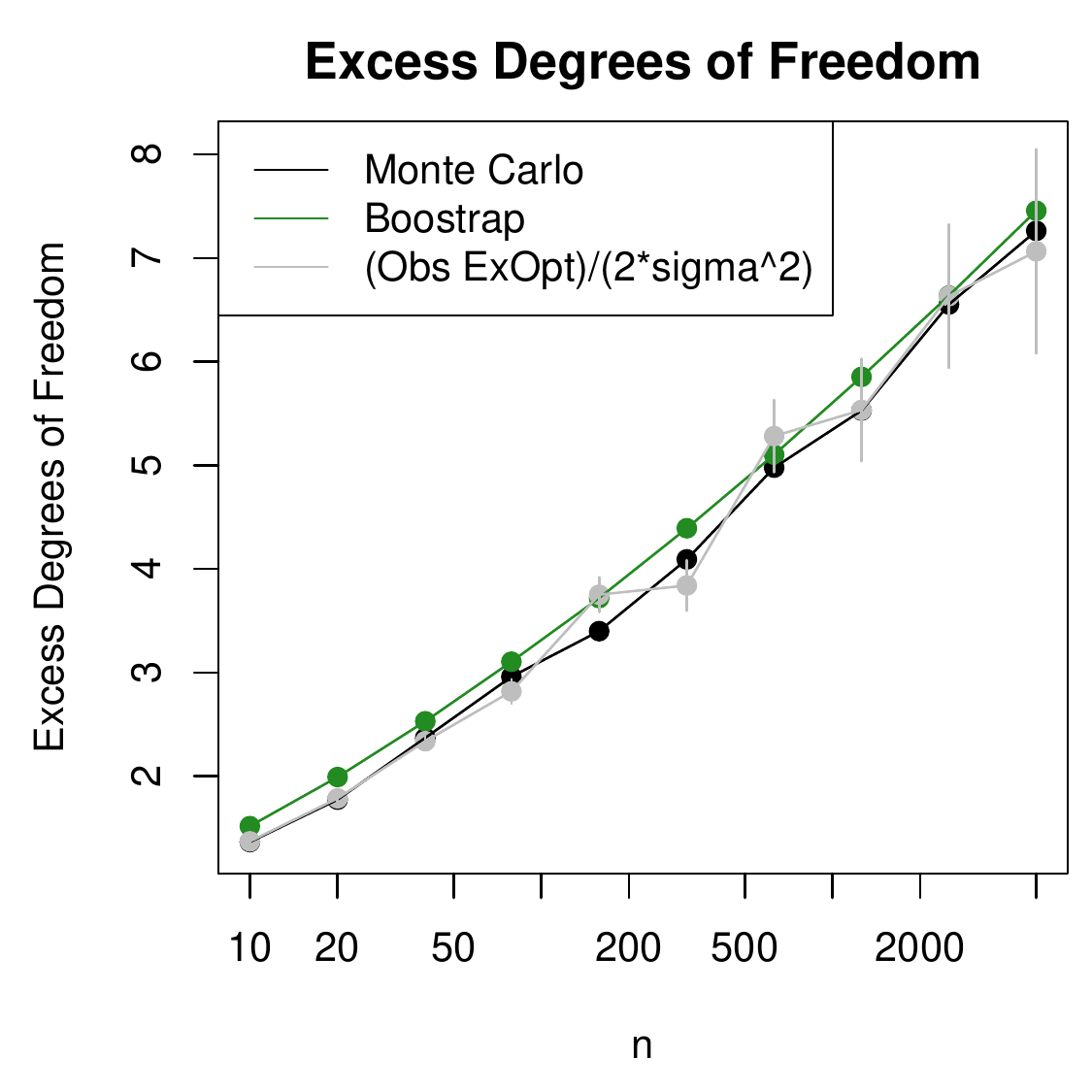}}
\hspace{-10pt} & 
\parbox[c]{\wid}{
\includegraphics[width=\wid]{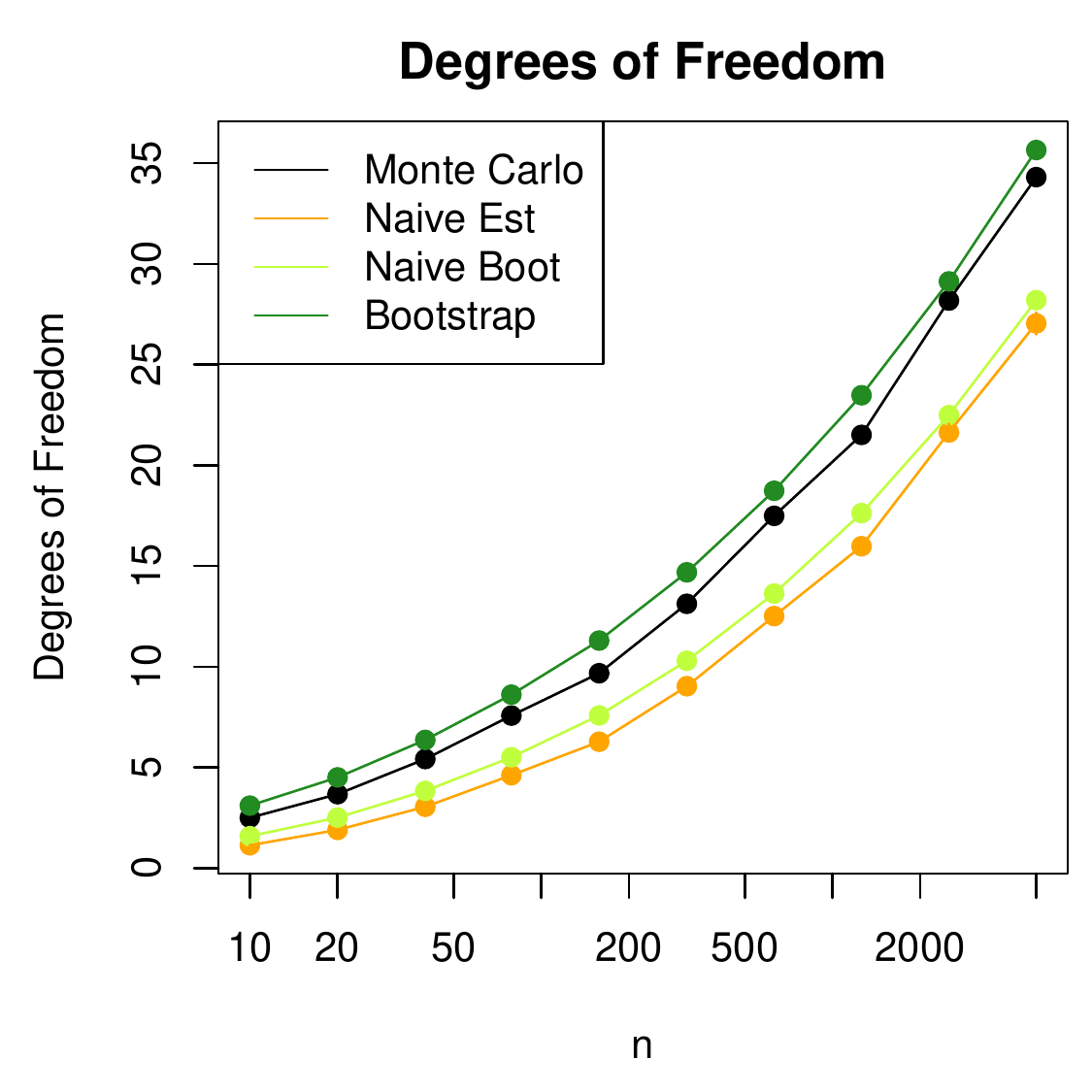}}
\hspace{-10pt} & 
\parbox[c]{\wid}{
\includegraphics[width=\wid]{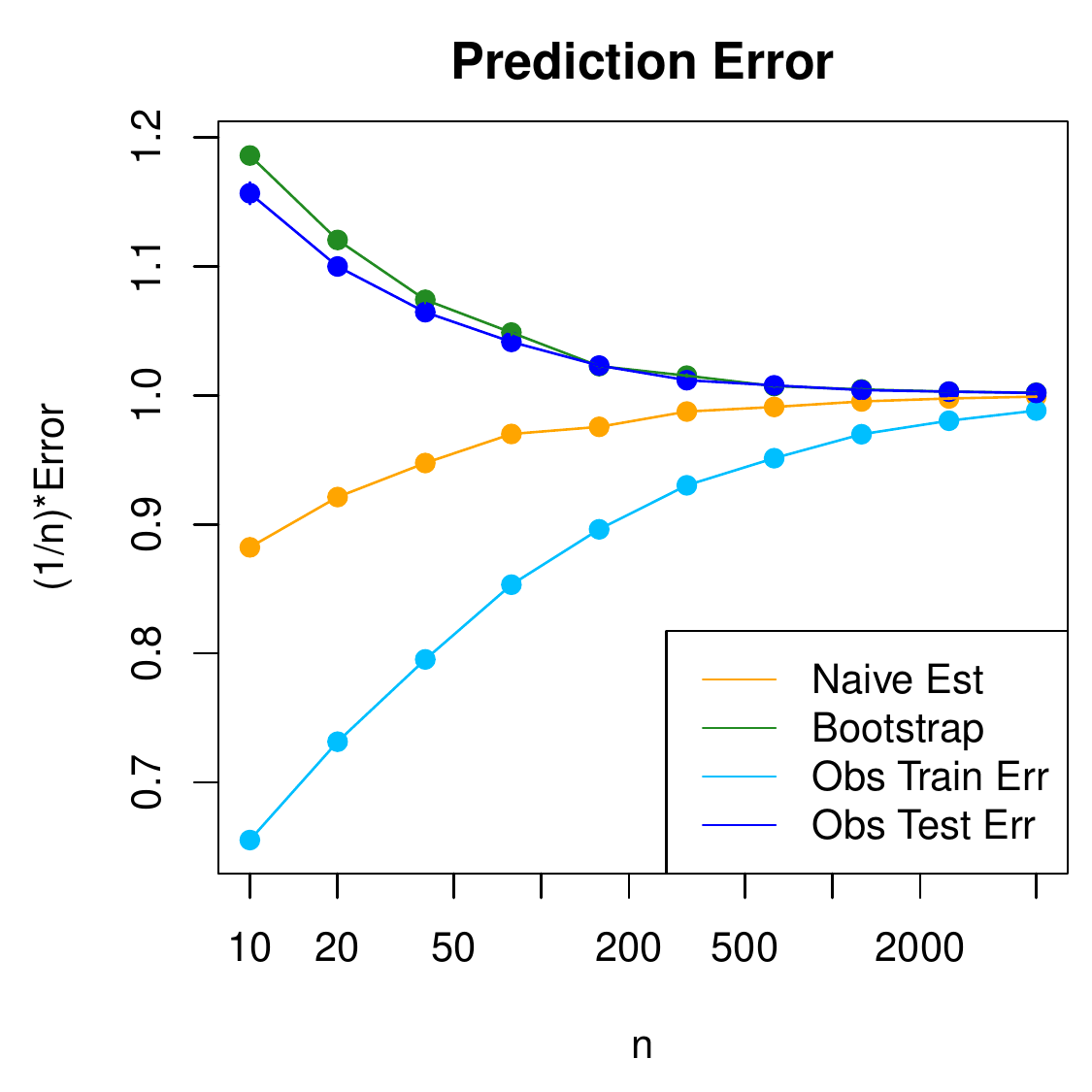}}
\hspace{-10pt} & 
\parbox[c]{0.5cm}{
\begin{turn}{-90} \small Null setting \end{turn}} \\
\parbox[c]{\wid}{
\includegraphics[width=\wid]{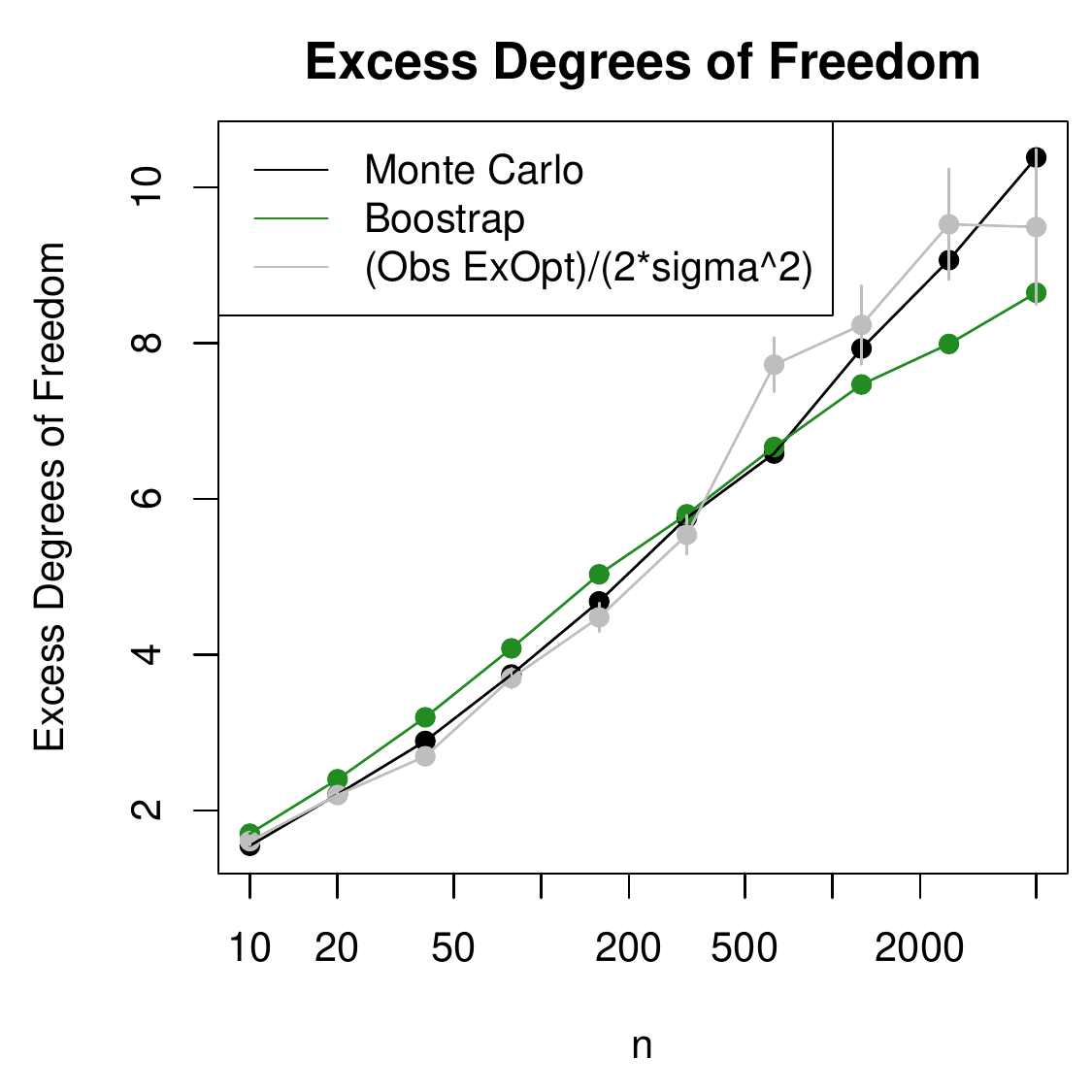}}
\hspace{-10pt} & 
\parbox[c]{\wid}{
\includegraphics[width=\wid]{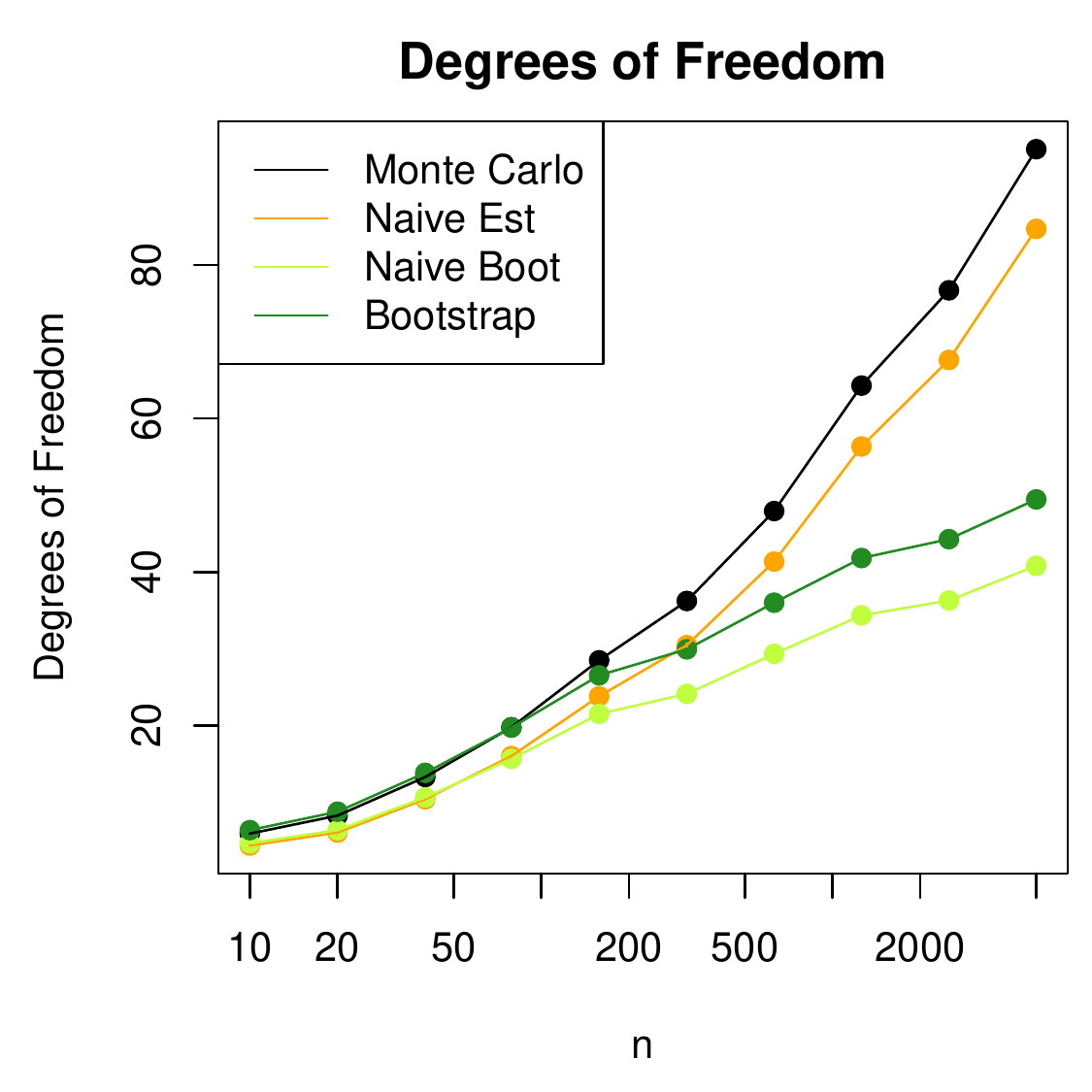}}
\hspace{-10pt} & 
\parbox[c]{\wid}{
\includegraphics[width=\wid]{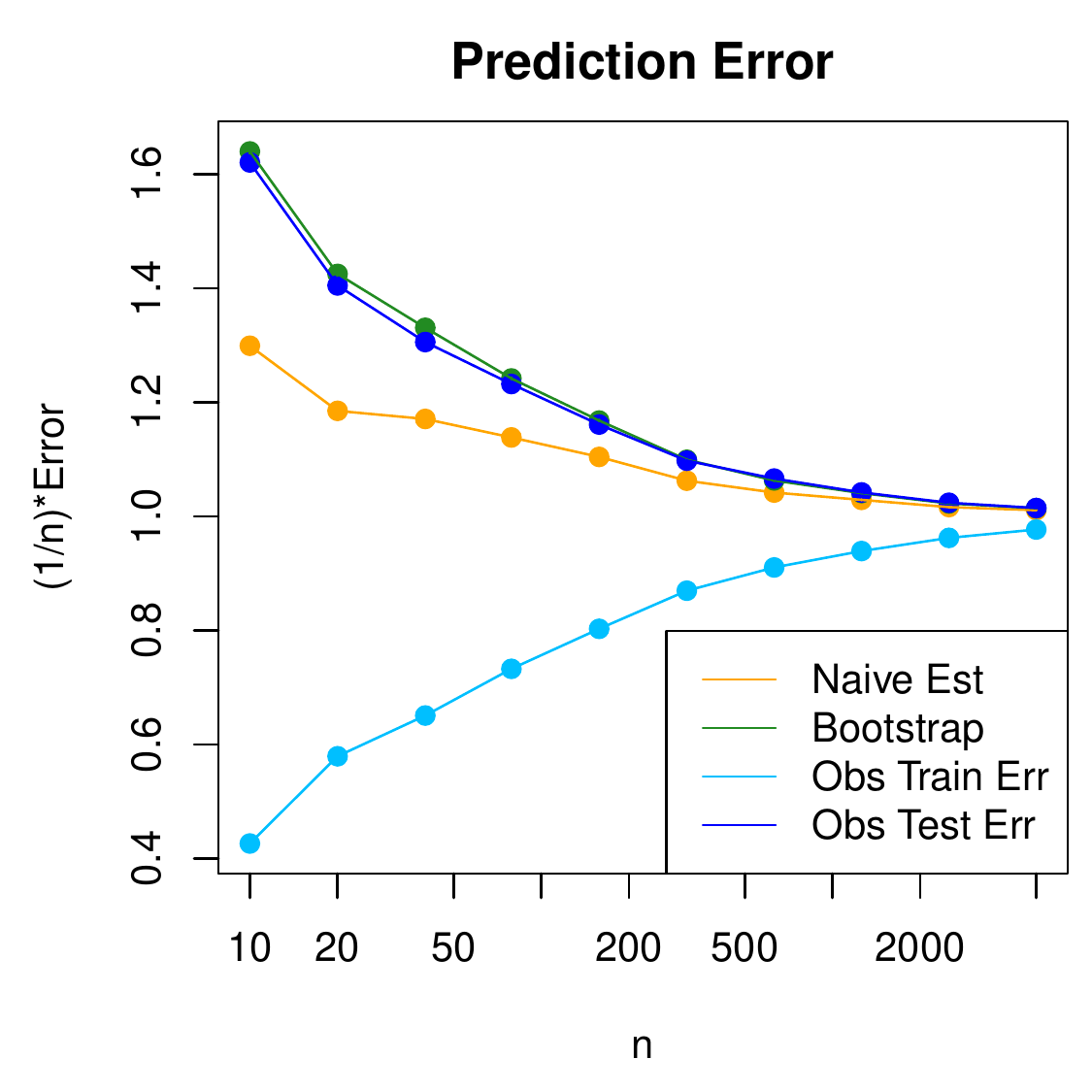}}
\hspace{-10pt} & 
\parbox[c]{0.5cm}{
\begin{turn}{-90} \small Strong sparsity setting \end{turn}}
\end{tabular}
\vspace{-10pt}
\caption{\it Simulation results for SURE-tuned soft-thresholding.} 
\label{fig:soft_thresh}
\end{figure}

\section{Discussion}
\label{sec:discussion}

We have proposed and studied a concept called excess optimism, in 
\eqref{eq:ex_opt}, which captures the added optimism of
a SURE-tuned estimator, beyond what is prescribed by SURE
itself.  By construction, an unbiased estimator of excess optimism
leads to an unbiased estimator of the prediction error of the rule
tuned by SURE.  Further motivation for the study of excess optimism
comes from its close connection to oracle estimation, as given in
Theorem \ref{thm:oracle_bd}, where we showed that the excess optimism
upper bounds the excess risk, i.e., the difference between the risk of
the SURE-tuned estimator and the risk of the oracle estimator.
Hence, if the excess optimism is shown to be sufficiently small 
next to the oracle risk, then this establishes the oracle inequality 
\eqref{eq:oracle_ineq} for the SURE-tuned estimator. 

Interestingly, excess optimism can be exactly characterized for a
family of shrinkage estimators, as studied in Section
\ref{sec:shrinkage}, where we showed that the excess optimism (and
hence the excess risk) of a class of shrinkage estimators---in both
simple normal means and regression settings---is at most $4\sigma^2$. 
For a family of subset regression estimators, such a precise
characterization is not possible, but we showed in Section
\ref{sec:subset_reg} that upper bounds on the excess optimism can be
formed that imply the oracle inequality \eqref{eq:oracle_ineq} for the
SURE-tuned (here, $C_p$-tuned) subset regression estimator. 

Characterizating excess optimism---equivalently excess degrees of
freedom, in \eqref{eq:ex_df}, which is just a constant multiple of the
former quantity---is a difficult task in general, due to 
discontinuities that may exist in the SURE-tuned estimator.  Such
discontinuities disallow the direct the use of Stein's formula
for estimating excess degrees of freedom, and in Section
\ref{sec:stein} we discussed recently developed extensions of Stein's  
formula to handle certain types of discontinuities. As an example 
application, we proved that one such extension could be used to 
bound the excess optimism of the SURE-tuned subset regression
estimator, over a family of nested subsets, by $20\sigma^2$,   
in the null case when $\theta_0=0$.
Finally, in Section \ref{sec:bootstrap}, we showed that estimation of 
excess degrees of freedom with the bootstrap is conceptually
straightforward, and appears to works reasonably well (but, it
tends to underestimate excess degrees of freedom in high-dimensional  
settings with nontrivial signal present in $\theta_0$).

We finish by noting an implication of some of our technical results
on the degrees of freedom of the best subset selection estimator, and 
discussing some extensions of our work on excess optimism to two 
related settings.   

\subsection{Implications for best subset selection} 
\label{sec:best_subset}

Our results in Sections \ref{sec:subset_reg_upper} and
\ref{sec:subset_reg_mh} have implications for the (Lagrangian version
of the) best subset selection estimator, namely, given a predictor
matrix $X \in \R^{n\times p}$, 
\begin{equation}
\label{eq:best_subset}
\hbeta^{\mathrm{subset}}_\lambda(Y) = 
\argmin_{\beta \in \R^p} \; \|Y-X\beta\|_2^2 + \lambda \|\beta\|_0, 
\end{equation}
where recall, the $\ell_0$ norm is defined by
$\|\beta\|_0=\sum_{j=1}^p 1\{\beta_j \not=0 \}$. Here $\lambda \geq
0$ is a tuning parameter. The best subset selection estimator 
in \eqref{eq:best_subset} can be seen as minimizing a SURE-like
criterion, cf.\ the SURE criterion in
\eqref{eq:sure_subset_reg}, where we define the collection $S$ to
contain all subsets of $\{1,\ldots,p\}$, and we replace the
multiplier 
$2\sigma^2$ in \eqref{eq:sure_subset_reg} with a generic parameter, 
$\lambda \geq 0$, used to weight the complexity penalty.  
Combining Lemma \ref{lem:chi_sq_max} (for the upper
bound) and Theorem \ref{thm:ex_df_subset_reg_mh} (for the lower
bound) provides the following result for best subset selection, whose 
proof is given in the appendix. 

\begin{theorem}
\label{thm:df_best_subset}
Assume that $Y \sim N(\theta_0,\sigma^2 I)$.  For any fixed value
of $\lambda \geq 0$, the degrees of freedom of the best subset
selection estimator in \eqref{eq:best_subset} satisfies 
\begin{equation}
\label{eq:df_best_subset}
\E \|\hbeta^{\mathrm{subset}}_\lambda (Y)\|_0 
\leq \df (X\hbeta^{\mathrm{subset}}_\lambda )
\leq \E \|\hbeta^{\mathrm{subset}}_\lambda (Y)\|_0 + 2.29p.
\end{equation}
\end{theorem}

In the language of \citet{tibshirani2015degrees}, the result in
\eqref{eq:df_best_subset} proves the search degrees of freedom of 
best subset selection---the difference between 
\smash{$\df(X\hbeta^{\mathrm{subset}}_\lambda)$} and 
\smash{$\E\|X\hbeta^{\mathrm{subset}}_\lambda(Y)\|_0$}---
 is nonnegative, and at most $2.29p$.  
Nonnegativity of search degrees of freedom here was conjectured by
\citet{tibshirani2015degrees} but not 
established in full generality (i.e., for general $X$); to be
fair, \citet{mikkelsen2016degrees} should be credited with
establishing this nonegativity, since, recall, the lower bound in 
\eqref{eq:df_best_subset} comes from Theorem
\ref{thm:ex_df_subset_reg_mh}, a result of these authors.
The upper bound in
\eqref{eq:df_best_subset}, as far as we can tell, is new.  Though it   
may seem loose, it implies that the degrees of 
freedom of the Lagrangian form of best subset 
selection is at most $3.29p$---in
comparison, \citet{janson2015effective} prove that best subset
selection in constrained form (for a specific configuration of the
mean particular $\theta_0$) has degrees of freedom approaching
$\infty$ as $\sigma \to 0$. This could be a reason to prefer the
Lagrangian formulation \eqref{eq:best_subset} over its constrained  
counterpart.  

\subsection{Heteroskedastic data models}
\label{sec:hetero}

Suppose now that $Y \in \R^n$, drawn from a heteroskedastic model 
\begin{equation}
\label{eq:data_model_hetero}
Y \sim F, \quad \text{where} \;
\E(Y) = \theta_0, \; \Cov(Y) = \diag(\sigma_1^2,\ldots,\sigma_n^2). 
\end{equation} 
where $\theta_0 \in \R^n$ is an unknown mean parameter,
and $\sigma^2_1,\ldots,\sigma_n^2>0$ are known variance 
parameters, now possibly distinct.  With the appropriate definitions
in place, essentially everything developed so far carries over to this 
setting.  

For an estimator \smash{$\htheta$} of $\theta_0$, define its
its prediction error, scaled by the variances, by
\begin{equation}
\label{eq:err_hetero}
\Err(\htheta) = \E\|\Sigma^{-1}(Y^*-\htheta(Y))\|_2^2 = 
\E\bigg[ \sum_{i=1}^n \frac{(Y_i^* - \htheta_i(Y))^2}{\sigma_i^2} 
\bigg], 
\end{equation}
where \smash{$\Sigma = \diag(\sigma_1^2,\ldots,\sigma_n^2)$}, and  
$Y^* \sim F$ is independent of $Y$. 
It is not hard to extend the optimism theorem and SURE, as described
in  \eqref{eq:opt_thm}, \eqref{eq:df}, \eqref{eq:opt_thm_2},
\eqref{eq:sure}, to the current heteroskedastic setting.  Similar
calculations reveal that the optimism 
\smash{$\Opt(\htheta) = \E\|\Sigma^{-1}(Y^*-\htheta(Y))\|_2^2
- \E\| \Sigma^{-1}(Y-\htheta(Y))\|_2^2$} can be expressed as 
\begin{equation}
\label{eq:opt_thm_3}
\Opt(\htheta) = 2 \tr \big(\Cov(\htheta(Y), \Sigma^{-1} Y) \big) = 
2\sum_{i=1}^n \frac{\Cov(\htheta_i(Y),Y_i)}{\sigma_i^2}. 
\end{equation}
Given an unbiased estimator \smash{$\hOpt$} of the optimism  
\smash{$\Opt(\htheta)$}, we can define an unbiased estimator 
\smash{$\hErr$} of prediction error \smash{$\Err(\htheta)$} by  
\begin{equation}
\label{eq:sure_hetero}
\hErr(Y) = \|\Sigma^{-1} (Y-\htheta(Y))\|_2^2 + \hOpt(Y), 
\end{equation}
which we will still refer to as SURE.  Assuming that \smash{$\htheta$}
is continuous and weakly differentiable, it is implied by Lemma 2 in 
\citet{stein1981estimation} that 
\begin{equation}
\label{eq:stein_formula_hetero}
\Opt(\htheta) = 2\E\bigg[ \sum_{i=1}^n 
\frac{\partial \htheta_i(Y)}{\partial Y_i} \bigg],
\end{equation}
i.e., \smash{$\hOpt(Y)= 2\sum_{i=1}^n 
\partial \htheta_i(Y)/\partial Y_i$} is an unbiased estimate of  
optimism.  

Sticking with our usual notation \smash{$\htheta_s, \hErr_s, \hOpt_s$}
to emphasize dependence on a tuning parameter $s \in S$, 
% When we consider an estimator
% \smash{$\htheta_s$} that depends on a tuning parameter $s \in S$,  
% and we use \smash{$\hErr_s$} to denote SURE in
% \eqref{eq:sure_hetero} for \smash{$\htheta_s$} and \smash{$\hOpt_s$}
% to denote an unbiased estimator of the optimism of
% \smash{$\htheta_s$}.    
we can define excess optimism in 
the current heteroskedastic setting just as before,
in \eqref{eq:ex_opt}.  An important note is that excess optimism still
upper bounds the excess prediction error, i.e., the result in
\eqref{eq:oracle_bd} of Theorem \ref{thm:oracle_bd} still holds.

We briefly sketch an example of an estimator that could be seen as an
extension of the simple shrinkage estimator in Section
\ref{sec:shrink_means} to the heteroskedastic setting.  In particular,
assuming normality in the model in \eqref{eq:data_model_hetero}, i.e.,  
\smash{$Y \sim F = N(\theta_0, \Sigma)$}, with 
\smash{$\Sigma=\diag(\sigma_1^2,\ldots,\sigma_n^2)$}, 
consider 
\begin{equation}
\label{eq:est_shrink_hetero}
\htheta_{s,i}(Y) = \frac{Y_i}{1 + \sigma_i^2 s}, \quad
i=1,\ldots,n, \quad \text{for} \; s \geq 0,
\end{equation}
For each $s \geq 0$, note that \smash{$\htheta_s$}  
is the Bayes estimator under the prior 
$\theta_0 \sim N(0,s^{-1}I)$.  The family in
\eqref{eq:est_shrink_hetero} of heteroskedastic
(nonuniform) shrinkage estimators is studied in \citet{xie2012sure}.   
It is easy to verify that SURE in \eqref{eq:sure_hetero} for this
family is  
\begin{equation}
\label{eq:sure_hetero_shrink} 
\hErr_s(Y) = \sum_{i=1}^n
\bigg(  Y_i^2 \frac{\sigma_i^2 s^2}{(1+ \sigma_i^2 s)^2}
+ \frac{2} {1+\sigma_i^2 s}\bigg).  
\end{equation}
(\citet{xie2012sure} arrive at a slightly different criterion 
because they study unscaled prediction error rather than
the scaled version we consider in \eqref{eq:err_hetero}.)

Unfortunately, the exact minimizer \smash{$\hs(Y)$} of the above
criterion cannot be written in closed form, as it could (recall
Lemma \ref{lem:sure_min_shrink}) in Section \ref{sec:shrink_means}. 
But, assuming that Assumptions \ref{ass:s_smooth} and
\ref{ass:g_smooth} of Section \ref{sec:stein_char} hold 
(we can directly check Assumption \ref{ass:theta_smooth} 
for the family of estimators in
\eqref{eq:est_shrink_hetero}), implicit differentiation can be used to
characterize the excess degrees of freedom of the 
SURE-tuned heteroskedastic shrinkage estimator
\smash{$\htheta_{\hs}$}. As before, this leads to 
\begin{equation}
\label{eq:ex_opt_implicit_hetero}
\ExOpt(\htheta_{\hs}) = -2\E\bigg[
\bigg(\frac{\partial^2 G}{\partial s^2} (Y,\hs(Y))\bigg)^{-1} 
\sum_{i=1}^n \bigg( \frac{\partial
  \hTheta_i}{\partial s} (Y,\hs(Y)) 
 \frac{\partial^2 G}{\partial Y_i \partial s} (Y,\hs(Y))\bigg)
 \bigg],
\end{equation}
where \smash{$\hTheta$} denotes the family in
\eqref{eq:est_shrink_hetero} as a function of $Y$ and $s$, and 
$G$ denotes the SURE criterion as a function of $Y$ and $s$. 
The above generalizes the result in \eqref{eq:ex_df_implicit} of
Theorem \ref{thm:ex_df_implicit} for the homoskedastic setting.
Computing \eqref{eq:ex_opt_implicit_hetero} for the heteroskedastic
shrinkage family in \eqref{eq:est_shrink_hetero} gives 
\begin{equation}
\label{eq:ex_opt_shrink_hetero}
\ExOpt(\htheta_{\hs}) = 
\E\left(\frac{\displaystyle
\sum_{i=1}^n 
\frac{4Y_i^2 \sigma_i^4 \hs(Y)}{(1+\sigma_i^2 \hs(Y))^5}} 
{\displaystyle
\sum_{i=1}^n \bigg[
\frac{\sigma_i^2}{(1+\sigma_i^2 \hs(Y))^2}  
\bigg( Y_i^2 - 
\frac{4 Y_i^2 \sigma_i^2 \hs(Y)}
{1+\sigma_i^2 \hs(Y)} +
\frac{3 Y_i^2 \sigma_i^4 \hs(Y)^2}
{(1+\sigma_i^2 \hs(Y))^2} +
\frac{2 \sigma_i^2}
{1+\sigma_i^2 \hs(Y)} 
\bigg) \bigg]}\right).
\end{equation}
We reiterate that the above hinges on Assumptions \ref{ass:s_smooth} 
and \ref{ass:g_smooth}.  It is not clear to us in what generality
these assumptions hold for the heteroskedastic shrinkage family
\eqref{eq:est_shrink_hetero} 
(clearly, when $\sigma_1^2=\ldots=\sigma_n^2$, these assumptions hold, 
since in this case the family reduces to the homoskedastic family in
\eqref{eq:est_shrink}, and then these assumptions can be easily
verified, 
as discussed previously).  Without Assumptions
\ref{ass:s_smooth} and \ref{ass:g_smooth}, there  
would need to be an additional term added to the right-hand side in
\eqref{eq:ex_opt_shrink_hetero} that accounts for discontinuities in
the SURE-tuned heteroskedastic shrinkage estimator
\smash{$\htheta_{\hs}$} (e.g., as specified in the second term on the
right-hand side in \eqref{eq:tibs_formula}).  Deriviation details for
\eqref{eq:ex_opt_shrink_hetero} are given in the appendix.  It can be
checked that \eqref{eq:ex_opt_shrink_hetero} is indeed equivalent to
\eqref{eq:ex_df_shrink} when all the variances are equal to
$\sigma^2$. 

Interestingly, as we now show, we can view ridge regression through
the lens of a heteroskedastic data setup as in
\eqref{eq:data_model_hetero}.
Given a predictor matrix $X \in \R^{n\times p}$, it is
well-known that the solution to the ridge regression problem 
\eqref{eq:ridge} is 
\smash{$\hbeta^{\mathrm{ridge}}_s(Y) = (X^T X + s I)^{-1} X^T
  Y$}, for any $s \geq 0$. Denote the singular value decomposition of
$X$ by $X=UDV^T$. If $Y$ follows the usual homoskedastic
distribution in \eqref{eq:data_model} with $F = N(\theta_0, I)$ (here
we have set $\sigma^2=1$ for simplicity, and without a loss of
generality), then a rotation and diagonal scaling gives 
\begin{equation*}
W \sim N(\alpha_0, D^{-2}),
\end{equation*}
where $W=D^{-1}U^T Y$, and $\alpha_0=D^{-1}U^T \theta_0$.  
Further, we can simply deal with (excess) optimism in this new  
coordinate system, since for any estimator \smash{$X\hbeta$} of
$\theta_0$, we have 
\begin{equation*}
\Opt(X\hbeta) = 2\tr\big(\Cov (X\hbeta(Y),Y)\big) =
2 \tr\big(\Cov(\halpha(W), D^2 W)\big) = 
\Opt(\halpha),
\end{equation*}
where \smash{$\halpha(W)=V^T \hbeta(Y)$}.
Thus, let us define \smash{$\halpha_s(W) = 
V^T \hbeta^{\mathrm{ridge}}_s(Y)$}, for $s\geq 0$. 
It is easy to see that
\smash{$\halpha_s(W) = (D^2+sI)^{-1} D^2 W$}, for $s \geq 0$, i.e.,    
\begin{equation}
\label{eq:est_shrink_hetero_2}
\halpha_{s,i}(W) = \frac{W_i}{1+d_i^{-2} s},  
\quad i=1,\ldots,r, \quad\text{for} \; s \geq 0,
\end{equation}
where $r$ is the rank of $X$, and $d_1 \geq \ldots \geq d_r > 0$ are
the diagonal elements of $D$.  Hence 
the setup in \eqref{eq:est_shrink_hetero_2} is exactly that in
\eqref{eq:est_shrink_hetero}. The result in
\eqref{eq:ex_opt_shrink_hetero} shows, under Assumptions
\ref{ass:s_smooth} and \ref{ass:g_smooth} (Assumption
\ref{ass:theta_smooth} can be checked directly), that the excess
optimism of the SURE-tuned ridge regression estimator is 
\begin{multline}
\label{eq:ex_opt_ridge}
\ExOpt(X\hbeta^{\mathrm{ridge}}_{\hs}) = \\ 
\E\left(\frac{\displaystyle
\sum_{i=1}^r 
\frac{4(u_i^T Y)^2 d_i^{-6} \hs(Y)}{(1+d_i^{-2} \hs(Y))^5}}    
{\displaystyle
\sum_{i=1}^r \bigg[
\frac{d_i^{-4}}{(1+d_i^{-2} \hs(Y))^2}  
\bigg( (u_i^T Y)^2 -  
\frac{4 (u_i^T Y)^2 d_i^{-2} \hs(Y)}
{1+d_i^{-2} \hs(Y)} +
\frac{3(u_i^T Y)^2 d_i^{-4} \hs(Y)^2} 
{(1+d_i^{-2} \hs(Y))^2} +
\frac{2 d_i^{-2}}
{1+d_i^{-2} \hs(Y)} 
\bigg) \bigg]}\right),
\end{multline}
where $u_1,\ldots,u_r \in \R^n$ are the columns of $U$. 
As before, we must stress that it is not at all clear to us in what
situations Assumption \ref{ass:s_smooth} and \ref{ass:g_smooth} will
hold for ridge regression, and so \eqref{eq:ex_opt_ridge} should be
seen as only one ``piece of the puzzle'' for ridge regression, as
it may be missing important terms (that account
for discontinuities in the SURE-tuned ridge estimator).  It could
still be interesting to work with the right-hand side in
\eqref{eq:ex_opt_ridge}, and derive bounds on this quantity
under various models for the decay of singular values of $X$.  This is
left to future work, along with a study of the
discontinuities of \smash{$X\hbeta^{\mathrm{ridge}}_{\hs}$}, and the  
resulting adjustments that need to be made to
\eqref{eq:ex_opt_ridge}. 
% TODO: Dijskstra work, should we mention that it gets a different
% answer, we think, because it makes a calculation error? 

\subsection{Efron's $Q$ measures}
\label{sec:efron_q}

We stick with the data model in \eqref{eq:data_model_hetero}.
Instead of the normality-inspired squared loss in
\eqref{eq:err_hetero}, let us consider a sequence of loss functions
$Q_i$, $i=1,\ldots,n$, and define the error metric 
\begin{equation}
\label{eq:err_q}
\Err(\htheta) = \E\bigg[\sum_{i=1}^n Q_i(\htheta_i(Y), Y_i^*) 
\bigg],  
\end{equation}
where $Y^* \sim F$, independent of $Y$.  We assume that, for each
$i=1,\ldots,n$, each $Q_i$ is the {\it tangency function} of a
differentiable, concave function $q_i$, i.e.,
\begin{equation*}
Q_i(u, v) = q_i(u) - q_i(v) + q_i'(u)(v-u),
\end{equation*}
where $q_i'$ denotes the derivative of $q_i$. We will refer to
$Q_i$ as one of {\it Efron's $Q$ measures}, in honor of
\citet{efron1986biased,efron2004estimation}, who developed an optimism
theorem in the current general setting.  (Our setup here is only a
very slight generalization of Efron's, in which we allow for different
loss functions $Q_i$, for $i=1,\ldots,n$.)

Some examples, as covered in \citet{efron1986biased}: when  
$q_i(u)=u(1-u)/\sigma_i^2$, we get the squared loss $Q_i(u,v) =
(v-u)^2/\sigma_i^2$, and \eqref{eq:err_q} 
recovers \eqref{eq:err_hetero}; when $q_i(u)=\min\{u,1-u\}$, we get
the 0-1 loss for $Q_i$; when $q_i(u) = -2(u\log{u}-(1-u)\log(1-u))$,
we get the binomial deviance for $Q_i$; in general, for any
exponential family distribution, there is a natural concave function
$q_i$ that can be defined that makes $Q_i$ the deviance. 

Now let us define
\begin{equation*}
\heta_i(Y) = -q_i'(\htheta_i(Y))/2, \quad i=1,\ldots,n. 
\end{equation*}
\citet{efron1986biased} derived the following beautiful  
generalization of the optimism theorem (with further discussion 
in \citet{efron2004estimation}): the optimism
\smash{$\Opt(\htheta) =
\E[\sum_{i=1}^n Q_i(\htheta_i(Y),Y^*_i)]-
\E[\sum_{i=1}^n Q_i(\htheta_i(Y),Y_i)]$}
can be alternatively expressed as
\begin{equation}
\label{eq:opt_thm_q}
\Opt(\htheta) = \sum_{i=1}^n \Cov(\heta_i(Y), Y_i). 
\end{equation}
Hence, given an estimator \smash{$\hOpt$} of optimism, we can define
an estimator \smash{$\hErr$} of the error \smash{$\Err(\htheta)$} by   
\begin{equation}
\label{eq:sure_q}
\hErr(Y) = \sum_{i=1}^n Q_i(\htheta_i(Y),Y_i) + \hOpt(Y),
\end{equation}
and \smash{$\hErr$} will be unbiased %(for \smash{$\Err(\htheta)$}) 
provided that \smash{$\hOpt$} is. %(for \smash{$\Opt(\htheta)$}).   

% When we consider an estimator
% \smash{$\htheta_s$} depending on a tuning parameter $s \in S$,   
% and as usual, write \smash{$\hErr_s$} and \smash{$\hOpt_s$} for the  
% error and optimism estimators, respectively, associated with 
% \smash{$\htheta_s$}, 
Keeping the usual notation \smash{$\htheta_s, \hErr_s, \hOpt_s$} to
mark the dependence on a tuning parameter $s \in S$,
we can define excess optimism for the current setting precisely as
before, in \eqref{eq:ex_opt}.  Assuming that 
\smash{$\hOpt_s$} is unbiased, an important realization is that the
result in \eqref{eq:oracle_bd} of Theorem \ref{thm:oracle_bd} 
holds as written, i.e., the excess optimism still upper bounds the
excess prediction error. 

In principle, this an exciting extension to pursue.  One problem 
is that it is difficult to form an unbiased estimator of the optimism
in \eqref{eq:opt_thm_q}, and therefore difficult to form an  
unbiased estimator of prediction error, as defined in
\eqref{eq:sure_q}.  By this, we mean specifically that it is difficult
to analytically construct an unbiased estimator of optimism
(the bootstrap can be used to give an approximately unbiased
estimator of optimism, as in Section
\ref{sec:bootstrap}).  Under appropriate smoothness conditions on
\smash{$\htheta$}, \citet{efron1986biased} proposed to use the
divergence 
\begin{equation}
\label{eq:efron_div}
\hOpt(Y) = 2\sum_{i=1}^n \frac{\partial \htheta_i}{\partial
  Y_i}(\htheta(Y))
\end{equation}
to estimate optimism.  Note the point of evaluation in 
\eqref{eq:efron_div}:  it is \smash{$\htheta(Y)$}, not
$Y$, as in the usual Stein divergence \eqref{eq:stein_div}.  
\citet{efron1986biased} showed the divergence estimator 
\smash{$\hOpt$} defined in \eqref{eq:efron_div} is approximately 
unbiased for \smash{$\Opt(\htheta)$}, where 
``approximately'' here means its expectation is
correct up to first-order in a Taylor expansion.  If we could
appropriately control the error in this approximation, under 
say an exponential family distribution for $Y$, then we might be able
to extend Theorem \ref{thm:oracle_bd} the results of Section
\ref{sec:subset_reg} on subset regression to generalized linear
models.   We leave this to future work.  

\appendix
\section{Proofs}

\subsection{Proof of Lemma \ref{lem:sure_min_shrink}} 

The function $g$ as defined is not convex, but it is smooth, so the
result follows from simply checking the image of its critical points,
and the boundary points of the contraint region.  As for the latter,
we note that $g(0) = 2b$ and $g(\infty) = a$.  As for the former, we
compute 
\begin{equation*}
g'(x) = \frac{2ax}{(1+x)^2} - \frac{2ax^2}{(1+x)^3} - \frac{2b}{(1+x)^2}.
\end{equation*}
Setting this equal to 0, and solving, yields the single critical
point
\begin{equation*}
x^* = \frac{b}{a-b}.
\end{equation*}
The image of this point is $g(x^*)=2b-b^2/a$, which is always
strictly less than $g(0)=2b$ as well as $g(\infty)=a$.
Hence $x^*$ is the constrained minimizer whenever $x^* \geq 0$,
i.e., whenever $a \geq b$.  If $a < b$, then either 0 or $\infty$ is
the minimizer, and as $a<b$ by assumption, the minimizer must be
$\infty$.

\subsection{Proof of Lemma \ref{lem:chi_sq_max}}

For each $s \in S$, the moment generating function for $W_s$ is  
\begin{equation*}
\E(e^{tW_s}) = (1-2t)^{-p_s/2} \quad 
\text{for} \; 0 \leq t < 1/2. 
\end{equation*}
Now using Jensen's inequality,
\begin{align*}
\exp\Big\{ t \E \Big[\max_{s \in S} \, (W_s - p_s) \Big] \Big\}  
&\leq \E \Big[ \exp \{ t \max_{s \in S} \, (W_s - p_s) \Big\} \Big] \\ 
&= \E \Big[ \max_{s \in S} \, \exp(t(W_s - p_s)) \Big] \\
&\leq \sum_{s \in S} \E(e^{tW_s}) e^{-t p_s} \\
&= \sum_{s \in S} ((1-2t)e^{-2t})^{-p_s/2}.
\end{align*}
Taking logs of both sides and dividing by $t$, then changing variables
to $\delta = 1-2t$, gives the result.

\subsection{Proof of Theorem \ref{thm:ex_df_subset_reg}}

Simply define $\delta_n=1-a_n$, $n=1,2,3,\ldots$.  By the first
assumption in 
\eqref{eq:ex_df_subset_reg_ass}, 
\begin{equation*}
\frac{1}{1-\delta_n} \frac{\log|S|}{\Risk(\htheta_{s_0})} \to 0.
\end{equation*}
Using a Taylor expansion of the function
$f(x)=\log(1/x)$ around $x=1$, for $n$ large enough, 
\begin{equation*}
0 \leq \frac{p_{\max}}{\Risk(\htheta_{s_0})}
\bigg(\frac{\log(1/\delta_n)}{1-\delta_n} - 1\bigg) \leq 
\frac{p_{\max}}{\Risk(\htheta_{s_0})}
\bigg(\frac{1-\delta_n}{\delta_n} -
\frac{1-\delta_n}{2\delta_n^2}\bigg) \to 0,
\end{equation*}
where the limit is implied by the second assumption in
\eqref{eq:ex_df_subset_reg_ass}.  This proves the result.

\subsection{Proof of Theorem \ref{thm:ex_df_soft_thresh_tibs}} 

The SURE criterion in \eqref{eq:sure_soft_thresh} is piecewise
quadratic in $s$, and is monotone for $s$ in between adjacent
(absolute) data values $|Y_i|$, $i=1,\ldots,n$, and so it must be
minimized at one of 
these data values or at 0 (this is a common observation, e.g., made in
\citet{donoho1995adapting}). Let us denote the order  
statistics of absolute values by $|Y|_{(1)} \geq \ldots \geq
|Y|_{(n)} \geq |Y|_{(n+1)}$, where we set $|Y|_{(n+1)}=0$ for
notational convenience.
We can reparametrize the family \eqref{eq:est_soft_thresh} of
soft-thresholding estimators so that our tuning parameter now becomes
an index $k=1,\ldots,n+1$, where a choice $k$ for the index 
corresponds to a choice $s=|Y|_{(k)}$ for the threshold level.
Accordingly, we can write SURE as
\begin{align}
\label{eq:sure_soft_thresh_k_1}
\hErr_k(Y) &= k |Y|_{(k)}^2 + \sum_{j=k+1}^n |Y|_{(j)}^2 +   
2\sigma^2 (k-1),
\end{align}
and we seek to minimize this over $k=1,\ldots,n+1$.

Letting $Y_i$ vary, and keeping all other coordinates $Y_{-i}$
fixed, we will track discontinuities in the $i$th component of the
SURE-tuned soft-thresholding estimator
\begin{equation*}
\htheta_{\hat{k}(\cdot,Y_{-i}),i}(\cdot,Y_{-i}) : \R \to \R.
\end{equation*}
Note that these discontinuities  
can only occur when the minimizer of \smash{$\hErr_k(Y)$} in
\eqref{eq:sure_soft_thresh_k_1} changes, over $k=1,\ldots,n+1$; also 
note, a relabeling of the order statistics does not induce such a
discontinuity, because all 
values \smash{$\hErr_k(Y)$}, $k=1,\ldots,n+1$ behave continuously as 
through a relabeling the order statistics due to ties. Therefore,
without a loss of generality, we may assume that  
$Y_1 \geq \ldots \geq Y_n \geq Y_{n+1} = 0$ (as we can 
repeat the same arguments inside each polyhedron over which the 
labeling of order statistics remains constant).   Let us rewrite 
\eqref{eq:sure_soft_thresh_k_1} as  
\begin{equation}
\label{eq:sure_soft_thresh_k_2}
\hErr_k(Y) = k Y_k^2  + \sum_{j=k+1}^n Y_j^2 
+ 2 \sigma^2(k-1).
\end{equation}
The minimizer \smash{$\hat{k}(Y)$} of \smash{$\hErr_k(Y)$},
$k=1,\ldots,n+1$ can only jump at an equality between two of these
SURE criterion values. As $Y_i$ varies, and $Y_{-i}$ remains fixed, 
such equalities
can only happen at a finite number of points.  This, along with the 
absolute continuity of the soft-thresholding operator at a fixed
threshold level, establishes p-almost differentiability of
\smash{$\htheta_{\hat{k}}$}.  

Now as $Y_i$ varies, let us analyze the rate of change of
\smash{$\hErr_k(Y)$}, $k=1,\ldots,n+1$, in three cases:
\begin{enumerate}
\item if $k>i$, then \smash{$\hErr_k(Y)$} does not change;
\item if $k=i$, then \smash{$\hErr_k(Y)$} changes at a linear rate
  with slope $2k$;  
\item if $k<i$, then  \smash{$\hErr_k(Y)$} changes at a linear rate
  with slope $2$.  
\end{enumerate}
We can see that as $Y_i$ increases, the minimizer
\smash{$\hat{k}(Y)$} can only jump from a value $\leq i$ to a value 
$>i$.  This means that the SURE-optimal threshold level can only
decrease as $Y_i$ increases (recalling the assumed ordering $Y_1 \geq
\ldots \geq Y_n \geq 0$), which proves \eqref{eq:jump_nonneg}.  
Finally, under normality, \smash{$\edf(\htheta_{\hs}) \geq 0$} and 
\eqref{eq:df_soft_thresh_lb} follow from \eqref{eq:tibs_formula} with
\smash{$\htheta=\htheta_{\hs}$}, and the observation that the
SURE-optimal threshold value \smash{$\hs(Y)$} is constant in $Y$ at
all nondiscontinuity points, so \smash{$\partial \hs(Y)/\partial
  Y_i=0$}, $i=1,\ldots,n$ almost everywhere.

% Without a loss of generality, assume that we are varying $Y_n$, and  
% that $Y_1 > \ldots > Y_{n-1} \geq 0$.  We can rewrite 
% \eqref{eq:sure_soft_thresh_k_1} as  
% \begin{multline}
% \label{eq:sure_soft_thresh_k_2}
% \hErr_k(Y) = k Y_k^2  + \sum_{i=k+1}^n Y_i^2 
% + 2 \sigma^2(k-1) + {} \\
% 1\{Y_n \leq Y_k\}(-kY_k^2+kY_n^2) + 
% 1\{Y_n \leq Y_{k+1}\}\big((k-1)Y_{k+1}^2-(k-1)Y_n^2\big).  
% \end{multline}
% This is a continuous, piecewise quadratic function of $Y_n$.   
% Let $j$ be the largest index among $1,\ldots,n-1$ such that 
% $Y_n \leq Y_j$, and 

\subsection{Proof of Theorem \ref{thm:ex_df_subset_reg_mh}}

This proof is essentially already found in 
\citet{mikkelsen2016degrees} (in their Section 5, where they study the
Lagrangian formulation of best subset selection). For completeness, we
recapitulate the arguments.

First note that, for any $s,t \in S$, we can express the difference
between SURE criterions  \eqref{eq:sure_subset_reg} for models $s$ and
$t$, each evaluated at an arbitrary point $y \in \R^n$, as
\begin{equation}
\label{eq:sure_subset_reg_diff}
\hErr_s(y)-\hErr_t(y) =
y^T (P_t-P_s) y + 2\sigma^2 (p_s-p_t).
\end{equation}
For $s \in S$, let us define $U_s$ to be the set of all
points $y \in \R^n$ such that the SURE criterion evaluated at $y$ is
strictly lower for model $s$ than for all other tuning parameter
values, i.e.,  
\begin{equation}
\label{eq:us_set}
U_s = \bigcap_{t \in S \setminus \{s\}}  
\Big\{ y \in \R^n : y^T (P_t-P_s) y + 2\sigma^2 (p_s-p_t) < 0
\Big\}.
\end{equation}
By construction \smash{$\htheta_{\hs}|_{U_s}
  = \htheta_s$}, which is a linear function and clearly Lipschitz.
It is clear that the sets $U_s$, $s \in S$ are regular open 
(a regular open set is one that is equal to the interior of its
closure) and that their closures cover $\R^n$. 
This proves that \smash{$\htheta_{\hs}$} is piecewise Lipschitz.

Now for any $s,t \in S$ and \smash{$y \in \bar{U}_s \cap
  \bar{U}_t$}, we will compute the tangent space to 
$\partial U_s$ at $y$.  This can be seen as the
collection of derivatives $\gamma'(0)$ of smooth curves $\gamma :
(-1,1) \to \partial U_s$ such that $\gamma(0)=y$.  We can compute such
derivatives by implicit differentiation.  Consider a smooth curve
$\gamma$ satisfying $\gamma(0)=y$ and $\gamma(x) \in \partial U_s 
\cap \partial U_t$ for $|x|$ sufficiently small.  Then for such $x$,
\smash{$\hErr_s(\gamma(x))=\hErr_t(\gamma(x))$}, which from
\eqref{eq:sure_subset_reg_diff}, can be written as
\begin{equation*}
\gamma(x)^T (P_t-P_s) \gamma(x) = 2\sigma^2 (p_t-p_s). 
\end{equation*}
Differentiating with respect to $x$, using the chain rule, and
evaluating this at $x=0$, gives  
\begin{equation*}
y^T (P_t-P_s) \gamma'(0) = 0,
\end{equation*}
which defines an $(n-1)$-dimensional subspace in which the derivative 
$\gamma'(0)$ must lie.  This shows us that the tangent space to
$\partial U_s$ at $y$ is $\{z \in \R^n : y^T (P_t-P_s) z = 0\}$, and
thus the outer unit normal vector to $\partial U_s$ at $y$ is
precisely as in \eqref{eq:unit_normal_subset_reg}.
(The orientation assigned to $\eta_s(y)$ in
\eqref{eq:unit_normal_subset_reg} is important: it is
oriented to point from $U_s$ to $U_t$, which can be verified by
examining the directional derivative of \smash{$\hErr_s-\hErr_t$} in
the direction of $\eta_s(y)$, evaluated at the point $y$, and checking 
that this is positive.)

Assuming normality of $Y$, the result in
\eqref{eq:ex_df_subset_reg_mh} is a direct application of 
\eqref{eq:mh_formula}.  For any $s,t \in S$ and \smash{$y \in
  \bar{U}_s \cap \bar{U}_t$}, it is immediate from
\eqref{eq:unit_normal_subset_reg} that  
\begin{equation*}
\Big\langle \htheta_t(y)-\htheta_s(y), \eta_s(y) \Big\rangle = 
\bigg\langle (P_t-P_s)y, \frac{(P_t-P_s)y}{\|(P_t-P_s)y\|_2}
\bigg\rangle =  \|(P_t-P_s)y\|_2,
\end{equation*}
which verifies \eqref{eq:ex_df_subset_reg_mh}. 

\subsection{Proof of Theorem \ref{thm:ex_df_subset_reg_nested}} 

Assume all models in $S$ are nested. For a pair $s,t \in S$ 
satisfying (say) $s \subseteq t$, i.e., $\col(X_s) \subseteq
\col(X_t)$, note that $P_t-P_s$ is itself a projection matrix (onto 
$\col(X_t) \setminus \col(X_s)$), and so for any \smash{$y \in 
  \bar{U}_s \cap \bar{U}_t$}, 
\begin{equation}
\label{eq:pt_minus_ps} 
\|(P_t-P_s)y\|_2^2 = y^T (P_t-P_s) y = 2\sigma^2 (p_t-p_s),
\end{equation}
where the first equality comes from idempotence and the second from
\eqref{eq:us_set}.  Plugging this into the result 
\eqref{eq:ex_df_subset_reg_mh} from Theorem
\ref{thm:ex_df_subset_reg_mh}, for all $s,t \in S$, verifies   
\eqref{eq:ex_df_subset_reg_mh_nested}.  

We work on bounding the integrals appearing in
\eqref{eq:ex_df_subset_reg_mh_nested}.  To rephrase
\eqref{eq:pt_minus_ps}, we know that for each $s,t \in S$ with $s
\subseteq t$, 
\begin{equation}
\label{eq:us_ut_bd_1}
\bar{U}_s \cap \bar{U}_t \subseteq \Big\{y \in \R^n :
\|(P_t-P_s)y\|_2^2 = 2\sigma^2 (p_t-p_s) \Big\}.
\end{equation} 
We could certainly integrate the normal density over the set on
the right-hand side above in order to bound its integral over
\smash{$\bar{U}_s \cap \bar{U}_t$}, but it turns out that the simple 
containment in 
\eqref{eq:us_ut_bd_1} is a bit too loose. In words, at each point
\smash{$y
  \in \bar{U}_s \cap \bar{U}_t$}, we know that the SURE criterions for 
models $s$ and $t$ must be equal, and this is precisely what 
is reflected on the right-hand side in \eqref{eq:us_ut_bd_1}; however, 
we are missing the fact that the SURE criterions for all 
other models $r$ must be no smaller than the common criterion value 
achieved by models $s,t$.  

To develop a more refined approach, we first note that each integral in  
\eqref{eq:ex_df_subset_reg_mh_nested}  
can be taken over \smash{$\bar{U}_s \cap \bar{U}_t 
  \cap \{y \in \R^n : \eta_s(y) \not= 0\}$} (rather than
\smash{$\bar{U}_s \cap \bar{U}_t$}), as in each term of
\eqref{eq:ex_df_subset_reg_mh} the integrand is zero whenever the  
outer unit normal vector is zero.  In our current setup (i.e., disjoint
regular open sets whose closures cover $\R^n$), it can be shown that 
the outer unit normal $\eta_s$ vanishes on \smash{$\bar{U}_s
  \cap \bar{U}_t \cap \bar{U}_r$}, when $s,t,r$ are distinct,
except on a set of $\cH^{n-1}$ measure zero (e.g., see Lemma A.2 of 
\citet{mikkelsen2016degrees}).  Therefore, we can exactly characterize  
\begin{multline}
\label{eq:us_ut_eq_2}
\bar{U}_s \cap \bar{U}_t \cap \{y \in \R^n : \eta_j(y)\not=0 \big\} 
= \cN \cup \Big\{y \in \R^n : 
\|(P_t-P_s)y\|_2^2 = 2\sigma^2 (p_t-p_s), \\
y^T(P_r-P_s)y + 2\sigma^2(p_s-p_r) < 0 \;\, \text{and} \;\, 
y^T(P_r-P_t)y + 2\sigma^2(p_t-p_r) < 0, \; 
\text{for all $r \not= s,t$} \Big\},
\end{multline}
where $\cN$ is a set of $\cH^{n-1}$ measure zero. 
Identifying (say) $s=\{1,\ldots,j\}$ and $t=\{1,\ldots,k\}$, we can
rewrite \eqref{eq:us_ut_eq_2} as
\begin{multline}
\label{eq:us_ut_eq_3}
\bar{U}_j \cap \bar{U}_k \cap \{y \in \R^n : \eta_j(y)\not=0 \big\}  
={} \\ \cN \cup \Big\{y \in \R^n :  
\|(P_k-P_j)y\|_2^2 = 2\sigma^2 (k-j), \;
\|(P_j-P_\ell)y\|_2^2 > 2\sigma^2(j-\ell), \; 
\text{for $\ell < j$}, \\  
\|(P_\ell-P_j)y\|_2^2 < 2\sigma^2(\ell-j), \; 
\text{for $j<\ell<k$}, \;
\|(P_\ell-P_k)y\|_2^2 < 2\sigma^2(\ell-k), \; 
\text{for $\ell>k$} \Big\}.
\end{multline}
% \begin{multline}
% \label{eq:us_ut_eq_3}
% \bar{U}_j \cap \bar{U}_k \cap \{y \in \R^n : \eta_j(y)\not=0 \big\}  
% = \cN \cup \Big\{y \in \R^n :  
% \|(P_k-P_j)y\|_2^2 = 2\sigma^2 (k-j), \\
% \|(P_j-P_\ell)y\|_2^2 > 2\sigma^2(j-\ell), \; 
% \text{for $\ell < j$}, \;  
% \|(P_\ell-P_j)y\|_2^2 < 2\sigma^2(\ell-j), \; 
% \text{for $\ell>j$, $\ell \not= k$}, \\ 
% \|(P_k-P_\ell)y\|_2^2 > 2\sigma^2(k-\ell), \; 
% \text{for $\ell<k$, $\ell \not= j$}, \;  
% \|(P_\ell-P_k)y\|_2^2 < 2\sigma^2(\ell-k), \; 
% \text{for $\ell>k$} \Big\}.
% \end{multline}
Let $v_1,\ldots,v_p \in \R^n$ be orthonormal basis vectors that span
$\col(X)$, constructed so that $v_i$ spans the column space
of $P_i-P_{i-1}$ for each $i=1,\ldots,p$ (where we take $P_0=0$ for
notational convenience), i.e., \smash{$v_i = P_{i-1}^\perp
  X_i/\|P_{i-1}^\perp X_i\|_2$}, 
$i=1,\ldots,p$ as in the theorem statement.  Then
\eqref{eq:us_ut_eq_3} becomes   
\begin{multline*} 
\bar{U}_j \cap \bar{U}_k \cap \{y \in \R^n : \eta_j(y)\not=0 \big\}  
={} \\ \cN \cup \bigg\{y \in \R^n :  
\sum_{i=j+1}^k (v_i^T y)^2 = 2\sigma^2 (k-j), \;
\sum_{i=\ell+1}^j (v_i^T y)^2 > 2\sigma^2(j-\ell), \;  
\text{for $\ell < j$}, \\  
\sum_{i=j+1}^\ell (v_i^T y)^2 < 2\sigma^2(\ell-j), \;  
\text{for $j<\ell<k$}, \; 
\sum_{i=k+1}^\ell (v_i^T y)^2 < 2\sigma^2(\ell-k), \;  
\text{for $\ell>k$} \bigg\}.
\end{multline*}
% \begin{multline*} 
% \bar{U}_j \cap \bar{U}_k \cap \{y \in \R^n : \eta_j(y)\not=0 \big\}  
% = \cN \cup \bigg\{y \in \R^n :  
% \sum_{i=j+1}^k (v_i^T y)^2 = 2\sigma^2 (k-j), \\ 
% \sum_{i=\ell+1}^j (v_i^T y)^2 > 2\sigma^2(j-\ell), \;  
% \text{for $\ell < j$}, \;  
% \sum_{i=j+1}^\ell (v_i^T y)^2 < 2\sigma^2(\ell-j), \;  
% \text{for $\ell>j$, $\ell \not= k$}, \\ 
% \sum_{i=\ell+1}^k (v_i^T y)^2 > 2\sigma^2(k-\ell), \;  
% \text{for $\ell<k$, $\ell \not= j$}, \;  
% \sum_{i=k+1}^\ell (v_i^T y)^2 < 2\sigma^2(\ell-k), \;  
% \text{for $\ell>k$} \bigg\}.
% \end{multline*}
Integrating the normal density over the set on the
right-hand side above, with respect to the appropriate
($(n-1)$-dimensional Hausdorff) measure, gives 
\begin{multline}
\label{eq:us_ut_eq_4}
\sigma\int_{\bar{U}_j \cap \bar{U}_k \cap \{\eta_j(y) \not=0 \}} 
\phi_{\theta_0,\sigma^2 I}(y) \, d\cH^{n-1}(y) ={} \\
\sigma\int_{\substack{
\sum_{i=\ell+1}^j (v_i^T y)^2 > 2\sigma^2 (j-\ell), \; \ell < j \\      
\sum_{i=j+1}^k (v_i^T y)^2 = 2\sigma^2(k-j), \;
\sum_{i=j+1}^\ell (v_i^T y)^2 < 2\sigma^2(\ell-j), \; j<\ell<k, \\ 
\sum_{i=k+1}^\ell (v_i^T y)^2 < 2\sigma^2(\ell-j), \; \ell > k}}  
\phi_{\theta_0,\sigma^2 I}(y) \, d\cH^{n-1}(y).
\end{multline}
We note that a sufficient condition for SURE at $y$ to be minimized
at one of $j,\ldots,k$, i.e., for $y$ to be an element of 
\smash{$\cup_{\ell=j}^k \bar{U}_\ell$}, is
\begin{equation*}
\sum_{i=\ell+1}^j (v_i^T y)^2 > 2\sigma^2 (j-\ell), 
\; \text{for $\ell <  j$}, \;
 \sum_{i=k+1}^\ell (v_i^T y)^2 < 2\sigma^2(\ell-k), \;  
\text{for $\ell>k$},
\end{equation*}
and so carrying on from \eqref{eq:us_ut_eq_4}, 
\begin{align*}
&\sigma\int_{\bar{U}_j \cap \bar{U}_k \cap \{\eta_j(y) \not=0 \}} 
\phi_{\theta_0,\sigma^2 I}(y) \, d\cH^{n-1}(y) \\
&\leq \sigma \P\Big(Y \in \cup_{\ell=j}^k \bar{U}_\ell \Big)  
\int_{\sum_{i=j+1}^k (v_i^T y)^2 = 2\sigma^2(k-j), \; 
\sum_{i=j+1}^\ell (v_i^T y)^2 < 2\sigma^2(\ell-j), \; j<\ell<k}  
\phi_{\theta_0,\sigma^2 I}(y) \, d\cH^{n-1}(y) \\
&= \sigma \P\Big(Y \in \cup_{\ell=j}^k \bar{U}_\ell \Big)  
\int_{\sum_{i=j+1}^k z_i^2 = 2\sigma^2(k-j), \; 
\sum_{i=j+1}^\ell z_i^2 < 2\sigma^2(\ell-j), \; j<\ell<k} 
\phi_{M^T \theta_0,\sigma^2 I}(z) \, d\cH^{n-1}(z),
\end{align*}
where $M \in \R^{n\times n}$ in defined to be an
orthogonal matrix whose first $p$ are given by
$v_1,\ldots,v_p$, i.e., given by the matrix $V \in \R^{n\times p}$ 
introduced in the theorem.  As the sets
\smash{$\bar{U}_\ell$}, $\ell=1,\ldots,d$ intersect on a set of
($n$-dimensional Lebesgue) measure zero, we can rewrite
the above as
\begin{multline}
\label{eq:us_ut_bd_5}
\sigma\int_{\bar{U}_j \cap \bar{U}_k \cap \{\eta_j(y) \not=0 \}} 
\phi_{\theta_0,\sigma^2 I}(y) \, d\cH^{n-1}(y) \leq{} \\
\sigma \sum_{\ell=j}^k \P(Y \in U_\ell) 
\int_{\sum_{i=j+1}^k z_i^2 = 2\sigma^2(k-j), \; 
\sum_{i=j+1}^\ell z_i^2 < 2\sigma^2(\ell-j), \; j<\ell<k} 
\phi_{M^T \theta_0,\sigma^2 I}(z) \, d\cH^{n-1}(z).
\end{multline}
In general, the integral in \eqref{eq:us_ut_bd_5} is difficult to
compute (though we will have luck in the case that $\theta_0=0$,
to be discussed shortly), so we can simply upper bound it 
by discarding the inequalities in the domain of integration, giving  
\begin{align*}
\sigma\int_{\bar{U}_j \cap \bar{U}_k \cap \{\eta_j(y) \not=0 \}} 
\phi_{\theta_0,\sigma^2 I}(y) \, d\cH^{n-1}(y) 
&\leq \sigma \sum_{\ell=j}^k \P(Y \in U_\ell) 
\int_{\sum_{i=j+1}^k z_i^2 = 2\sigma^2(k-j)} 
\phi_{M^T \theta_0,\sigma^2 I}(z) \, d\cH^{n-1}(z) \\
&= \sum_{\ell=j}^k \P(Y \in U_\ell) \; 
\Lambda_{k-j} \Big( B_{k-j} \big(\mu_{(j+1):k},    
\sqrt{2(k-j)} \big) \Big),
\end{align*}
where the last line used the definition of
Gaussian surface area, recalling the notation 
$\mu = V^T \theta_0/\sigma$ as in the
theorem.  Summing the above bound over all pairs $j<k$ with separation
$k-j=d$ gives
\begin{align*}
\sigma \sum_{j=1}^{p-d} 
\int_{\bar{U}_j \cap \bar{U}_{j+d} \cap \{\eta_j(y) \not=0 \}}  
\phi_{\theta_0,\sigma^2 I}(y) \, d\cH^{n-1}(y) 
&\leq \sum_{j=1}^{p-d} 
\sum_{\ell=j}^{j+d} \P(Y \in U_\ell) \; 
\Lambda_d \Big( B_d \big(\mu_{(j+1):(j+d)},    
\sqrt{2d} \big) \Big) \\
&\leq (d+1) \max_{j=1,\ldots,d} \,
\Lambda_d \Big( B_d \big(\mu_{(j+1):(j+d)}, \sqrt{2d} \big) \Big), 
\end{align*}
where in the last line, we recognized that each index $\ell$ appears 
in the double sum $d+1$ times. An upper bound on the full sum (over all
pairs $j,k$) in \eqref{eq:ex_df_subset_reg_mh_nested} is given
by multiplying the last line above by \smash{$\sqrt{2d}$}, and
summing this over $d=1,\ldots,p$, which establishes
\eqref{eq:ex_df_subset_reg_mh_balls}.  

When the balls in \eqref{eq:ex_df_subset_reg_mh_balls} are
all centered at the origin, i.e., when $\theta_0=0$ (or more
generally, this would happen in a nested family $S$ such that $P_s
\theta_0=\theta_0$ for all $s \in S$), we can upper bound 
\eqref{eq:ex_df_subset_reg_mh_balls} by invoking  
%standard chi-squared probability bounds and 
known results on the Gaussian surface area of balls.  Importantly,
though, it turns out to be more fruitful to return to an earlier step
along the way to deriving \eqref{eq:ex_df_subset_reg_mh_balls},
namely, the integral on the right-hand side in \eqref{eq:us_ut_bd_5},
which recall we upper bounded in the general $\theta_0$ case by 
dropping the inequality constraints in the domain of integration. Let
us write this integral as 
\begin{equation}
\label{eq:key_integral_1}
\P\Bigg(\sum_{i=j+1}^\ell W_i < 2 (\ell-j), \;
\text{for $j<\ell<k$} \; \Bigg| \;
\sum_{i=j+1}^k W_i = (k-j) \Bigg)
\Lambda_{k-j} \Big( B_{k-j} \big(0, \sqrt{2(k-j)} \big) \Big), 
\end{equation}
where $W_i$, $i=j+1,\ldots,k$ are i.i.d.\ $\chi^2_1$ random
variates. To simplify notation, we denote $k-j=d$ and relabel 
these random variates as $W_1,\ldots,W_d$.  Because $W_1,\ldots,W_d$
are i.i.d., they are still i.i.d.\ conditional on their sum being
equal to $2d$, and when we further condition on 
$(W_1,\ldots,W_d)$ being equal to $(w_1,\ldots,w_d)$ up to a
circular permutation, any ones of the $d$ options 
\begin{equation*}
(w_1,w_2,\ldots,w_d), \; (w_d, w_1, \ldots,w_{d-1}), \; \ldots, \; 
(w_2,w_3,\ldots,w_1)
\end{equation*}
is equally likely.  Now we recall and apply the following classic
result in combinatorics. 

\begin{proposition}[The gas stations problem]
\label{prop:gas_stations}
Let $w_1,\ldots,w_d$ be nonnegative numbers that sum to $2d$. Then
there exists exactly one circular permutation of $(w_1,\ldots,w_d)$,
call it \smash{$(w_{i_1},\ldots,w_{i_d})$}, such that 
\begin{equation*}
w_{i_1}+\ldots+w_{i_q} \leq 2q, \; \text{for all $q=1,\ldots,d$}. 
\end{equation*}
\end{proposition}

By Proposition \ref{prop:gas_stations} and the discussion preceding
it, we see that \eqref{eq:key_integral_1} becomes simply
\begin{equation}
\label{eq:key_integral_2}
\frac{1}{d}
\Lambda_d \big( B_d (0, \sqrt{2d}) \big),
\end{equation}
and by following the exact same steps leads up to
\eqref{eq:ex_df_subset_reg_mh_balls}, we obtain the sharper upper
bound that is given by the first inequality of 
\eqref{eq:ex_df_subset_reg_mh_balls_null}. 

For the Gaussian surface area of an origin-centered ball, 
\citet{ball1993reverse} gave the formula  
\begin{equation*}
\Lambda_d\big(B_d(0,r)\big) = \frac{r^{d-1} e^{-r^2/2}}
{2^{d/2-1} \Gamma(d/2)},
\end{equation*}
in any dimension $d$ (see \citet{klivans2008learning} for a simple, 
direct proof). Plugging this formula into the first inequality in 
\eqref{eq:ex_df_subset_reg_mh_balls_null} gives
\begin{equation*}
\sum_{d=1}^p \sqrt{2d}\bigg(1+\frac{1}{d}\bigg) 
\Lambda_d \big( B_d (0, \sqrt{2d}) \big) \leq
2 \sum_{d=1}^p
\bigg(1 + \frac{1}{d}\bigg) 
\frac{d^{d/2} e^{-d}}
{\Gamma(d/2)}.
\end{equation*}
Continuing on with the chain of upper bounds, we apply the  
following Stirling-type bound for the gamma function (e.g.,
\citet{jameson2015simple}),   
\begin{equation*}
\frac{x^{x-1/2} e^{-x}}{\Gamma(x)} \leq \frac{1}{\sqrt{2\pi}} 
\quad \text{for all $x > 0$},
\end{equation*}
which yields
\begin{equation}
\label{eq:d_sum_0}
2 \sum_{d=1}^p
\bigg(1 + \frac{1}{d}\bigg) 
\frac{d^{d/2} e^{-d}}
{\Gamma(d/2)} \leq
%\frac{2}{\sqrt{2\pi}} \sum_{d=1}^p 
%\bigg(\sqrt{d}  + \frac{1}{\sqrt{d}}\bigg) 
%2^{d/2-1/2} e^{-d/2} = 
\frac{1}{\sqrt\pi} \sum_{d=1}^p 
\bigg(\sqrt{d}  + \frac{1}{\sqrt{d}}\bigg) 
\bigg(\frac{2}{e}\bigg)^{d/2}.
\end{equation}
We split the right-hand side above into two sums and bound each
individually. Consider first
\begin{equation}
\label{eq:d_sum_1}
\frac{1}{\sqrt\pi} \sum_{d=1}^p 
\sqrt{d} \bigg(\frac{2}{e}\bigg)^{d/2} \leq 
\frac{1}{\sqrt\pi} \sum_{d=1}^\infty 
\sqrt{d} \bigg(\frac{2}{e}\bigg)^{d/2} \leq 
\frac{1}{\sqrt\pi} \sum_{d=1}^N 
\sqrt{d} \bigg(\frac{2}{e}\bigg)^{d/2} +
\frac{1}{\sqrt\pi} \sum_{d=N+1}^\infty
d \bigg(\frac{2}{e}\bigg)^{d/2}.
\end{equation}
The second term on the right-hand side above can be calculated as  
\begin{align}
\nonumber
\frac{1}{\sqrt\pi} \sum_{d=N+1}^\infty
d \bigg(\frac{2}{e}\bigg)^{d/2} &= 
\sqrt{\frac{2}{\pi e}} \sum_{d=N+1}^\infty d 
\bigg(\sqrt{\frac{2}{e}}\bigg)^{d-1} \\
\nonumber
&= \sqrt{\frac{2}{\pi e}} \frac{d}{dx}
\bigg(\sum_{d=N+1}^\infty x^d\bigg)\bigg|_{x=\sqrt{2/e}} \\ 
\label{eq:d_sum_2}
&= \frac{1}{\sqrt{\pi}} \frac{\sqrt{2/e}^{N+1}}{1-\sqrt{2/e}}  
\bigg(N+1-\frac{\sqrt{2/e}}{1-\sqrt{2/e}}\bigg).
\end{align}
Thus we can upper bound the right-hand side in \eqref{eq:d_sum_1} by
computing the first sum with $N=1000$ numerically and computing the
second via \eqref{eq:d_sum_2}, which gives 
\begin{equation}
\label{eq:d_sum_3}
\frac{1}{\sqrt\pi} \sum_{d=1}^{1000} 
\sqrt{d} \bigg(\frac{2}{e}\bigg)^{d/2} +
\frac{1}{\sqrt\pi} \sum_{d=1001}^\infty
d \bigg(\frac{2}{e}\bigg)^{d/2} < 8.21.
\end{equation}
It remains to consider
\begin{equation}
\label{eq:d_sum_4}
\frac{1}{\sqrt\pi} \sum_{d=1}^p 
\frac{1}{\sqrt{d}} \bigg(\frac{2}{e}\bigg)^{d/2} \leq  
\frac{1}{\sqrt\pi} \sum_{d=1}^\infty 
\frac{1}{\sqrt{d}} \bigg(\frac{2}{e}\bigg)^{d/2} \leq  
\frac{1}{\sqrt\pi} \sum_{d=1}^N 
\frac{1}{\sqrt{d}} \bigg(\frac{2}{e}\bigg)^{d/2} + 
\frac{1}{\sqrt\pi} \sum_{d=N+1}^\infty
\bigg(\frac{2}{e}\bigg)^{d/2}.
\end{equation}
As before, the second term in \eqref{eq:d_sum_4} we can compute as  
\smash{$(1/\sqrt{\pi})\sqrt{2/e}^{N+1}(1-\sqrt{2/e})^{-1}$}, and the
first term we can evaluate numerically at $N=1000$, which gives  
\begin{equation}
\label{eq:d_sum_5}
\frac{1}{\sqrt\pi} \sum_{d=1}^{1000} 
\frac{1}{\sqrt{d}} \bigg(\frac{2}{e}\bigg)^{d/2} +
\frac{1}{\sqrt\pi} \sum_{d=1001}^\infty
\bigg(\frac{2}{e}\bigg)^{d/2} < 1.75.
\end{equation}
Putting \eqref{eq:d_sum_3} and \eqref{eq:d_sum_5} together,
we can upper bound the right-hand side in \eqref{eq:d_sum_0} by
$8.21+1.75=9.96<10$, which establishes the second inequality in
\eqref{eq:ex_df_subset_reg_mh_balls_null}, and completes the proof. 

\subsection{Proof of Theorem \ref{thm:df_best_subset}}

For the lower bound, we note that an argument analogous to that given
in the proof of Theorem \ref{thm:ex_df_subset_reg_mh} shows that the
excess degrees of freedom of subset selection, i.e., the quantity
\begin{equation*}
\df(X\hbeta^{\mathrm{subset}}_\lambda) -
  \E\|X\hbeta^{\mathrm{subset}}_\lambda(Y)\|_0,
\end{equation*}
is exactly equal to the right-hand side in
\eqref{eq:ex_df_subset_reg_mh}, where the sum is taken over all pairs
of subsets.  See Section 5 of
\citet{mikkelsen2016degrees}. Nonnegativity of the integrand in each
term of the sum therefore proves the lower bound in
\eqref{eq:df_best_subset}.  

Meanwhile, the search degrees of freedom is
upper bounded by the quantity considered in \eqref{eq:chi_sq_max} of 
Lemma \ref{lem:chi_sq_max}, where $S$ is the set of all subsets of
$\{1,\ldots,p\}$.  The upper bound is thus
\begin{align*}
\min_{\delta \in [0,1)} \; \frac{2}{1-\delta} \log \sum_{s \in S}
(\delta e^{1-\delta})^{-p_s/2} &= 
\min_{\delta \in [0,1)} \; \frac{2}{1-\delta} \log \sum_{k=0}^p
{p \choose k} \Big((\delta e^{1-\delta})^{-1/2}\Big)^k \\
&= \min_{\delta \in [0,1)} \; \frac{2p}{1-\delta}
\log\Big(1+(\delta e^{1-\delta})^{-1/2}\Big),
\end{align*}
where the last step used the binomial theorem.  Sraightforward
numerical calculation shows that
\begin{equation*}
\min_{\delta \in [0,1)} \; \frac{2}{1-\delta}
\log\Big(1+(\delta e^{1-\delta})^{-1/2}\Big) < 1.145,
\end{equation*}
completing the proof.

\subsection{Derivation details for \eqref{eq:ex_opt_shrink_hetero}}

First, we compute
\begin{equation*}
\frac{\partial \hTheta_i}{\partial s} (Y,s)) = 
-\frac{Y_i \sigma_i^2}{(1+\sigma_i^2 s)^2}. 
\end{equation*}
Next, 
\begin{equation*}
\frac{\partial G}{\partial s} (Y,s) =
\sum_{i=1}^n \bigg(
\frac{2 Y_i^2 \sigma_i^2 s}
{(1+\sigma_i^2 s)^2} -
\frac{2 Y_i^2 \sigma_i^4 s^2}
{(1+\sigma_i^2 s)^3} -
\frac{2 \sigma_i^2}
{(1+\sigma_i^2 s)^2}
\bigg). 
\end{equation*}
Then,
\begin{equation*}
\frac{\partial^2 G}{\partial Y_i \partial s} (Y,s) = 
\frac{4Y_i \sigma_i^2 s}{(1+\sigma_i^2 s)^2}
\bigg(1 - \frac{\sigma_i^2 s}{1+\sigma_i^2 s} \bigg) 
= \frac{4Y_i \sigma_i^2 s}{(1+\sigma_i^2 s)^3}.
\end{equation*}
Finally, 
\begin{align*}
\frac{\partial^2 G}{\partial s^2} (Y,s) &=
\sum_{i=1}^n \bigg(
\frac{2 Y_i^2 \sigma_i^2}
{(1+\sigma_i^2 s)^2} -
\frac{4 Y_i^2 \sigma_i^4 s}
{(1+\sigma_i^2 s)^3} -
\frac{4 Y_i^2 \sigma_i^4 s}
{(1+\sigma_i^2 s)^3} +
\frac{6 Y_i^2 \sigma_i^6 s^2}
{(1+\sigma_i^2 s)^4} +
\frac{4 \sigma_i^4}
{(1+\sigma_i^2 s)^3}
\bigg) \\
&= \sum_{i=1}^n \bigg[
\frac{2\sigma_i^2}{(1+\sigma_i^2 s)^2} 
\bigg( Y_i^2 - 
\frac{4 Y_i^2 \sigma_i^2 s}
{1+\sigma_i^2 s} +
\frac{3 Y_i^2 \sigma_i^4 s^2}
{(1+\sigma_i^2 s)^2} +
\frac{2 \sigma_i^2}
{1+\sigma_i^2 s} 
\bigg) \bigg].
\end{align*}
Therefore, plugging the relevant quantities into
\eqref{eq:ex_opt_implicit_hetero}, we get
\eqref{eq:ex_opt_shrink_hetero}. 

\bibliographystyle{plainnat}
\bibliography{ryantibs}

\end{document}